\documentclass[11pt]{amsart}
\usepackage[T1]{fontenc}
\usepackage[utf8]{inputenc}
\usepackage{lmodern,microtype}
\usepackage[a4paper,margin=1.15in]{geometry}
\usepackage{amsmath,amssymb,amsthm}
\usepackage{booktabs,enumitem,array}
\usepackage{xcolor}
\usepackage[colorlinks=true,linkcolor=blue!45!black,citecolor=blue!45!black,urlcolor=blue!45!black]{hyperref}
\makeatletter
\renewcommand\paragraph{\@startsection{paragraph}{4}{\z@}{.4\baselineskip}{-\fontdimen2\font}{\normalfont\itshape}}
\makeatother
\newtheorem{theorem}{Theorem}[section]
\newtheorem{lemma}[theorem]{Lemma}
\newtheorem{proposition}[theorem]{Proposition}
\newtheorem{corollary}[theorem]{Corollary}
\newtheorem{compres}[theorem]{Computational result}
\theoremstyle{definition}
\newtheorem{definition}[theorem]{Definition}
\theoremstyle{remark}
\newtheorem{remark}[theorem]{Remark}
\newcommand{\cA}{\mathcal A}
\newcommand{\cB}{\mathcal B}
\newcommand{\cC}{\mathcal C}
\newcommand{\cO}{\mathcal O}
\newcommand{\cI}{\mathcal I}
\newcommand{\cR}{\mathcal R}
\newcommand{\PP}{\mathbb P}
\newcommand{\CC}{\mathbb C}
\newcommand{\QQ}{\mathbb Q}
\newcommand{\ZZ}{\mathbb Z}
\newcommand{\T}[1]{T\langle #1\rangle}
\DeclareMathOperator{\mdr}{mdr}
\DeclareMathOperator{\PGL}{PGL}
\DeclareMathOperator{\Aut}{Aut}
\newcommand{\length}{\operatorname{length}}
\hypersetup{pdftitle={Generic splitting types and Terao's conjecture for seventeen and eighteen lines},pdfauthor={Cristian Anghel and Rares Stan}}

\title[Generic splitting types and Terao's conjecture for 17 and 18 lines]{Generic splitting types and
Terao's conjecture\\ for seventeen and eighteen lines}

\author{Cristian Anghel}
\address{Simion Stoilow Institute of Mathematics of the Romanian Academy, P.O.\ Box 1-764, RO-014700 Bucharest, Romania}
\email{cristian.anghel@imar.ro}

\author{Rare\c{s} Stan}
\address{Simion Stoilow Institute of Mathematics of the Romanian Academy, P.O.\ Box 1-764, RO-014700 Bucharest, Romania}
\email{rstan@imar.ro}

\date{November 2026}
\subjclass[2020]{Primary 14N20; Secondary 32S22, 14F06, 37F75, 37F10, 52C35, 05B35}
\keywords{Line arrangements, free arrangements, Terao's conjecture, generic splitting type, logarithmic
vector fields, foliations, inflection divisor, critically fixed rational maps, rank $3$ matroids,
realization spaces}

\begin{document}

\begin{abstract}
We prove that the generic splitting type of the logarithmic bundle of an arrangement of at most eighteen
lines in the complex projective plane is determined by the intersection lattice.  In particular, Terao's
freeness conjecture holds for seventeen and eighteen lines.  For any number of lines, it holds
for free arrangements whose smaller exponent is at most eight.  By the theorem for at most
sixteen lines, the problem reduces to a uniformity condition for lattices of seventeen lines: if all
realizations have minimal Jacobian-syzygy degree at least seven, then all or none have degree seven.  The
proof has two parts.  First, a realization of minimal degree seven has defect at most four, the defect
being the length of the zero scheme of a minimal section.  This uses the residual inflection quartic of the
associated degree-seven foliation, local inflection identities, and fixed-point theory for critically fixed
rational maps of degree seven with exact finite certificates.  Second, the uniformity condition holds at
these defects.  For multiplicities at most six, a critical realization contains a twelve-line core of
minimal weight, and an exact search through extensions of the cores shows that the lattice has a single
Galois orbit of realizations; the only family left at defect four has no critical member.  For a point of
multiplicity seven we use the zero scheme of the pencil derivation and a classification of transversal
structures.  Terao's conjecture for seventeen and eighteen lines needs only defects one, two and three
and, for seventeen lines, this classification at multiplicity seven.  The proof is computer-assisted.
\end{abstract}

\maketitle
\tableofcontents

\section{Introduction and main results}\label{sec:intro}
Let $\cA$ be an arrangement of $d$ lines in $\PP^2=\PP^2_\CC$ and $\T\cA$ its logarithmic bundle,
normalized by $c_1(\T\cA)=1-d$.  On a general line $\ell$,
\[
 \T\cA|_\ell\cong\cO_\ell(-k)\oplus\cO_\ell(k+1-d),\qquad k=k(\cA)=\min\{\mdr(\cA),\lfloor(d-1)/2\rfloor\},
\]
where $\mdr(\cA)$ is the least degree of a nonzero Jacobian syzygy (see \cite{AD18,GST16}; when
$2\mdr(\cA)\ge d-1$ the bundle is semistable and this is the Grauert--M\"ulich theorem \cite{OSS}).  The
Chern classes of $\T\cA$ are determined by the intersection lattice, but $\mdr$ is not: for Ziegler's pair
of arrangements of $9$ lines with $6$ triple points \cite{Zie89a} it is $5$ or $6$ according to whether the
triple points lie on a conic.  One may ask whether the generic splitting type is determined by the
lattice; this question was asked by Cook, Harbourne, Migliore and Nagel \cite[Question~7.12]{CHMN18} and
stated as a conjecture by Dimca \cite[Conjecture~3.9]{Dim24} (see also \cite{DP26}).  For $d\le16$ it was
proved in \cite{GST16}.

\begin{theorem}[Generic splitting type through eighteen lines]\label{thm:main-gst}
If two arrangements of at most eighteen lines in $\PP^2_\CC$ have isomorphic intersection lattices, their
logarithmic bundles have the same generic splitting type.
\end{theorem}

Terao's conjecture \cite[Conjecture~4.138]{OT92} asks whether freeness \cite{Sai80} of a complex
hyperplane arrangement is determined by its intersection lattice; the exponents of a free arrangement are
determined by its lattice \cite{Ter81}.  For arrangements of lines it was known for at most twelve lines \cite{FV14}, for thirteen
\cite{DIM19} and fourteen lines \cite{BBJKL21} (see also \cite{BK23}), and for fifteen and sixteen lines
\cite{AS}.  For seventeen lines it was proved recently by Monteiro, S\'a Earp and Silva \cite{MSS},
independently of the present work, with a reduction to thirteen-line cores instead of the twelve-line cores
used here; the sixteen lattices of \cite[Table~3]{MSS} are those of Section~\ref{sec:class41}
(Remark~\ref{rem:mss}).  It is also
known for free arrangements whose smaller exponent is at most $5$ \cite[Corollary~5.5]{Abe14}, at most $7$
\cite[Theorem~C]{GST16}, and at most $8$ by Corollary~\ref{cor:anyd} below, and for free
arrangements with a point of multiplicity at least the smaller exponent \cite{DimcaNF} (see also
\cite{DKP26}), and divisionally free arrangements, a combinatorially defined class, are free
\cite{Abe16}.

\begin{theorem}[Terao for seventeen and eighteen lines]\label{thm:terao1718}
Let $\cA$ and $\cB$ be arrangements of $d\in\{17,18\}$ lines in $\PP^2_{\CC}$ with isomorphic intersection
lattices.  If $\cA$ is free, then $\cB$ is free.
\end{theorem}

Theorem~\ref{thm:terao1718} follows from Theorem~\ref{thm:main-gst}: in a counterexample the free member has
exponents $(8,8)$, respectively $(8,9)$, hence $k=8$, while the non-free member has $\mdr=7$
(Propositions~\ref{red:17} and~\ref{red:18}).  It also has a direct proof, given in
Sections~\ref{sec:reductions}--\ref{sec:18}, which needs much less: only Theorem~\ref{thm:main} below, its
analogue for defects two and three, the results of \cite{GST16} and, for seventeen lines, the
classification of Section~\ref{sec:sept}.  The statements are over $\CC$; no assertion in positive characteristic is made.

\subsection{The finite statements}
Let $\mathcal W$ be the class of lattices on seventeen lines satisfying
\begin{equation}\label{eq:window}
 b(\Lambda)=80,\qquad r_H\le8\text{ for every line},\qquad
 m_P\le6\text{ for every vertex}.
\end{equation}
A realization is called \emph{critical} if its minimal syzygy degree is seven.
\begin{theorem}[The defect-one window]\label{thm:main}
If $\Lambda\in\mathcal W$ has a critical realization, then $\Lambda$ has
one Galois orbit of labelled realizations. Consequently every realization
is critical and none is free.
\end{theorem}
Theorem~\ref{thm:nu} proves the same assertion with $b=81$ or $82$.  For seventeen lines, defect one
suffices, together with the classification of Section~\ref{sec:sept} for a point of multiplicity seven.
For eighteen lines, Proposition~\ref{prop:18} selects a seventeen-line deletion of defect two or three whose
multiplicities are at most six.  This is why the defects four and higher are not needed for
Theorem~\ref{thm:terao1718}, and the septuple strata are not needed for its eighteen-line part.

The eighteen-line case is not a formal consequence of the seventeen-line case.  In an eighteen-line
counterexample every line has at most eight vertices (Proposition~\ref{red:18}), so every deletion $\cA_H$
of the free member has $c_2(\T{\cA_H})=73-r_H\ge65$, whereas a free arrangement of seventeen lines has
$c_2=d_1d_2\le64$ \cite{Ter81}.  Hence no deletion of the free member is free, and Terao's conjecture for seventeen lines
says nothing about these deletions.  We use instead the stronger seventeen-line statement of
Theorem~\ref{thm:nu}: a lattice in $\mathcal W_2$ or $\mathcal W_3$ does not have both a critical and a
non-critical realization.  By Proposition~\ref{red:18}, every deletion of an eighteen-line counterexample is a
$7$-candidate with a critical and a non-critical realization, and Proposition~\ref{prop:18} chooses one in
$\mathcal W_2\cup\mathcal W_3$.  (The reduction of eighteen lines to defects two and three is
\cite[Proposition~8.2(2)]{GST16}; Proposition~\ref{prop:18} adds the bound $m_P\le6$.)  Theorem~\ref{thm:nu} is proved with the same census and the same search as
Theorem~\ref{thm:main}; eighteen lines cost two more defects of the same computation.

For Theorem~\ref{thm:main-gst}, by \cite[Theorems~A, B and~3.5]{GST16} it suffices to prove a uniformity
condition $(C_7)$ for the lattices of seventeen lines with at most eight vertices on each line, the
\emph{$7$-candidates} (Section~\ref{sec:c7}).  We prove:

\begin{theorem}[Defect bound]\label{thm:defect}
Every critical realization of a $7$-candidate has defect $\nu\le4$.
\end{theorem}

\begin{theorem}\label{thm:c7-intro}
Condition $(C_7)$ holds for every $7$-candidate of defect at most $4$, and for every $7$-candidate with a
point of multiplicity at least $7$.
\end{theorem}

Since the reduction of \cite{GST16} works for any number of lines, $(C_7)$ also has consequences for
arrangements of more than eighteen lines.

\begin{corollary}[Any number of lines]\label{cor:anyd}
Let $\cA$ and $\cB$ be arrangements of $d$ lines with isomorphic intersection lattices, where $d$ is
arbitrary.
\begin{enumerate}[label=\textup{(\alph*)}]
\item If $\min\{k(\cA),k(\cB)\}\le7$, then $k(\cA)=k(\cB)$.  In particular, if $\mdr(\cA)\le7$ and
$2\mdr(\cA)\le d-3$, then $\mdr(\cB)=\mdr(\cA)$.
\item If $\cA$ is free with exponents $(d_1,d_2)$, $d_1\le d_2$, and $d_1\le8$, then $\cB$ is free.
\end{enumerate}
\end{corollary}

Part (b) extends the case $d_1\le5$ of Terao's conjecture \cite[Corollary~5.5]{Abe14} and the case $d_1\le7$
of \cite[Theorem~C]{GST16} to $d_1\le8$; part (a) proves, for $\mdr\le7$, the conjecture
\cite[Conjecture~1.2]{DP26} that $\mdr$ is determined by the lattice when $\mdr<(d-2)/2$.

\paragraph{Logical structure.}
Table~\ref{tab:structure} shows where the main statements are proved.  For the generic splitting type,
\cite[Theorem~3.5]{GST16} shows that if $(C_s)$ holds for all $s\le S$, then the generic splitting type is
combinatorial for arrangements of at most $2S+4$ lines, and Theorems~A and~B of \cite{GST16} give $(C_s)$
for $s\le6$ (Section~\ref{sec:c7}).  Hence the single condition $(C_7)$ for all $7$-candidates, a statement
about lattices on seventeen lines, gives the generic splitting type for $17=2\cdot7+3$ and $18=2\cdot7+4$
lines at once; it is proved case by case, as in Table~\ref{tab:c7}.  For Terao's conjecture the two cases
are proved separately, and the eighteen-line proof does not use the seventeen-line one.  Neither uses the
appendices or Section~\ref{sec:nu4}; the seventeen-line proof uses Theorem~\ref{thm:sept}
(Section~\ref{sec:sept}) to exclude points of multiplicity seven, the eighteen-line proof does not.  Both use the theorem for at most sixteen lines
of \cite{GST16} (Theorems~A and~B there), through the consequences \ref{red:G1} and~\ref{red:G2} of
Section~\ref{sec:reductions}.

\begin{table}[ht]
\centering\small
\caption{Where the main statements are proved.}\label{tab:structure}
\begin{tabular}{@{}>{\raggedright\arraybackslash}p{2.75cm}>{\raggedright\arraybackslash}p{5.4cm}>{\raggedright\arraybackslash}p{3.05cm}>{\raggedright\arraybackslash}p{2.25cm}@{}}
\toprule
Statement & Reduction & Finite statement & Conclusion\\\midrule
Terao, $17$ lines\newline(Theorem~\ref{thm:terao1718}) & a counterexample gives a lattice of $\mathcal W$
with a critical and a non-critical realization (Proposition~\ref{red:17}, with Theorem~\ref{thm:sept}) & Theorem~\ref{thm:main}
(Sections~\ref{sec:crit}--\ref{sec:search}) & Corollary~\ref{cor:terao17}\\[2pt]
Terao, $18$ lines\newline(Theorem~\ref{thm:terao1718}) & a counterexample has a deletion in
$\mathcal W_2\cup\mathcal W_3$ with a critical and a non-critical realization
(Propositions~\ref{red:18} and~\ref{prop:18}) & Theorem~\ref{thm:nu}\newline(Section~\ref{sec:nu}) &
Corollary~\ref{cor:terao18}\\[2pt]
Generic splitting type, $17$ and $18$ lines (Theorem~\ref{thm:main-gst}) & $(C_7)$ for all
$7$-candidates suffices \cite[Theorems~A, B and~3.5]{GST16} & Theorem~\ref{thm:c7}(c)\newline
(Table~\ref{tab:c7}) & Section~\ref{sec:proof}\\
\bottomrule
\end{tabular}
\end{table}

\begin{table}[ht]
\centering\small
\caption{The cases of $(C_7)$ for a $7$-candidate (Theorem~\ref{thm:c7}).}\label{tab:c7}
\begin{tabular}{@{}>{\raggedright\arraybackslash}p{4.3cm}>{\raggedright\arraybackslash}p{6.6cm}>{\raggedright\arraybackslash}p{2.6cm}@{}}
\toprule
$7$-candidate & Reason & Where\\\midrule
defect $\nu\le0$ & automatic, by Bogomolov's inequality & Section~\ref{sec:sept}\\
a point of multiplicity $8$ & then $\nu=0$ & Lemma~\ref{lem:septcount}\\
a point of multiplicity $7$ & if critical, one of $21$ lattices, all of whose realizations are critical &
Theorem~\ref{thm:sept}\\
multiplicities $\le6$, $\nu=1$ & if it has a critical realization, it has one Galois orbit of
realizations, so all of them are critical (Theorem~\ref{thm:main}); proof: a twelve-line core of minimal
weight and an exact search through the realization spaces of the extensions &
Sections~\ref{sec:crit}--\ref{sec:search}\\
multiplicities $\le6$, $\nu\in\{2,3\}$ & the same (Theorem~\ref{thm:nu}), with the same census and the
same search &
Section~\ref{sec:nu}\\
multiplicities $\le6$, $\nu=4$ & the same (Theorem~\ref{thm:nu4}); the search leaves one
lattice, $\Lambda_F$, which has no critical realization (Lemma~\ref{lem:LF}) & Section~\ref{sec:nu4}\\
defect $\nu\ge5$ & no critical realization (Theorem~\ref{thm:defect}) & Appendices~\ref{app:A}
and~\ref{app:B}\\
\bottomrule
\end{tabular}
\end{table}

\paragraph{Methods.}
For multiplicities at most six and defects one to four, a critical realization contains a twelve-line
core of minimal weight.  Minimality of the core constrains the five remaining lines, and an exact search
through realization spaces stops at empty or Galois-rigid subarrangements
(Sections~\ref{sec:crit}--\ref{sec:reliability}).  At defect four the search leaves a single seventeen-line
lattice $\Lambda_F$ whose realizations form a one-parameter family; none of them is critical, by
determinants of the degree-seven Ziegler multirestrictions (Section~\ref{sec:nu4}).  For a point of
multiplicity seven, the pencil derivation of the point has a zero scheme off the pencil, and criticality
becomes a collinearity condition on that scheme; this leads to a classification of transversal structures
and to exactly $21$ critical lattices (Section~\ref{sec:sept}).  A failed or incomplete \textsc{Singular}
calculation of the general realization engine is never read as an empty realization space or as a single
Galois orbit: the lattice is then treated as open, and the search continues below it.  Section~\ref{sec:reliability} states what has been verified, and how.

Theorem~\ref{thm:defect} is proved in Appendices~\ref{app:A} and~\ref{app:B}.  A minimal derivation of
degree seven defines a saturated foliation of degree seven for which all lines of $\cA$ are invariant, so its
inflection polynomial \cite{Per01} is $Q_\cA R$ with a quartic $R$, the \emph{residual quartic}, which lies in the square
of the ideal of the zero scheme $Z$ \cite{GST16}.  We show that $Z$ is curvilinear, that every
non-invariant line meets $Z$ in length at most four (a Gauss-map argument), that repeated factors of $R$
are non-invariant new lines, that a doubled smooth conic is impossible (fixed-point indices of a
second-intersection map), and that a reduced residual has $\nu\le5$.  Local inflection identities exclude a
fourth-power residual and bound the other nonreduced ones.  For a reduced residual at defect five only three
configurations remain; an invariant residual line then carries a critically fixed rational map of degree
seven, resonances of the foliation become conditions on its multipliers, and we exclude the three
configurations with the Schwarz--Pick lemma, an exact catalog of critically fixed maps, a two-line
incidence certificate and exact elimination for maps with one free critical point.

\subsection{Notation}
All arrangements are arrangements of lines in $\PP^2=\PP^2_{\CC}$.  For an arrangement $\cA$ of $d$ lines
with defining polynomial $Q$ in $S=\CC[x,y,z]$, let $D_0(\cA)$ be the graded module of Jacobian syzygies of
$Q$ (derivations killing $Q$, graded by the degree of the coefficients), $\T\cA$ the associated rank-two
bundle and $\mdr(\cA)$ the least degree of a nonzero element of $D_0(\cA)$, so that
$D_0(\cA)_e=H^0(\T\cA(e))$.  For a point $P$ let $m_P$ be the number of lines through it; the
\emph{vertices} are the points with $m_P\ge2$.  For a line $H\in\cA$ let $r_H$ be the number of vertices on
$H$.  Put
\[
 b(\cA)=\sum_P(m_P-1),\qquad c_1(\T\cA)=1-d,\qquad c_2(\T\cA)=b(\cA)-(d-1).
\]
The arrangement is \emph{free} if $D_0(\cA)$ is a free $S$-module, that is, if $\T\cA$ splits as
$\cO(-d_1)\oplus\cO(-d_2)$ \cite{Sai80,OT92}; then $d_1+d_2=d-1$ and $d_1d_2=c_2(\T\cA)$ \cite{Ter81}, and
$(d_1,d_2)$ are the \emph{exponents}.  Pencils are free and are excluded whenever we discuss a
counterexample.  For a non-pencil arrangement, the
\emph{lattice} of $\cA$ is the simple matroid of rank three on the set of lines whose rank-two flats are the
vertices; $b$, $r_H$ and $m_P$ depend only on the lattice.  If $\cB\subset\cA$ and $L\in\cA$, we write
$r_{\cB}(L)$ for the number of points of $L$ on the lines of $\cB$ other than $L$ (so $r_\cB(L)$ is the
number of vertices of $\cB$ on $L$ when $L\in\cB$).  Adding a line $L$ to $\cB$ raises $b$ by $r_\cB(L)$,
and deleting a line $H$ from $\cA$ lowers $b$ by $r_H$.

For a lattice $\Lambda$ let $\cR(\Lambda)$ be the set of its complex realizations modulo $\PGL_3(\CC)$; a
realization is \emph{labelled}: it assigns a line to every element of $\Lambda$, and two realizations are
identified only by a projective transformation compatible with the labels.  The group $\Aut(\CC)$ acts on
$\cR(\Lambda)$.  We say that $\Lambda$ has \emph{one Galois orbit of realizations} if $\cR(\Lambda)$ is
nonempty and $\Aut(\CC)$ acts transitively on it.  Freeness and the dimensions of the graded pieces of $D_0$
are constant on such an orbit.

\paragraph{Organization.}  Section~\ref{sec:reductions} contains the reductions.
Sections~\ref{sec:crit}--\ref{sec:prop} establish the core bounds and reconstruction,
Sections~\ref{sec:census}--\ref{sec:search} describe the exact calculation, and Sections~\ref{sec:checks}
and~\ref{sec:reliability} its checks and the scope of the verification; Section~\ref{sec:class41} records a
by-product.
Sections~\ref{sec:nu} and~\ref{sec:18} treat defects two and three and eighteen lines.
Sections~\ref{sec:nu4} and~\ref{sec:sept} treat defect four and points of multiplicity seven, and
Section~\ref{sec:proof} completes the proof.  Section~\ref{sec:verification} lists the computations.  The
defect bound is proved in Appendices~\ref{app:A} and~\ref{app:B}.

\section{Reduction of Terao's conjecture to the seventeen-line windows}\label{sec:reductions}

We give the reductions to the seventeen-line windows explicitly.  They refine
\cite[Proposition~8.2]{GST16}, which reduces Terao's conjecture for seventeen
lines to $(C_7)$ at defect one and for eighteen lines to $(C_7)$ at defects two
and three.  The result of \cite{GST16} used here is the generic
splitting-type theorem for at most sixteen lines, together with its general
consequences \cite[Proposition~8.1(a),(c)]{GST16}; Proposition~\ref{red:17} also
uses the classification of Section~\ref{sec:sept} (Theorem~\ref{thm:sept}).  In
particular, the reductions do not use Terao's conjecture for seventeen lines as
an intermediate input to the eighteen-line argument.

\subsection{The input from sixteen lines and elementary tools}

For an arrangement $\cC$ of $d$ complex projective lines, set
$K_d=\lfloor(d-1)/2\rfloor$.  The following two consequences of the theorem for
at most sixteen lines apply to arrangements of \emph{any} cardinality:
\begin{enumerate}[label=\textup{(G\arabic*)},leftmargin=*]
\item\label{red:G1} If $K_d>6$ and two arrangements of $d$ lines have isomorphic
lattices, then one has $\mdr\le6$ if and only if the other does; when this
occurs their minimal degrees agree.
\item\label{red:G2} If a free arrangement has smaller exponent at most $7$,
then every arrangement with the same lattice is free.
\end{enumerate}
These are precisely \cite[Proposition~8.1(a),(c)]{GST16} with $S=6$, whose
hypotheses are supplied by \cite[Theorems~A and~B]{GST16}.  Stating the
cardinality-independent consequences is essential: a deletion of an
eighteen-line arrangement has seventeen lines, so the theorem restricted to
cardinalities at most sixteen would not by itself justify that step.

We use $D(\cC)$ for the module of logarithmic derivations, with Euler
derivation $E$.  The quotient $D(\cC)/SE$ identifies with $D_0(\cC)$.
We recall the line-restriction facts in a convenient form.

\begin{lemma}\label{red:restriction}
Let $\cC$ be an arrangement, let $H\in\cC$, and write
$\cC_H=\cC\setminus\{H\}$.  Then:
\begin{enumerate}[label=\textup{(\roman*)},leftmargin=*]
\item $\mdr(\cC_H)\le\mdr(\cC)\le\mdr(\cC_H)+1$.
\item If $\theta\in D(\cC)_e\setminus SE$ has nonzero projective restriction
to $H$, then $r_H\le e+1$.  If its projective restriction is zero, then
$\theta=\alpha_H\theta_0+hE$, with
$\theta_0\in D(\cC_H)_{e-1}\setminus SE$.
\item Let $L\notin\cC$ meet the lines of $\cC$ in $r_{\cC}(L)$ distinct
points.  If $|\cC|-r_{\cC}(L)>e$, then every derivation in $D(\cC)_e$ is
tangent to $L$.  In particular, a non-Euler derivation of degree $e$ on
$\cC$ then gives one of the same degree on $\cC\cup\{L\}$.
\end{enumerate}
\end{lemma}
\begin{proof}
For (i), use $D(\cC)\subseteq D(\cC_H)$ in one direction and multiply a
derivation on $\cC_H$ by $\alpha_H$ in the other, always working modulo
Euler derivations.  For (ii), the nonzero tangential restriction is a section
of $\cO_H(e+1)$ vanishing at every vertex on $H$.  If instead it vanishes
identically, choose coordinates with $H=\{z=0\}$.  The restriction of
$\theta$ is a polynomial multiple of the restriction of $E$, so subtracting
an Euler multiple makes every coefficient divisible by $z$.  The quotient
is tangent to every remaining line, since its equation is coprime to $z$;
the quotient cannot be an Euler multiple.  Finally, the normal restriction
of a degree-$e$ derivation to $L$ is a section of $\cO_L(e)$ which vanishes
to order at least $m_p(\cC)-1$ at each incidence point $p$ on $L$.
These orders sum to $|\cC|-r_{\cC}(L)$, proving (iii).  These are also
\cite[Lemmas~2.5--2.6 and Proposition~2.9]{GST16}.
\end{proof}

\begin{lemma}\label{red:chern-free}
An arrangement $\cC$ with $\mdr(\cC)\ge8$ is free in either of the following
two situations:
\[
 (|\cC|,c_2(\T\cC))=(17,64)\quad\text{or}\quad(18,72).
\]
Its exponents are respectively $(8,8)$ and $(8,9)$.
\end{lemma}
\begin{proof}
Put $F=\T\cC(8)$.  The assumption on $\mdr$ gives $H^0(F(-1))=0$.
For seventeen lines, $c_1(F)=c_2(F)=0$, and
\cite[Lemma~5.1]{GST16} gives $F\cong\cO^2$.
For eighteen lines, $c_1(F)=-1$ and $c_2(F)=0$.
Riemann--Roch gives $\chi(F)=1$, while Serre duality and
$F^*\cong F(1)$ give $h^2(F)=h^0(F(-2))=0$.
Thus $F$ has a nonzero section.  It has no divisorial zero, since such a
zero would produce a nonzero section of $F(-1)$.  Its finite zero scheme
has length $c_2(F)=0$, and consequently
\[
 0\longrightarrow\cO\longrightarrow F\longrightarrow\cO(-1)
 \longrightarrow0.
\]
This sequence splits because $H^1(\cO(1))=0$.
Untwisting proves the two assertions.  Equivalently, both are instances
of the Chern-class freeness criterion in \cite[Proposition~2.2]{GST16}.
\end{proof}

\subsection{Seventeen lines and the exclusion of multiplicity seven}

\begin{proposition}\label{red:17}
Suppose that $\cA$ and $\cB$ are
arrangements of seventeen complex projective lines with isomorphic
intersection lattices $\Lambda$, where $\cA$ is free and $\cB$ is not free.
Then $\cA$ has exponents $(8,8)$, and
\[
 b(\Lambda)=80,\qquad \mdr(\cB)=7,\qquad
 r_H\le8\ \text{for every }H,\qquad
 m_P\le6\ \text{for every }P.
\]
In particular, $\Lambda\in\mathcal W$ and $\cB$ is critical.
\end{proposition}
\begin{proof}
Write the exponents of $\cA$ as $(d_1,d_2)$ with $d_1\le d_2$.
They satisfy $d_1+d_2=16$.  By \ref{red:G2}, $d_1\le7$ is impossible;
hence $(d_1,d_2)=(8,8)$ and $c_2(\T\cA)=64$.  Chern classes are
combinatorial, so $c_2(\T\cB)=64$ and $b(\Lambda)=64+16=80$.
By \ref{red:G1}, $\mdr(\cB)\ge7$, while
Lemma~\ref{red:chern-free} excludes $\mdr(\cB)\ge8$.
Thus $\mdr(\cB)=7$.

Fix a line $H$ of the lattice and denote the corresponding deletions by
$\cA_H$ and $\cB_H$.  Lemma~\ref{red:restriction}(i) gives
$\mdr(\cB_H)\in\{6,7\}$.  If it were $6$, \ref{red:G1} would give
$\mdr(\cA_H)=6$, and then $\mdr(\cA)\le7$, a contradiction.
Consequently $\mdr(\cB_H)=7$ for every $H$.  A degree-seven derivation
of $\cB$ cannot have zero projective restriction to $H$, since the
factorization in Lemma~\ref{red:restriction}(ii) would give a degree-six
non-Euler derivation of $\cB_H$.  The same lemma therefore gives $r_H\le8$.

Next, $\cB$ is nearly free (this is used only in Remark~\ref{red:17-dependency}).  A nonzero minimal
section of
$\T\cB(7)$ has no divisorial zero and its zero scheme $Z$ satisfies
\[
 \operatorname{length}Z=c_2(\T\cB)-7(16-7)=64-63=1
\]
by \cite[Proposition~2.2]{GST16}.  Thus $Z$ is a single reduced point $p$,
and there is an exact sequence
\[
 0\longrightarrow\cO(-7)\longrightarrow\T\cB
 \longrightarrow\cI_p(-9)\longrightarrow0.
\]
Pulling this extension back along the standard resolution
$0\to\cO(-11)\to\cO(-10)^2\to\cI_p(-9)\to0$ gives
\[
 0\longrightarrow\cO(-11)\longrightarrow
 \cO(-7)\oplus\cO(-10)^2\longrightarrow\T\cB\longrightarrow0.
\]
Indeed, the intermediate extension splits because $H^1(\cO(3))=0$.
Taking all graded global sections is exact here, since
$H^1(\cO(t))=0$ on $\PP^2$ for all $t$.  Therefore
\[
 0\longrightarrow S(-11)\longrightarrow
 S(-7)\oplus S(-10)^2\longrightarrow D_0(\cB)\longrightarrow0
\]
is a minimal free resolution: all entries in its first map have positive
degree.  This proves that $\cB$ is nearly free with exponents $(7,10)$.

Finally, $\Lambda$ is not a pencil (a pencil of seventeen lines has exponents $(0,16)$) and $r_H\le8$
for every $H$, so $\Lambda$ is a $7$-candidate of defect $b(\Lambda)-79=1$, and $\cB$ is a critical
realization of it.  By Lemma~\ref{lem:septcount}(b) its multiplicities are at most $8$, and a point of
multiplicity $8$ would give defect $0$.  If $\Lambda$ had a point of multiplicity $7$, then by
Theorem~\ref{thm:sept} every realization of $\Lambda$ would be critical, in particular $\cA$; but
$\mdr(\cA)=8$.  Hence every multiplicity is at most $6$.
\end{proof}

\begin{remark}\label{red:17-dependency}
The last step uses Theorem~\ref{thm:sept}, whose proof (Section~\ref{sec:sept}) does not use this
proposition.  Since $\cB$ is nearly free with exponents $(7,10)$, the near-freeness criterion
\cite[Corollary~1.7]{DimcaNF} (a nearly free line arrangement with exponents $(d'_1,d'_2)$ and a point of
multiplicity at least $d'_1$ has only nearly free lattice-mates) would also exclude a point of
multiplicity at least $7$; the argument above keeps the proof within this paper.  Neither the lower bound
$r_H\ge4$ nor the general residual conditions for counterexamples are required for this reduction.
\end{remark}

\subsection{Eighteen lines: the reduction holds for every deletion}

\begin{proposition}\label{red:18}
Suppose that $\cA$ and $\cB$ are
arrangements of eighteen complex projective lines with isomorphic
intersection lattices $\Lambda$, where $\cA$ is free and $\cB$ is not free.
Then $\cA$ has exponents $(8,9)$, $b(\Lambda)=89$,
$\mdr(\cB)=7$, and every line has at most eight vertices.
For \emph{every} line $H$ of $\Lambda$, the two deletions satisfy
\[
 \mdr(\cB_H)=7,\qquad \mdr(\cA_H)\ge8,
 \qquad b(\Lambda\setminus H)=79+\nu_H,
 \qquad \nu_H=10-r_H.
\]
Moreover, every realization of $\Lambda\setminus H$ has minimal degree
at least $7$, and every line of this deletion has at most eight vertices.
Thus $\Lambda\setminus H$ is a $7$-candidate (Section~\ref{sec:c7}) with a critical and a
non-critical realization.
\end{proposition}
\begin{proof}
The exponents of $\cA$ sum to $17$.  By \ref{red:G2} their smaller member
cannot be at most $7$, so they are $(8,9)$ and $c_2(\T\cA)=72$.
It follows that $b(\Lambda)=72+17=89$.
Condition~\ref{red:G1} excludes $\mdr(\cB)\le6$, and
Lemma~\ref{red:chern-free} excludes $\mdr(\cB)\ge8$.
Consequently $\mdr(\cB)=7$.

For any $H$, Lemma~\ref{red:restriction}(i) gives
$\mdr(\cB_H)\in\{6,7\}$.  If it were $6$, \ref{red:G1}, applied to
the two \emph{seventeen-line} deletions, would give
$\mdr(\cA_H)=6$ and therefore $\mdr(\cA)\le7$, a contradiction.
Hence $\mdr(\cB_H)=7$.  As in the preceding proof, a degree-seven
derivation of $\cB$ must have nonzero projective restriction to every
$H$, and Lemma~\ref{red:restriction}(ii) gives $r_H\le8$.

Lemma~\ref{red:restriction}(i) also gives $\mdr(\cA_H)\ge7$.
If equality held, a degree-seven non-Euler derivation of $\cA_H$ would
be tangent to $H$ by Lemma~\ref{red:restriction}(iii), since
\[
 |\cA_H|-r_{\cA_H}(H)=17-r_H\ge9>7.
\]
It would then be a degree-seven non-Euler derivation of $\cA$, again
a contradiction.  This proves $\mdr(\cA_H)\ge8$.

Deleting $H$ lowers $b$ by $r_H$, whence
\[
 b(\Lambda\setminus H)=89-r_H=79+(10-r_H),\qquad
 c_2(\T{\cA_H})=73-r_H=63+(10-r_H).
\]
Deleting a line cannot increase the number of vertices on another line.
The deletion is not a pencil because it has a realization of minimal
degree $7$.  Finally, a realization of $\Lambda\setminus H$ with minimal
degree at most $6$ would, by \ref{red:G1}, force the same for $\cB_H$.
This proves all the assertions.
\end{proof}

\begin{remark}\label{red:18-dependency}
Proposition~\ref{red:18} does not assert that every deletion has defect
$2$ or $3$, or maximal multiplicity at most $6$.  Those properties require
the separate choice-of-line argument.  Its role is to choose $H$ with
$r_H\in\{7,8\}$ and with all multiplicities of $\Lambda\setminus H$ at most
$6$.  The proposition above then supplies the required incompatible
realizations in $\mathcal W_2$ or $\mathcal W_3$.  No case of Terao's
conjecture for seventeen lines and no case of generic splitting-type
invariance for seventeen lines is used in this reduction.
\end{remark}
\subsection{Candidates and condition \texorpdfstring{$(C_7)$}{(C7)}}\label{sec:c7}
Following \cite[Definition~3.2]{GST16}, a \emph{$7$-candidate} is a lattice on $17$ lines, not a pencil,
with at most $8$ vertices on every line.  Its \emph{defect} is $\nu=c_2-63=b-79$.  A realization is
\emph{critical} if its minimal degree is $7$.  Condition $(C_7)$ for a candidate says: if all its
realizations have $\mdr\ge7$, then all of them are critical or none is.  By
\cite[Theorems~A and~B]{GST16} the generic splitting type is combinatorial for at most $16$ lines; this gives
$(C_s)$ for all $s\le6$, since for an $s$-candidate all of whose realizations have $\mdr\ge s$ the generic
splitting type $k=\min\{\mdr,s+1\}$ equals $s$ exactly on the critical realizations ($(C_s)$ for $s\le6$ is
also established directly in the proofs of \cite[Theorems~A and~B]{GST16}, and used in this form in the
proof of \cite[Theorem~C]{GST16}).  Hence, by the induction \cite[Theorem~3.5]{GST16}, if $(C_7)$ holds for all
$7$-candidates, then the generic splitting type is combinatorial for $17$ and for $18$ lines.  For $\nu\ge1$, $\mathcal W_\nu$ denotes the set of
$7$-candidates of defect $\nu$ whose multiplicities are at most $6$; thus $\mathcal W_1=\mathcal W$.

\section{The critical member}\label{sec:crit}

Throughout this section $\cA$ is a critical realization of a lattice in $\mathcal W$; only
$b(\cA)=80$, $r_H\le8$ and $\mdr(\cA)=7$ are used.

\begin{lemma}\label{lem:gen}
There is $\theta\in D_0(\cA)_7$ such that every element of $D_0(\cA)$ of degree at most $9$ is a
polynomial multiple of $\theta$.  For every $H\in\cA$ the vector field $\theta$ is not everywhere
proportional to the Euler field $E$ along $H$.
\end{lemma}
\begin{proof}
Let $\theta\ne0$ be of degree $7$.  As a section of $F=\T\cA(7)$ it has no divisorial zeros, since
$H^0(F(-1))=0$.  We have $c_1(F)=-2$ and $c_2(F)=64-7\cdot16+49=1$, so
$0\to\cO\to F\to\cI_p(-2)\to0$ for a reduced point $p$.  For $e\le9$, $H^0(\cI_p(e-9))=0$, which
proves the first assertion.

Let $H=\{z=0\}$ and write $\theta(z)=zg$.  Then $\delta=\theta-gE$ is a derivation tangent to $\cA$
with $\delta(z)=0$, congruent to $\theta$ modulo $E$; in particular $\delta$ is not a multiple of $E$.  Suppose $\theta$, hence $\delta$, is proportional
to $E$ along $H$.  Then $(\delta(x),\delta(y))|_H=g'\cdot(x,y)$ with $g'$ a binary form of degree $6$.
By Ziegler's restriction theorem \cite{Zie89}, $\delta|_H$ is a derivation of the multiarrangement on $H$ whose
points are the vertices $P\in H$, with multiplicities $m_P-1$.  If $\beta_P$ is the linear form of
$P$ on $H$, this gives $\beta_P^{m_P-2}\mid g'$ for all $P$.  Since
$\sum_{P\in H}(m_P-2)=16-r_H\ge8>6$, $g'=0$.  Then $z$ divides $\delta$, and $\delta/z$ is tangent to
all lines of $\cA$, of degree $6$, and not a multiple of $E$.  Its component in $D_0(\cA)$ is nonzero
of degree $6$, a contradiction.
\end{proof}

\begin{lemma}\label{lem:del}
Let $S\subset\cA$ with $|S|=k\ge3$.  Then $\mdr(\cA\setminus S)\ge10-k$.
\end{lemma}
\begin{proof}
Let $0\ne\psi\in D_0(\cA\setminus S)$ have degree $e\le9-k$, and let $Q_S$ be the product of the lines
of $S$.  Then $Q_S\psi$ is tangent to every line of $\cA$, so $Q_S\psi=\theta'+fE$ with
$\theta'\in D_0(\cA)_{e+k}$.  Since $e+k\le9$, Lemma~\ref{lem:gen} gives $\theta'=q\theta$ with
$\deg q=e+k-7\le2$ ($q=0$ if $e+k<7$).  Along a line $H\in S$ the left side vanishes, so
$q\theta=-fE$ along $H$; by Lemma~\ref{lem:gen}, $q|_H=0$.  Hence $Q_S\mid q$, and $q=0$ because
$\deg q\le2<k$.  Then $Q_S\psi=fE$; as $x,y,z$ have no common factor, $Q_S\mid f$ and $\psi$ is a multiple
of $E$, so $\psi=0$ in $D_0$.
\end{proof}

\begin{proposition}\label{prop:need}
Let $\cC\subset\cA$ be a subarrangement of $n$ lines.  Then
\[
 b(\cC)\ \ge\ \beta(n):=\begin{cases}7n-43,&10\le n\le13,\\ 56,&n=14.\end{cases}
\]
If $n=12$ and $b(\cC)=41$, then $r_{\cC}(H)\le7$ for every $H\in\cC$.
\end{proposition}
\begin{proof}
Let $r=n-7$ and $G=\T\cC(r)$.  By Lemma~\ref{lem:del}, $H^0(G(-1))=0$.  One computes $c_1(G)=n-13$
and $c_2(G)=b(\cC)-(n-1)(r+1)+r^2=b(\cC)-7n+43$.  If $H^0(G)\ne0$, a section has no divisorial
zeros and its zero scheme has length $c_2(G)\ge0$.  If $H^0(G)=0$, then
$h^2(G)=h^0(G(10-n))$, which is $0$ for $n\ge11$ and equals $h^0(G)=0$ for $n=10$; so
$\chi(G)=2+c_1(c_1+3)/2-c_2\le0$, which gives $c_2(G)\ge1$ in all four cases.  For $n=14$ we have
$\mdr(\cC)\ge7$; for $G'=\T\cC(6)$, $c_1=-1$, $c_2=b(\cC)-55$, $H^0(G')=0$ and
$h^2(G')=h^0(G'(-2))=0$, so $\chi(G')=1-c_2\le0$.  The last assertion follows from
$b(\cC\setminus H)=41-r_\cC(H)\ge\beta(11)=34$.  (For $n=13,14$ the bounds also follow from
$r_H\le8$; the content of the proposition is in the cases $n\le12$.)
\end{proof}

\section{Cores of minimal weight}\label{sec:cores}

Let $\Lambda\in\mathcal W$ have a critical realization.  For a set $S$ of five lines let $f(S)$ be the
number of vertices all of whose lines belong to $S$, and
$\varphi(S)=\sum_{H\in S}(r_H-7)-f(S)$.  Deleting lines one at a time gives
\begin{equation}\label{eq:delid}
 b(\Lambda\setminus S)=80-\sum_{H\in S}r_H+f(S)=45-\varphi(S).
\end{equation}

\begin{proposition}[Twelve-line core]\label{prop:core}
Every $12$-line sublattice of $\Lambda$ has $b\ge41$, and some $12$-line sublattice has $b\le46$.
\end{proposition}
\begin{proof}
The first statement is Proposition~\ref{prop:need}.  Let $t_m$ be the number of points of multiplicity
$m$ and $W=t_4+3t_5+6t_6$.  From $\sum\binom m2t_m=136$ and $\sum(m-1)t_m=80$ we get
$t_3=56-3W+3t_5+8t_6$, $t_2=3W-32-t_5-3t_6$, hence $W\ge11$, and $\sum_H(r_H-7)=W-15$.  A line lies in
$\binom{16}4=1820$ five-element sets, and a vertex of multiplicity $2,3,4,5$ is contained in
$455,91,13,1$ of them.  Therefore
\begin{align*}
 \sum_{|S|=5}\varphi(S)&=1820(W-15)-455t_2-91t_3-13t_4-t_5\\
 &=715W+220t_5+715t_6-17836\ \ge\ -9971\ >\ -2\binom{17}5 ,
\end{align*}
so some $S$ has $\varphi(S)\ge-1$, and $b(\Lambda\setminus S)\le46$ by \eqref{eq:delid}.
\end{proof}

\begin{definition}
A \emph{minimal core} of $\Lambda$ is a $12$-line sublattice $M=\Lambda\setminus S$ with $b(M)$ minimal
among all $12$-line sublattices.  We write $b(M)=41+\delta$; by Proposition~\ref{prop:core},
$0\le\delta\le5$.  For a line $L\in S$ put $v(L)=r_M(L)$.
\end{definition}

\begin{proposition}\label{prop:min}
Let $M=\Lambda\setminus S$ be a minimal core.
\begin{enumerate}[label=\textup{(\alph*)}]
\item $\sum_{L\in S}v(L)\le39-\delta$.
\item For $L\in S$ and $H\in M$: $r_M(H)+\epsilon(H,L)\le v(L)$, where $\epsilon(H,L)=0$ if $L$ meets $H$
in a vertex of $M$ and $1$ otherwise.
\item $r_M(H)\le\lfloor(39-\delta)/5\rfloor$ for every $H\in M$; that is, $r_M(H)\le7$, and
$r_M(H)\le6$ if $\delta=5$.
\item For every sublattice $N$ of $\Lambda$ with $n$ lines, $10\le n\le14$: $b(N)\ge\beta(n)$; and
$b(N)\ge41+\delta$ if $n=12$.
\end{enumerate}
\end{proposition}
\begin{proof}
(a) Sum $m_P-1$ over the points of the plane.  For a point on a line of $M$ the contribution is
$(m_M(P)-1)+s(P)$, where $s(P)$ is the number of lines of $S$ through $P$; the other points contribute
$s(P)-1\ge0$ when $s(P)\ge1$.  Hence $80=b(M)+\sum_{L\in S}v(L)+\varepsilon$ with $\varepsilon\ge0$.
(b) $M\cup\{L\}\setminus\{H\}$ is a $12$-line sublattice, and
$b(M\cup L\setminus H)=b(M)+v(L)-r_{M\cup L}(H)$ with $r_{M\cup L}(H)=r_M(H)+\epsilon(H,L)$; use
minimality.  (c) By (a) some $L$ has $v(L)\le\lfloor(39-\delta)/5\rfloor$; apply (b).  (d) is
Proposition~\ref{prop:need} and minimality.
\end{proof}

\section{Reconstruction and rigid sublattices}\label{sec:prop}

\begin{lemma}\label{lem:join}
Let $\Lambda\in\mathcal W$, let $N\subset\Lambda$ be a sublattice and $L\in\Lambda\setminus N$.
\begin{enumerate}[label=\textup{(\alph*)}]
\item If $|N|\ge13$, then $L$ contains at least two vertices of $N$.
\item If $|N|=12$ and $L$ contains at most one vertex of $N$, then $r_N(L)=8$.
\end{enumerate}
Consequently, if $\Lambda$ has a critical realization, $M$ is a minimal core and the lines of $S$ are
ordered as $L_1,\dots,L_5$ with $v(L_1)\le\dots\le v(L_5)$, then $L_j$ contains at least two vertices
of $N_{j-1}=M\cup\{L_1,\dots,L_{j-1}\}$ for every $j$.
\end{lemma}
\begin{proof}
If $L$ contains no vertex of $N$, then $r_N(L)=|N|\ge12$.  If it contains exactly one, of multiplicity
$m\le5$ in $N$, then $r_N(L)=1+|N|-m\ge|N|-4$.  Since $r_N(L)\le r_\Lambda(L)\le8$, this proves (a) and (b).
For the last statement, by Proposition~\ref{prop:min}(a) not all $v(L)$ are equal to $8$, so
$v(L_1)\le7$ and (b) applies to $L_1$; for $j\ge2$ use (a).
\end{proof}

\begin{proposition}\label{prop:rigid}
Let $\Lambda\in\mathcal W$ and let $N\subset\Lambda$ be a sublattice with at least $13$ lines, or a
$12$-line sublattice such that some line of $\Lambda\setminus N$ contains two vertices of $N$.  If $N$ has
one Galois orbit of realizations, then $\Lambda$ has at most one Galois orbit of realizations.
\end{proposition}
\begin{proof}
Let $\cA,\cA'$ be realizations of $\Lambda$ and $\cB,\cB'$ their subarrangements with lattice $N$.
There are $\sigma\in\Aut(\CC)$ and $g\in\PGL_3(\CC)$ with $\cB'=g\sigma(\cB)$, compatibly with the
labels.  By Lemma~\ref{lem:join}(a) the remaining lines can be added so that each of them is the join of two
vertices of the arrangement constructed before it (start with a line through two vertices of $N$ if
$|N|=12$); these vertices are determined by the labelled lattice, and joins commute with $g$ and $\sigma$.  Hence $\cA'=g\sigma(\cA)$.
\end{proof}

\section{The census of cores}\label{sec:census}

Let $\mathcal K_\delta$ be the set of lattices $M$ on $12$ lines with multiplicities at most $6$,
$b(M)=41+\delta$ and at most $7$ points on every line (at most $6$ if $\delta=5$).  By
Propositions~\ref{prop:core} and~\ref{prop:min}, a minimal core of a lattice in $\mathcal W$ with a
critical realization belongs to some $\mathcal K_\delta$, $0\le\delta\le5$.

The lattices on $12$ lines with multiplicities at most $6$, $b=41+\delta$ and at most $8$ points on
every line were generated with the programs of \cite{AS} (canonical augmentation \cite{McK98} of the
configurations of points of multiplicity at least $4$, completion by triple points, canonical forms with
\texttt{nauty} \cite{MP14}).  Their numbers are
\[
 767\ (r\le7),\quad 4\,025,\quad 18\,569,\quad 68\,427,\quad 211\,826,\quad 540\,173\qquad(\delta=0,\dots,5).
\]
By the number of sextuple points (none and one) they split as $696+70$, $3\,733+291$, $17\,453+1\,115$,
$66\,106+2\,321$, $208\,830+2\,996$, $538\,293+1\,879$; the remaining four lattices have two sextuple
points.
Restricting the number of points per line as in the definition of $\mathcal K_\delta$ leaves
\[
 |\mathcal K_\delta|=767,\ 3\,650,\ 14\,226,\ 40\,142,\ 88\,669,\ 5 .
\]
The completeness of these lists is part of the computation.  Two sectors can be checked by hand.  A
lattice with two sextuple points consists of two pencils with a common line and one transversal, or of
two disjoint pencils; this gives exactly four lattices, in the classes $\delta=0,1,2,5$, and only the
first two belong to $\mathcal K_\delta$.  A lattice of $\mathcal K_5$ has $20$ triple points, five on
every line, and six double points; it is a Steiner triple system on $13$ points with one point deleted,
and the two such systems have $1+4$ orbits of points.

\section{Realization spaces}\label{sec:real}

For a lattice $N$ we compute $\cR(N)$ exactly.  A frame of four lines in general position is fixed; the
other lines are placed one at a time.  A line through two already constructed vertices is their join; a
line through one constructed vertex $P$ is $u+tv$, where $u,v$ are two placed lines through $P$ and $t$
is a new parameter (the missing member $v$ of the pencil is a line of the arrangement); a line through no
constructed vertex is $(a:b:1)$ with two new parameters (it cannot pass through the vertex $(0:0:1)$
of the first two frame lines).  The coordinates are polynomials over $\ZZ$ in the parameters.  The
realization space is the set of parameters where the determinants of the dependent triples vanish and
those of the independent triples, and the common factors removed during the construction, do not.

\emph{Fast engine} (the program of \cite{AS,GST16}, at most one parameter).  With no parameter the test
is rational arithmetic.  With one parameter, let $G$ be the gcd of the equations.  If $G=0$ the
realization space is a cofinite subset of the line (status \textsf F), or empty if an independent triple
is identically dependent.  Otherwise the realizations are the roots of the irreducible factors of $G$
dividing no forbidden polynomial; the status is \textsf E (empty), \textsf R (one factor, i.e.\ one
Galois orbit) or \textsf M (several factors).

\emph{General engine} (new; any number of parameters).  The equations are factored and the factors
dividing a forbidden polynomial are discarded.  The minimal primes of the ideal generated by the cleaned
equations are computed with \textsc{Singular} (\texttt{minAssGTZ}); a prime containing an irreducible
factor of a forbidden polynomial is discarded.  The status is \textsf E if no prime remains, \textsf R
if exactly one prime remains and it has dimension $0$, and \textsf O otherwise; the dimensions of the
components are recorded.
A zero-dimensional prime over $\QQ$ has a finite field extension of $\QQ$
as its coordinate ring. Its embeddings in $\CC$ form a single orbit under
$\Aut(\CC)$, and none of the retained nonvanishing conditions removes
one conjugate without removing them all. Thus status \textsf R means one
Galois orbit of labelled realizations. Several zero-dimensional components
are not assigned this status merely because each is geometrically rigid.

\begin{compres}\label{cr:status}
The statuses of the lattices of $\mathcal K_\delta$ are as follows (\textsf O$_k$: one component, of
dimension $k$; \textsf O$_{0,0}$: two orbits).
\begin{center}\small
\begin{tabular}{crrrrl}\toprule
$b$&$|\mathcal K_\delta|$&\textsf E&\textsf R&\textsf F&other\\\midrule
41&767&718&41&5&\textsf O$_1$: 2, \textsf O$_3$: 1\\
42&3\,650&3\,367&225&41&\textsf O$_1$: 8, \textsf O$_2$: 8, \textsf O$_4$: 1\\
43&14\,226&12\,703&1\,156&278&\textsf O$_1$: 34, \textsf O$_2$: 50, \textsf O$_3$: 5\\
44&40\,142&34\,656&4\,427&866&\textsf O$_1$: 84, \textsf O$_2$: 97, \textsf O$_3$: 10, \textsf O$_4$: 1, \textsf O$_{0,0}$: 1\\
45&88\,669&73\,781&12\,119&2\,257&\textsf O$_1$: 331, \textsf O$_2$: 173, \textsf O$_3$: 6, \textsf O$_4$: 2\\
46&5&5&0&0&\\\bottomrule
\end{tabular}
\end{center}
In particular no lattice of $\mathcal K_5$ is realizable: twelve lines cannot have $20$ triple points
with five on every line.
\end{compres}

\section{The extension search}\label{sec:search}

Let $\Lambda\in\mathcal W$ have a critical realization, $M$ a minimal core with $b(M)=41+\delta$ and
$L_1,\dots,L_5$ the other lines, ordered by $v$ (Lemma~\ref{lem:join}).  Put $N_j=M\cup\{L_1,\dots,L_j\}$.
If $M$ has status \textsf R, then $\Lambda$ has one Galois orbit of realizations by
Proposition~\ref{prop:rigid}.  Status \textsf E is impossible.  So assume $M$ is one of the
$4\,261$ realizable lattices of Computational result~\ref{cr:status} whose status is not \textsf R (the
\emph{roots}).

A \emph{node} is a pair $(N,M)$ of a lattice $N$ on $n$ lines, $12\le n\le16$, and a marked $12$-line
sublattice $M$.  Its \emph{children} are the pairs $(N\cup\{L\},M)$, where the new line $L$ is given by
its set $A$ of \emph{anchors}, the vertices of $N$ on $L$ (pairwise without common line; the lines of $N$
through no anchor meet $L$ in new double points), subject to the following conditions, all of which hold
for $N=N_{j}$, $L=L_{j+1}$:
\begin{enumerate}[label=\textup{(C\arabic*)},leftmargin=*]
\item $|A|\ge2$; the anchors have multiplicity at most $5$ in $N$; every line of $N\cup L$ has at most
$8$ points (Lemma~\ref{lem:join}, \eqref{eq:window});
\item $v(L)\ge v(L')$ for the lines $L'$ added before, and
$\sum_{L'}v(L')+(5-j)\,v(L)\le39-\delta$ (ordering and Proposition~\ref{prop:min}(a));
\item $r_M(H)+\epsilon(H,L)\le v(L)$ for all $H\in M$ (Proposition~\ref{prop:min}(b));
\item $b(N\cup L)\ge\beta(n+1)$, with $\beta(15)=64$, $\beta(16)=72$, $\beta(17)=80$ (at most $8$ points
on each later line); $r_N(L)\ge v_0$ and $b(N\cup L)+(16-n)\,v_0\le80$, where $v_0$ is the maximum of the numbers $r_M(H)$, $H\in M$, and $v(L')$, $L'$ a line added before $L$
(by (C2) the largest $v(L')$ is that of the line added last, so $v_0$ depends only on $N$ with its marked
sublattice $M$ and its marked added lines, not on their order); and $b(N\cup L)=80$ if $n+1=17$ (by
Proposition~\ref{prop:min}(b) and the ordering, every line $L''$ added from now on has
$r(L'')\ge v(L'')\ge v_0$ at the time it is added);
\item every sublattice of $N\cup L$ with $10\le s\le14$ lines has $b\ge\beta(s)$, and $b\ge41+\delta$ if
$s=12$ (Proposition~\ref{prop:min}(d)).
\end{enumerate}
Here $v(L)$ is computed from the anchors: it is the number of anchors lying on a line of $M$, plus the
number of lines of $M$ through no anchor.  Two children are identified if they are isomorphic as
lattices with a marked sublattice and a marked set of added lines.

The status of each child is computed exactly.  If the parent has status \textsf F, with chart
$\ell_i(t)$, the new line is the join of two anchors, the conditions are the collinearity of the other
anchors and the non-incidence with the other vertices, and the status is obtained from the gcd of the
collinearity equations as for the fast engine; the same is done modulo each irreducible factor for a
parent of status \textsf M.  For the other parents the engines of Section~\ref{sec:real} are applied to the
child.  A child of status \textsf E is discarded.  A child of status \textsf R is terminal: by
Proposition~\ref{prop:rigid} every lattice of $\mathcal W$ containing it has at most one Galois orbit
of realizations.  The other children are called \emph{open} and are expanded in turn.  (The root of
status \textsf O$_{0,0}$ and its descendants are treated by the general engine.)

\begin{compres}\label{cr:search}
For every root, the search ends before reaching $17$ lines: no node on $16$ lines is open.  The numbers of
open nodes are
\begin{center}\small
\begin{tabular}{crrrrr}\toprule
$b(M)$&roots&open on $13$ lines&on $14$&on $15$&on $16$\\\midrule
41&8&7&13&4&0\\
42&58&140&155&52&0\\
43&367&789&325&1&0\\
44&1\,059&771&108&0&0\\
45&2\,769&675&25&7&0\\\bottomrule
\end{tabular}
\end{center}
About $1.6$ million children were examined; no child of status \textsf M occurred.
\end{compres}

\begin{proof}[Proof of Theorem~\ref{thm:main}]
Let $\Lambda\in\mathcal W$ have a critical realization, and choose $M$, $L_1,\dots,L_5$ as above.  The
lattice $M$ belongs to $\mathcal K_\delta$ and is realizable.  If it has status \textsf R we are done.
Otherwise it is a root.  Every $N_j$ is realizable, so its status is not \textsf E, and
$(N_{j+1},M)$ satisfies (C1)--(C5), so it is (isomorphic to) a child of $(N_j,M)$.  If none of
$N_1,\dots,N_4$ had status \textsf R, then by induction $(N_1,M),\dots,(N_4,M)$ would be open nodes, and
$(N_4,M)$ would be an open node on $16$ lines, contradicting Computational result~\ref{cr:search}.  So
some $N_j$ has status \textsf R, and Proposition~\ref{prop:rigid} applies.
\end{proof}

\begin{corollary}\label{cor:terao17}
Terao's conjecture holds for arrangements of seventeen lines in $\PP^2_{\CC}$.
\end{corollary}
\begin{proof}
A counterexample would, by Proposition~\ref{red:17}, give a lattice in
$\mathcal W$ with a free realization and a critical realization.
Theorem~\ref{thm:main} excludes this possibility.
\end{proof}

\section{Validation checks}\label{sec:checks}
The programs, the records and the logs of the following checks are in the ancillary directory
\texttt{gst1718\_computation} (directories \texttt{check} and \texttt{logs}; the script
\texttt{check/run\_checks.sh} runs them).
\begin{itemize}[leftmargin=*]
\item \emph{Census.}  The two sectors of Section~\ref{sec:census} (the lattices with two sextuple
points and the class $\mathcal K_5$) were checked by hand.  The completeness of the census rests on the
generator (Section~\ref{sec:reliability}).
\item \emph{Statuses of the cores.}  The fast and the general engine were compared on all $146\,160$
lattices of $\mathcal K_0,\dots,\mathcal K_4$ which have a construction with at most one parameter, after
a random relabelling of the lines (which changes the frame and the construction): no disagreement.  For
the $1\,294$ other lattices the general engine was run three times with different labellings, with
identical statuses and dimensions.  The five lattices of $\mathcal K_5$ were recomputed with six
labellings, and an exhaustive search finds no realization over $\mathbb F_p$ for
$p=13,17,19,23$.
\item \emph{Finite fields.}  An independent backtracking search for realizations over $\mathbb F_{101}$,
$\mathbb F_{103}$ and $\mathbb F_{107}$ (the program of \cite{AS}) was run on the lattices of $\mathcal K_0$,
$\mathcal K_1$, $\mathcal K_2$ with a one-parameter construction in the frame chosen by that program
($18\,392$ lattices, $16\,673$ of them of status \textsf E), on a quarter of $\mathcal K_3$ ($9\,844$ lattices, $8\,572$ of status \textsf E) and on an
eighth of $\mathcal K_4$ ($10\,591$ lattices, $9\,061$ of status \textsf E).  No lattice of
status \textsf E has a realization over any of the three fields, and every lattice of status \textsf F,
\textsf M or \textsf O has one over each of them.  The lattices of status \textsf R have one over at least
one of the fields, except $157$ lattices whose Galois orbit is defined by an irreducible polynomial of degree
$3$, $4$ or $6$ without roots modulo the three primes.  The same search was run on a random sample of
$4\,147$ children of status \textsf E of the extension search ($13$ to $16$ lines); the $4\,000$ of them
with a one-parameter construction were searched, and none has a realization over any of the three
fields.
\item \emph{Incremental statuses.}  For all roots with $b(M)=41,42$ and for random samples of $60$, $60$
and $80$ roots with $b(M)=43,44,45$, every status computed on a chart of the parent was compared with the
status computed from scratch by the engines of Section~\ref{sec:real}: no disagreement.
\item \emph{Independent implementation.}  A second program for the search of Section~\ref{sec:search} was
written from the description (C1)--(C5) alone, with its own generation of the children, its own
implementation of the conditions and of the identification of isomorphic children, and with the status of
every child computed from scratch.  On the $4\,261$ roots it examined $1\,084\,865$ children (after
identification) and found the same numbers of open nodes as in Computational result~\ref{cr:search}
on $13$, $14$, $15$ and $16$ lines, in every class.  (Both programs use the same two engines.)
\item \emph{Positive control.}  The test list \texttt{check/sub24.g6} contains $23$ known realizable lattices; thirteen have multiplicities at most $6$ and are critical
(\texttt{check/critical\_d0.py}).  For each of them, every
sublattice satisfies the bounds of Proposition~\ref{prop:need}; their minimal cores (all of them, or a
random sample of $25$ when there are more, with up to four admissible orders of the other lines) belong to
the census, are never of status \textsf E, satisfy Proposition~\ref{prop:min}, and along each of these orders of the other lines the conditions
(C1)--(C5) hold and the path reaches a node of status \textsf R on
$13$, $14$ or $15$ lines.  Twelve of these lattices have $\delta=0$, and all their minimal cores have
status \textsf R, so they do not exercise the extension search itself.  The lattice obtained from
$\cA(G(6,6,3))$ by deleting a line has $\delta=4$, two sextuple points and $1\,960$ minimal cores; the
sample of $25$ contains $18$ cores of status \textsf F or \textsf O, and for these the search produces the
corresponding paths.
\end{itemize}

\section{A by-product: the lattices with a minimal core of weight 41}\label{sec:class41}
If the search of Section~\ref{sec:search} is continued through the nodes of status \textsf R, and all
realizable lattices of $\mathcal K_\delta$ are taken as roots, it produces every lattice of $\mathcal W$
which satisfies the necessary conditions of a critical realization and has a minimal core in
$\mathcal K_\delta$.  We did this for $\delta=0$ ($49$ roots).  The result is a list of exactly $16$
lattices on $17$ lines.  Each has one Galois orbit of realizations, and all of them are critical: they are not free (exact
computation), and $\dim D_0=0,0,1,3,6$ in degrees $5,\dots,9$ modulo two primes (which can only increase
these dimensions; \texttt{check/critical\_d0.py}), so $\mdr=7$ by Lemma~\ref{red:chern-free}.  Hence they are
nearly free with exponents $(7,10)$: the sequence $0\to\cO\to\T\cA(7)\to\cI_p(-2)\to0$ of the proof of
Lemma~\ref{lem:gen}, with $p$ a reduced point, gives the resolution
$0\to S(-11)\to S(-7)\oplus S(-10)^2\to D_0(\cA)\to0$.  In particular the lattices of $\mathcal W$ with a critical realization and a minimal core of
weight $41$ are exactly these sixteen.  Twelve of them are the lattices $\Gamma_1,\dots,\Gamma_{12}$ of
\cite[Table~3]{MSS} (Remark~\ref{rem:mss}); they are realized over $\QQ(\sqrt{-3})$.
The other four are realized over $\QQ$; their equations are:
\begin{itemize}[leftmargin=*]\small
\item $(t_2,\dots,t_6)=(20,18,2,2,2)$: $x+y$, $z$, $x-z$, $y+z$, $x+2y+z$, $2x+3y+z$, $3x+4y+z$, $x+y+z$,
$2x+2y+z$, $3x+3y+z$, $4x+4y+z$, $x$, $y$, $2y+z$, $3x+2y+z$, $x+2y$, $3x+2y$;
\item $(23,15,3,2,2)$: $2y+3z$, $3x+2y+3z$, $2x+2y+3z$, $x+2y+3z$, $3x+4y+6z$, $x$, $y$, $x+y+2z$, $z$,
$y+z$, $x+y+z$, $x+z$, $y+3z$, $4y+3z$, $3x+4y+3z$, $3x+5y+6z$, $3x+2y+6z$;
\item $(26,8,6,5,0)$: $x+y-z$, $2x+2y+z$, $2x-y+z$, $2x+z$, $x+y$, $z$, $y-z$, $x+z$, $2x+y+z$, $x+y+z$,
$x$, $y$, $y+z$, $3x+y+z$, $2y-z$, $x-y$, $x+2y$;
\item $(23,10,8,2,1)$: $x+y-z$, $2x+2y+z$, $2x-y+z$, $2x+z$, $x+y$, $z$, $y-z$, $x+z$, $2x+y+z$, $x+y+z$,
$x$, $y$, $2x+y$, $x+2y$, $x-y$, $3x+y+z$, $2y-z$.
\end{itemize}
This also shows that the search reaches $17$ lines when realizable lattices exist: the twelve lattices
$\Gamma_1,\dots,\Gamma_{12}$ are recovered without being used as input.

\begin{remark}\label{rem:mss}
These sixteen lattices are, up to isomorphism, the sixteen lattices of \cite[Table~3]{MSS}.  We
compared canonical forms (computed with \texttt{nauty}) of the line--vertex incidence graphs with the
public data of \cite{MSS}: the twelve lattices over $\QQ(\sqrt{-3})$ are their lattices
$\Gamma_1,\dots,\Gamma_{12}$, and the four rational lattices above are, in the order listed,
$R^b_{17}$, $N_1$, $N_2$ and $R^a_{17}$.  In \cite{MSS} these sixteen are shown to be the only
lattices with a non-free realization in their window.  That window is narrower than $\mathcal W$:
it also requires, among other conditions, that two vertices of multiplicities $6,5$ or $6,6$ lie on a
common line, a condition derived from the free member of a counterexample.  The lattice of
$\cA(G(6,6,3))$ minus a line, used as a positive control in Section~\ref{sec:checks}, lies in
$\mathcal W$ and has a critical realization; it has two sextuple points not on a common line, and its
minimal cores have weight $45$.  The comparison was made after our computation, and the theorem of
\cite{MSS} is not used in our proofs.
\end{remark}

\section{Defects two and three}\label{sec:nu}

We use the candidates, the defect and the classes $\mathcal W_\nu$ of Section~\ref{sec:c7}; thus
$\mathcal W_1=\mathcal W$.

\begin{theorem}\label{thm:nu}
Let $\nu\in\{2,3\}$ and $\Lambda\in\mathcal W_\nu$.  If $\Lambda$ has a critical realization, then
$\Lambda$ has one Galois orbit of realizations.  In particular $(C_7)$ holds for every $7$-candidate of
defect $2$ or $3$ with multiplicities at most $6$.
\end{theorem}

The proof is that of Theorem~\ref{thm:main}, with the following changes.  Let $\Lambda\in\mathcal W_\nu$
have a critical realization $\cA$.

\paragraph{The critical member.}
Now $F=\T\cA(7)$ has $c_1=-2$ and $c_2=\nu$, and the section $\theta$ gives
$0\to\cO\to F\to\cI_Z(-2)\to0$ with $Z$ of length $\nu\ge1$; again $H^0(\cI_Z(e-9))=0$ for $e\le9$.  The
rest of the proof of Lemma~\ref{lem:gen} uses only $r_H\le8$ and $\mdr(\cA)=7$.  So
Lemmas~\ref{lem:gen}, \ref{lem:del} and Proposition~\ref{prop:need} hold for every $\nu\ge1$.  Moreover,
deleting a line lowers $b$ by the number of its points, which is at most $8$; hence
\begin{equation}\label{eq:triv}
 b(N)\ \ge\ 79+\nu-8\,(17-n)\qquad\text{for every sublattice $N$ of $\Lambda$ with $n$ lines.}
\end{equation}
Put $\beta_\nu(n)=\max\{\beta(n),\,79+\nu-8(17-n)\}$ for $10\le n\le14$ and $\beta_\nu(n)=79+\nu-8(17-n)$
for $n\ge15$.  For $\nu=2,3$: $\beta_\nu(10)=27$, $\beta_\nu(11)=34$, $\beta_\nu(12)=39+\nu$,
$\beta_\nu(13)=47+\nu$, $\beta_\nu(14)=55+\nu$.

\paragraph{Cores.}
For a set $S$ of five lines, $b(\Lambda\setminus S)=44+\nu-\varphi(S)$.  With $W=t_4+3t_5+6t_6$ one finds
$t_2=3W-35+3\nu-t_5-3t_6\ge0$, so $W\ge10$ for $\nu=2$ and $W\ge9$ for $\nu=3$,
$\sum_H(r_H-7)=W+3\nu-18$, and
\[
 \sum_{|S|=5}\varphi(S)=715W+220t_5+715t_6+4186\nu-22022\ \ge\
 \begin{cases}-6500,&\nu=2,\\-3029,&\nu=3.\end{cases}
\]
Since $\binom{17}5=6188$, some $S$ has $\varphi(S)\ge-1$ if $\nu=2$ and $\varphi(S)\ge0$ if $\nu=3$; in both
cases $b(\Lambda\setminus S)\le47$.  So a minimal core $M=\Lambda\setminus S$ has $b(M)=41+\delta$ with
$0\le\delta\le6$.

\paragraph{Minimal cores.}
Proposition~\ref{prop:min} becomes:
\begin{enumerate}[label=\textup{(\alph*)}]
\item $\sum_{L\in S}v(L)\le38+\nu-\delta$;
\item $r_M(H)+\epsilon(H,L)\le v(L)$ for $L\in S$, $H\in M$;
\item $r_M(H)\le\rho(\nu,\delta):=\min\{8,\lfloor(38+\nu-\delta)/5\rfloor\}$, and $r_M(H)\le7$ if
$\delta=0$;
\item $b(N)\ge\beta_\nu(n)$ for every sublattice $N$ with $10\le n\le14$ lines, and $b(N)\ge41+\delta$
if $n=12$;
\item $v(L)\ge6+\nu-\delta$ for every $L\in S$; consequently $\delta\ge\nu-2$.
\end{enumerate}
Here (e) is \eqref{eq:triv} for $M\cup\{L\}$, and then $5(6+\nu-\delta)\le38+\nu-\delta$ by (a).

\paragraph{Reconstruction.}
Lemma~\ref{lem:join}(a), (b) and Proposition~\ref{prop:rigid} hold as stated.  The lines of $S$ can
still be ordered so that $v(L_1)\le\dots\le v(L_5)$ and every $L_j$ contains two vertices of $N_{j-1}$.
Indeed, a line of $S$ with at most one vertex of $M$ has $v=8$ and passes through exactly one vertex
of $M$, of multiplicity $5$ in $M$; two such lines pass through different vertices, and $k$ points of
multiplicity $5$ lie on at least $5k-\binom k2$ lines, which is more than $12$ for $k=4,5$ (the only larger possibilities since $|S|=5$).  So at most
three lines of $S$ are of this kind.  If some $v(L)$ is at most $7$, take for $L_1$ a line with minimal
$v$; otherwise all $v(L)$ are equal to $8$ and we take for $L_1$ a line with two vertices of $M$.  For
$j\ge2$ use Lemma~\ref{lem:join}(a).

\paragraph{Census.}
Let $\mathcal K^\nu_\delta$ be the set of lattices on $12$ lines with multiplicities at most $6$,
$b=41+\delta$ and at most $\rho(\nu,\delta)$ points on every line (at most $7$ if $\delta=0$).  A
minimal core of $\Lambda$ belongs to $\mathcal K^\nu_\delta$ with $\max(0,\nu-2)\le\delta\le6$.
\begin{itemize}[leftmargin=*]
\item $\nu=2$: the bound on the number of points of a line is $7$ for $\delta\le5$, so
$\mathcal K^2_\delta=\mathcal K_\delta$ for $\delta\le4$, and
$\mathcal K^2_5$ consists of the $143\,699$ lattices with $b=46$ and at most $7$ points on every line.
For $\delta=6$, $\rho=6$ and $\mathcal K^2_6$ is empty: a lattice on $12$ lines with $b=47$ and at most
$6$ points on every line has $\sum_Pm_P\le72$, hence at most $25$ vertices and $\sum_P(m_P-2)\ge22$,
whereas $\sum_P\binom{m_P-1}2=66-47=19$ and $m-2\le\binom{m-1}2$.
\item $\nu=3$: $\delta\ge1$; $\rho=8$ for $\delta=1$ and $\rho=7$ for $2\le\delta\le6$.  So
$\mathcal K^3_1$ consists of all $4\,025$ lattices with $b=42$ of Section~\ref{sec:census},
$\mathcal K^3_\delta=\mathcal K_\delta$ for $\delta=2,3,4$, $\mathcal K^3_5=\mathcal K^2_5$, and
$\mathcal K^3_6$ consists of the $155\,619$ lattices with $b=47$ and at most $7$ points on every line.
\end{itemize}
The two new lists were generated with the programs of Section~\ref{sec:census}; the list for $b=46$
coincides with the sublist of the $540\,173$ lattices of Section~\ref{sec:census} with at most $7$ points
on every line.  The statuses of the lattices not covered by Computational result~\ref{cr:status} are as
follows (the row $b=46$ also contains the five lattices of $\mathcal K_5$, which have status \textsf E):
\begin{center}\small
\begin{tabular}{lrrrrl}\toprule
class&lattices&\textsf E&\textsf R&\textsf F&other\\\midrule
$b=42$, $r=8$ on some line&375&364&8&0&\textsf O: 3\\
$b=46$, $r\le7$&143\,699&115\,649&23\,162&3\,221&\textsf O$_1$: 1\,431, \textsf O$_2$: 231, \textsf O$_3$: 4, \textsf O$_8$: 1\\
$b=47$, $r\le7$&155\,619&115\,143&34\,379&1\,839&\textsf M: 4, \textsf O$_1$: 3\,962, \textsf O$_2$: 279,\\
&&&&&\textsf O$_3$: 8, \textsf O$_4$: 1, \textsf O$_{0,0}$: 4\\\bottomrule
\end{tabular}
\end{center}

\paragraph{The search.}
The nodes, the children and the statuses are as in Section~\ref{sec:search}, with $80$ replaced by
$79+\nu$, $39-\delta$ by $38+\nu-\delta$, $\beta$ by $\beta_\nu$, and $v_0$ by
$\max\{v_0,6+\nu-\delta\}$; in (C2) one adds $v(L)\ge6+\nu-\delta$.  The roots are the lattices of
$\mathcal K^\nu_\delta$ of status \textsf F, \textsf M or \textsf O.

\begin{compres}\label{cr:nu}
For $\nu=2$ and $\nu=3$ and every root, no node on $17$ lines is open: every child on $17$ lines of an
open node on $16$ lines has status \textsf E.  (For $\nu=2$ there is no open node on $16$ lines either.)
The numbers of open nodes are
\begin{center}\small\upshape
\begin{tabular}{ccrrrrrrc}\toprule
$\nu$&$b(M)$&$|\mathcal K^\nu_\delta|$&roots&open on $13$ lines&on $14$&on $15$&on $16$&on $17$\\\midrule
2&41&767&8&28&113&38&0&0\\
2&42&3\,650&58&137&791&181&0&0\\
2&43&14\,226&367&789&325&20&0&0\\
2&44&40\,142&1\,059&771&108&0&0&0\\
2&45&88\,669&2\,769&675&51&7&0&0\\
2&46&143\,699&4\,888&521&0&0&0&0\\
\midrule
3&42&4\,025&61&1\,369&5\,211&500&7&0\\
3&43&14\,226&367&789&1\,857&79&0&0\\
3&44&40\,142&1\,059&771&108&0&0&0\\
3&45&88\,669&2\,769&675&51&7&0&0\\
3&46&143\,699&4\,888&521&9&0&0&0\\
3&47&155\,619&6\,097&219&0&0&0&0\\
\bottomrule
\end{tabular}
\end{center}

About $3.6$ million children were examined for $\nu=2$ and $11.6$ million for $\nu=3$, from $9\,149$ and
$15\,241$ roots.
\end{compres}

\begin{proof}[Proof of Theorem~\ref{thm:nu}]
Let $\Lambda\in\mathcal W_\nu$ have a critical realization, and choose $M$ and $L_1,\dots,L_5$ as above.
If $M$ has status \textsf R we are done; otherwise $M$ is a root.  Every $N_j$ is realizable, and
$(N_{j+1},M)$ is a child of $(N_j,M)$ for $j=0,\dots,4$.  If none of $N_1,\dots,N_4$ has status \textsf R,
then $(N_1,M),\dots,(N_4,M)$ are open nodes and $(\Lambda,M)=(N_5,M)$ is a child of $(N_4,M)$ on $17$
lines; it is realizable, whereas by Computational result~\ref{cr:nu} every such child has status
\textsf E, a contradiction.  Hence some $N_j$ with $j\le4$ has status \textsf R, and
Proposition~\ref{prop:rigid} applies.
\end{proof}

\paragraph{Checks.}
(i) \emph{Regression.}  For $\nu=1$ the program used here reproduces the numbers of
Computational result~\ref{cr:search} in all five classes.

(ii) \emph{Independent implementation.}  The second program of Section~\ref{sec:checks}, with the changes
of constants depending on $\nu$ (but without the bound $v(L)\ge6+\nu-\delta$, so that it examines a
slightly larger tree), was run on all roots of the class $(\nu,b(M))=(2,41)$ and on random samples of ten
($\nu=2$) or eight ($\nu=3$) roots of every other class: $106$ roots, $29$ of them of status \textsf O.
For every root the numbers of open nodes and of nodes of status \textsf R agree at every level.

(iii) \emph{Finite fields.}  The search of Section~\ref{sec:checks} over $\mathbb F_{101}$,
$\mathbb F_{103}$, $\mathbb F_{107}$ was run on the lattices with $b=42$ and a line with $8$ points and on
one lattice in forty of the classes $b=46$ and $b=47$ ($5\,619$ lattices with a one-parameter
construction, $4\,634$ of status \textsf E).  No lattice of status \textsf E has a realization over any
of the three fields, the lattices of status \textsf F, \textsf M, \textsf O have one over each of them, and
those of status \textsf R over at least one, except $113$ lattices whose Galois orbit is defined by an
irreducible polynomial of degree $3$ to $6$ without roots modulo the three primes.

(iv) \emph{Statuses of the new classes.}  For every lattice of the three classes of the table above
which has a construction with at most one parameter ($357$ lattices with $b=42$, $140\,058$ with $b=46$
and $137\,878$ with $b=47$), the status given by the fast engine was recomputed with the general engine
after a random relabelling of the lines: the statuses agree in all $278\,293$ cases.  The general-engine
statuses of the $21\,400$ lattices of these classes whose construction needs more than one parameter (and
of the $33$ such lattices of $\mathcal K_1$, already checked in Section~\ref{sec:checks}) were
recomputed after two random relabellings, with identical statuses and dimensions.

(v) \emph{Positive control.}  The extended Hesse arrangement and $\cA(G(6,1,3))$ have $21$ lines and are
free with exponents $(1,7,13)$ (so $D_0$ has exponents $(7,13)$), the coexponents of the reflection groups $G_{26}$ and $G(6,1,3)$
\cite[Chapter~6 and Appendix~C]{OT92}, so every arrangement obtained from them by deleting four lines has
$\mdr\le7$, and it is critical as soon as $D_0$ vanishes in degree $6$.  These deletions give exactly twelve
isomorphism classes of critical lattices of defect $2$ or $3$ with multiplicities at most $6$ (nine of
defect $2$, three of defect $3$; \texttt{lib/known17.py}, with $D_0$ in degree $6$ computed modulo a prime,
which can only increase its dimension).  For each of them the conditions of the window hold, no minimal
core tested has status \textsf E (all minimal cores, or twelve chosen at random when there are more), and
along up to three admissible orders of the other lines the path from the core reaches a node of status
\textsf R on $13$ or $14$ lines.

\section{Eighteen lines}\label{sec:18}

Proposition~\ref{red:18} supplies the reduction for every deletion of
a hypothetical counterexample. We now choose one deletion whose point
multiplicities are at most six.

\begin{proposition}\label{prop:18}
In the situation of Proposition~\ref{red:18} there is a line $H$ with $r_H\in\{7,8\}$ such that all points of
$\Lambda\setminus H$ have multiplicity at most $6$.  Thus $\Lambda\setminus H$ belongs to
$\mathcal W_2\cup\mathcal W_3$ and has a critical and a non-critical realization.  If $\Lambda$ has no
point of multiplicity $7$, one can take $r_H=8$, so that $\Lambda\setminus H\in\mathcal W_2$.
\end{proposition}

\begin{proof}
For a point $P$ of multiplicity $m$, a vertex is \emph{off the pencil} of $P$ if it lies on no line through
$P$; let $U_P$ be the sum of $m_q-1$ over these vertices.  Each of the $18-m$ lines not through $P$ meets
the lines through $P$ in $m$ distinct points, so $89=(m-1)+m(18-m)+U_P$, and such a line contains at most
$8-m$ vertices off the pencil.

\emph{No point has multiplicity $8$:} then $U_P=89-7-80=2$, but no line can contain a vertex off the pencil.
(Multiplicities larger than $8$ are excluded by $r_H\le8$ on a line outside the pencil; such a line exists because the arrangement is not a pencil.)

\emph{Two points $P,P'$ of multiplicity $7$ on a line $H_0$ of $\Lambda$ do not exist.}  Here $U_P=6$.
Let $V_1,\dots,V_6$ and $W_1,\dots,W_6$ be the other lines through $P$ and $P'$, and $C_1,\dots,C_5$ the
remaining lines.  If $C_i$ contains $k_i$ of the points $V_a\cap W_b$, it meets $H_0$, the $V_a$ and the
$W_b$ in $13-k_i$ points, so $5\le k_i\le6$.  A vertex off the pencil of $P$ lies on exactly one line
$W_b$ and on $w\ge1$ lines $C_i$, and contributes $w$, or it lies only on $u\ge2$ lines $C_i$ and
contributes $u-1$.  The contributions of the first kind add up to $\sum_i(6-k_i)$, the number of $i$ with
$k_i=5$.  A line with $k_i=5$ has already $8$ points, so the vertices of the second kind involve only
lines with $k_i=6$, each in at most one of them, and contribute at most the number of these lines.  Hence
$U_P\le5$, a contradiction.

\emph{Two points $P,P'$ of multiplicity $7$ not on a line of $\Lambda$.}  Let $V_1,\dots,V_7$ and
$W_1,\dots,W_7$ be the lines through $P$ and $P'$, and $C_1,\dots,C_4$ the others.  The point $P'$ is off
the pencil of $P$ and contributes $6=U_P$; so there is no other vertex off the pencil of $P$.  Hence every
point $C_i\cap W_b$ lies on some $V_a$. Applying the same identity with
$P$ and $P'$ interchanged shows that $C_i\cap C_j$ lies on both a $V_a$
and a $W_b$. Take $PP'$
as the line at infinity: the $V_a$ are the lines $x=\xi$, $\xi\in X$, the $W_b$ the lines $y=\eta$,
$\eta\in Y$, with $|X|=|Y|=7$, and $C_i$ is the graph of an affine map $g_i$ with $g_i(X)=Y$.  For
$i\ne j$ the affine map $g_i^{-1}g_j$ permutes $X$ and has exactly one fixed point, which belongs to $X$:
the point $C_i\cap C_j$ lies on some line $x=\xi$ and some line $y=\eta$, which meet the line at infinity
only in $P$ and $P'$, so $C_i\cap C_j=(\xi,\eta)$ is an affine point and $g_i(\xi)=g_j(\xi)$; there is no
other fixed point since $C_i\ne C_j$.
These maps generate a finite group of affine maps of $\CC$ (it acts faithfully on $X$), which fixes the
barycentre $c$ of $X$.  So all of them are rotations about $c\in X$, by elements of a group $\mu_k$ which
acts freely on the six points of $X\setminus\{c\}$ and contains the four distinct elements
$g_1^{-1}g_j$, $j=1,\dots,4$ (the first is the identity); hence $k=6$ and $X=\{c\}\cup(c+t\mu_6)$.  After affine changes of $x$ and of $y$,
$X=Y=\{0\}\cup\mu_6$ and $C_i=\{y=\zeta_ix\}$ with $\zeta_i\in\mu_6$.  So every realization of $\Lambda$ is
projectively equivalent to
\[
 \cA_I:\quad x\,(x^6-z^6)\;y\,(y^6-z^6)\prod_{\zeta\in I}(y-\zeta x)=0,\qquad I\subset\mu_6,\ |I|=4 ,
\]
The derivation $\theta=x(x^6-z^6)\partial_x+y(y^6-z^6)\partial_y$ is tangent to every line of
$\cA_I$: indeed $\theta(y-\zeta x)=(y^7-(\zeta x)^7)-z^6(y-\zeta x)$ for $\zeta^6=1$.  It has degree $7$ and
is not a multiple of the Euler derivation, so $\mdr(\cA_I)\le7$, whereas a free arrangement of $18$ lines
with $b=89$ has exponents $(8,9)$.  So no realization of $\Lambda$ is free, which contradicts the existence
of $\cA$.

\emph{Choice of $H$.}  By the above, $\Lambda$ has at most one point of multiplicity $7$ and no point of
higher multiplicity.  Suppose first that all multiplicities are at most $6$.  From $\sum_P(m_P-1)=89$ and
$\sum_P\binom{m_P}2=153$ one gets $t_2=3t_4+8t_5+15t_6-39\ge0$ and
$\sum_H(r_H-7)=t_4+3t_5+6t_6-12\ge\tfrac13\cdot39-12=1$.  So some line $H$ has $8$ points, and
$\Lambda\setminus H\in\mathcal W_2$.

Now let $P$ be the point of multiplicity $7$; then $U_P=6$, and the off-pencil vertices $q$ lie on
pairwise disjoint sets of the $11$ lines not through $P$.  For a line $H_i$ through $P$ let
$g_i=r_{H_i}-1$ be the number of its vertices other than $P$, and $e_i$ the number of pairs of lines not
through $P$ meeting on $H_i$; then $\sum_ie_i+\sum_q\binom{m_q}2=\binom{11}2=55$.  If $g_i\le5$, then
$e_i\ge7$: distributing eleven lines among at most five points minimizes the number of pairs at $11=3+2+2+2+2$, giving $\binom32+4\binom22=7$.  If this held for all seven lines, we would get
$\sum_q\binom{m_q}2\le6=\sum_q(m_q-1)$, so all $m_q=2$, and six disjoint pairs among $11$ lines, which is
absurd.  So some line $H$ through $P$ has $r_H\ge7$; the point $P$ has multiplicity $6$ in
$\Lambda\setminus H$, which therefore belongs to $\mathcal W_2$ if $r_H=8$ and to $\mathcal W_3$ if
$r_H=7$.
\end{proof}

\begin{corollary}\label{cor:terao18}
Terao's conjecture holds for arrangements of $18$ lines in $\PP^2_\CC$.
\end{corollary}
\begin{proof}
Otherwise Proposition~\ref{prop:18} gives a lattice in $\mathcal W_2\cup\mathcal W_3$ with a critical and
a non-critical realization, contradicting Theorem~\ref{thm:nu}.
\end{proof}
The eighteen-line argument does not use Corollary~\ref{cor:terao17}.  Defect $3$ is needed only for the lattices
$\Lambda$ with exactly one point of multiplicity $7$ through which no line has $8$ points.

\begin{proof}[Proof of Theorem~\ref{thm:terao1718}]
Combine Corollaries~\ref{cor:terao17} and~\ref{cor:terao18}.
\end{proof}
\section{Reliability of the computation}\label{sec:reliability}

\paragraph{The general engine.}  The general engine accepts the output of \textsc{Singular} only as a
complete, versioned transcript: the number of minimal primes, then every minimal prime with its ideal and
either its dimension or the forbidden factor that it contains, then a final summary; the process must
terminate normally, with no error output and within the time limit.  Anything else (an error, a malformed
or incomplete output, a nonzero exit status, a timeout, or a transcript with no minimal prime at all,
which \texttt{minAssGTZ} never returns), and also a construction with more than eight parameters, gives the
status \emph{undecided}.  An undecided lattice is treated like a lattice of status \textsf O: it is never
discarded, and the search continues below it (the programs of the searches also count such nodes
separately).  In particular a failed calculation
can never be read as an empty realization space, nor as a single Galois orbit, and the conclusions of the
searches do not depend on whether undecided nodes occurred.  No lattice of the census has status
undecided.

\paragraph{The searches.}  The census, the statuses of the cores, all searches and all checks were computed
with the programs in the ancillary material.  This material also contains the census, the statuses and, for
every root, the record of the search (the numbers of children, of children sent to each engine and of nodes of
each status at every level) together with the list of all open nodes.  It does not contain the transcripts
of the individual \textsc{Singular} calls: when the environment variable
\texttt{GST17\_SINGULAR\_AUDIT\_DIR} is set, the general engine keeps the input, the output and a receipt
(exit status, running time, hashes) of every call, so that the transcripts of any part of the computation
can be produced by running it again.  Canonical forms are computed with \texttt{nauty}; a failure of
\texttt{nauty} (a nonzero exit status, or an output other than one canonical string per input) stops the
program, so it is never turned into a canonical form.  The program \texttt{check/validate\_records.py}
checks that the census lists consist of canonical forms without repetitions, that the statuses are aligned
with them, that the records of the general engine exist for exactly the lattices which the searches need,
all with a certified status, that the records of all classes of the searches are present (compared with
an explicit list), and that the records of every search contain each root of its class exactly once.  The numbers of roots and of
children examined are
\begin{center}\small
\begin{tabular}{crr}\toprule
Defect & Roots & Children examined\\\midrule
$1$ & $4\,261$ & $1\,604\,862$\\
$2$ & $9\,149$ & $3\,614\,809$\\
$3$ & $15\,241$ & $11\,556\,224$\\
$4$ & $15\,471$ & $18\,151\,840$\\\bottomrule
\end{tabular}
\end{center}
These counts refer to rooted search paths and include repetitions across roots; they are not counts of
non-isomorphic lattices.

\paragraph{Scope of the verification.}  For each computational input of the proof we list how it has been
checked.
\begin{itemize}[leftmargin=*]
\item \emph{Census} (Section~\ref{sec:census}).  Generated with the programs of \cite{AS}; two sectors are
checked by hand.  In the searches with the reduction of Section~\ref{sec:nu4} ($\delta=2,3$), every
$12$-line sublattice of weight $41+\delta$ that was examined was found in the census list.  The
completeness of the census rests on the generator.
\item \emph{Fast engine} (Section~\ref{sec:real}).  Its statuses were compared with those of the general
engine, after a random relabelling of the lines, for every lattice with a one-parameter construction of
$\mathcal K_0,\dots,\mathcal K_4$ (Section~\ref{sec:checks}), of the classes $b=46$ and $b=47$ with at
most $7$ points on every line (Section~\ref{sec:nu}), and of the classes $b=42$ and $b=43$ with a line
with $8$ points (Sections~\ref{sec:nu} and~\ref{sec:nu4}); it was also compared with a search over three
finite fields on large samples.  No disagreement.
\item \emph{General engine.}  A failed calculation gives the status undecided, which is never read as
\textsf E or \textsf R (see above); all its statuses on the cores were recomputed after
two further random relabellings (Sections~\ref{sec:checks}, \ref{sec:nu} and~\ref{sec:nu4}).
\item \emph{Extension searches.}  A second implementation of the search, written from its description,
was run on all roots for defect one, on $106$ roots for defects two and three and on two roots for defect
four (Sections~\ref{sec:checks}, \ref{sec:nu} and~\ref{sec:nu4}).  The two implementations share the
census and the two engines.  The incremental statuses were compared with statuses computed from scratch
(Section~\ref{sec:checks}).
\item \emph{Exact computations} (Lemma~\ref{lem:LF}, Section~\ref{sec:18}, Appendices~\ref{app:A}
and~\ref{app:B}).  Deterministic symbolic or modular computations, listed in Section~\ref{sec:verification}.
\item \emph{Inputs not verified here}: the census generator beyond the two sectors checked by hand, and
the results of \cite{GST16}, \cite{AS}, \cite{Zie89} and \cite{Abe} that are quoted.
\end{itemize}
The two routes have different dependencies.  The generic splitting type (Theorem~\ref{thm:main-gst}) uses
\cite[Theorems~A, B and~3.5]{GST16} and Theorem~\ref{thm:c7}.  The reductions for Terao's conjecture are
proved in Section~\ref{sec:reductions}; besides elementary arguments they use only \cite{GST16} (which
relies on \cite{AS}; both are arXiv preprints, themselves computer-assisted) and, for seventeen lines,
Theorem~\ref{thm:sept}.

\section{Defect four}\label{sec:nu4}

We use the setting of Section~\ref{sec:nu} with $\nu=4$: Lemmas~\ref{lem:gen}, \ref{lem:del} and
Proposition~\ref{prop:need}, the bounds $\beta_\nu$, the properties (a)--(e) of a minimal core, the
reconstruction and the search with the conditions of Section~\ref{sec:search}.

\begin{theorem}\label{thm:nu4}
Let $\Lambda\in\mathcal W_4$.  If $\Lambda$ has a critical realization, then $\Lambda$ has one Galois
orbit of realizations.  In particular $(C_7)$ holds for every $7$-candidate of defect $4$ with
multiplicities at most $6$.
\end{theorem}
In the notation of Sections~\ref{sec:cores} and~\ref{sec:nu}, for $\nu=4$ $t_2=3W-23-t_5-3t_6\ge0$ gives $W\ge8$ and
$\sum_{|S|=5}\varphi(S)=715W+220t_5+715t_6-5278>0$, so some $S$ has $\varphi(S)\ge1$ and
$b(\Lambda\setminus S)=48-\varphi(S)\le47$; with (e), $2\le\delta\le6$.  The bounds are
$\beta_4(10)=27$, $\beta_4(11)=35$, $\beta_4(12)=43$, $\beta_4(13)=51$, $\beta_4(14)=59$.  The classes of
cores are: for $\delta=2$ ($\rho=8$) all $18\,569$ lattices with $b=43$ of the census of
Section~\ref{sec:census}, and for $\delta=3,\dots,6$ ($\rho=7$) the classes $\mathcal K_3$, $\mathcal K_4$,
$\mathcal K^2_5$, $\mathcal K^3_6$ of Sections~\ref{sec:census} and~\ref{sec:nu}.  The $4\,343$ lattices with $b=43$ and a line with $8$ points have statuses
\textsf E: $3\,614$, \textsf R: $438$, \textsf F: $214$, \textsf O: $77$.

For $\delta=2$ and $\delta=3$ the search was run with the following additional reduction.  Number the
lattices of the class by their positions in the census list of lattices with $b=41+\delta$ and multiplicities
at most $7$ (\texttt{census/m7\_X23\_lo3.g6} and \texttt{census/m7\_X22\_lo4.g6}; any total order would
do).  Among the minimal cores of $\Lambda$ let $M$ be one with the largest
number, and let $N_1,\dots,N_4$ be the nodes of the path of $\Lambda$ from $M$.  Every $12$-line
sublattice of weight $41+\delta$ of a node $N_j$ is a minimal core of $\Lambda$; so it belongs to
$\mathcal K^\nu_\delta$ (property~(c) holds for every minimal core), in particular to the census list, and
its number is at most that of $M$.  Hence a node with such a sublattice of larger number lies on no such
path, and it is discarded.  A node with such a sublattice which is not in the census list lies on no such
path either; for it the program stops with an error, so that an incomplete census list could not go
unnoticed, and this does not occur.

\begin{compres}\label{cr:nu4}
For $\nu=4$ exactly one node on $17$ lines is open.  It occurs below the root with number $13\,396$ of
the class $b(M)=43$, and its lattice $\Lambda_F$ (Lemma~\ref{lem:LF}) has status \textsf F.  The numbers
of open nodes are (the classes marked $^*$ were run with the reduction above):
\begin{center}\small\upshape
\begin{tabular}{ccrrrrrrc}\toprule
$\nu$&$b(M)$&$|\mathcal K^\nu_\delta|$&roots&open on $13$ lines&on $14$&on $15$&on $16$&on $17$\\\midrule
4&43$^*$&18\,569&658&8\,323&11\,663&514&9&1\\
4&44$^*$&40\,142&1\,059&123&44&0&0&0\\
4&45&88\,669&2\,769&675&51&16&0&0\\
4&46&143\,699&4\,888&521&9&0&0&0\\
4&47&155\,619&6\,097&219&0&0&0&0\\
\bottomrule
\end{tabular}
\end{center}

About $18.2$ million children were examined, from $15\,471$ roots.  Six lattices on $17$ lines have
status \textsf R (below the roots $3\,281$ and $13\,396$ of the class $b(M)=43$); each of them has one
Galois orbit of realizations, defined over a quadratic field.
\end{compres}

\begin{lemma}\label{lem:LF}
Let $\Lambda_F$ be the lattice on the lines $H_0,\dots,H_{16}$ whose points of multiplicity at least $3$
are
\[
\begin{gathered}
\{0,1,2,7,8\},\ \{0,3,4,6,9\},\ \{5,6,7,12,13\},\ \{5,8,9,14,15\},\ \{2,3,12,14,16\},\\
\{1,4,13,15,16\},\ \{10,11,13,14\},\ \{0,10,12,15\},\ \{0,5,11,16\},\\
\{1,3,10\},\ \{2,4,10\},\ \{6,8,10\},\ \{7,9,10\},\\
\{2,6,11\},\ \{3,7,11\},\ \{4,8,11\},\ \{1,9,11\}
\end{gathered}
\]
(indices of the lines), with $34$ double points; so $b=83$, $H_0$ has $6$ points and the other lines
have $8$.  Then the realizations of $\Lambda_F$ form a one-parameter family, all of them have
$\mdr\ge8$, and the generic one has $\mdr=9$.  In particular no realization of $\Lambda_F$ is critical.
\end{lemma}

\begin{proof}
\emph{Realizations.}  The lines $H_{11},H_6,H_7,H_8$ are in general position in $\Lambda_F$, so every
realization is projectively equivalent to one with $H_{11}=\{x=0\}$, $H_6=\{y=0\}$, $H_7=\{z=0\}$,
$H_8=\{x+y+z=0\}$.  The other lines are then determined one by one, each by two points already
constructed, except one line, which lies in a pencil through a point already constructed; its
parameter $t$ gives the family
\begin{gather*}
H_0\colon x+y+tz,\quad H_1\colon(2t-1)x+(2t-1)y+t^2z,\quad H_2\colon x+y,\quad H_3\colon x+tz,\\
H_4\colon x+ty+tz,\quad H_5\colon y+tz,\quad H_9\colon tx+(2t-1)y+t^2z,\quad
H_{10}\colon tx+(2t-1)y+tz,\\
H_{12}\colon y-tz,\quad H_{13}\colon(2t-1)y+tz,\quad H_{14}\colon 2tx+(2t-1)y+tz,\\
H_{15}\colon 2tx+(3t-1)y+(t^2+t)z,\quad H_{16}\colon 2x+y+tz.
\end{gather*}
The excluded member of the pencil is a line already constructed.  All concurrences of $\Lambda_F$
hold identically in $t$, and a further concurrence or coincidence occurs exactly for $t$ in the set $B$
of roots of $t(t-1)(t+1)(t-2)(2t-1)(3t-1)(3t-2)(2t^2-2t+1)$.  So the realization space is
$\CC\setminus B$, an irreducible curve.

\emph{A line with eight points.}  Let $\cA$ be a realization and $H\in\cA$ a line with $r_H=8$; we show
$\mdr(\cA)\ge7$.  Let $\theta\in D(\cA)_e\setminus SE$ with $e=\mdr(\cA)$ and suppose $e\le6$.  If the
projective restriction of $\theta$ to $H$ is nonzero, Lemma~\ref{red:restriction}(ii) gives $8\le e+1$, a
contradiction.  Otherwise $\theta=\alpha_H\theta_0+hE$ with $\theta_0\in D(\cA_H)_{e-1}\setminus SE$; since
$|\cA_H|-r_{\cA_H}(H)=16-8>e-1$, Lemma~\ref{red:restriction}(iii) shows that $\theta_0$ is tangent to $H$,
so $\theta_0\in D(\cA)_{e-1}\setminus SE$, contrary to the minimality of $e$.

\emph{Generic $\mdr$.}  At $t=3\notin B$ an exact computation over $\QQ$ gives $D_0$ of dimension
$0,0,3$ in degrees $7,8,9$ (\texttt{check/lambdaF/witness.py}, which also checks the realization space
above); since dimensions can only grow under specialization, the
generic realization has $\mdr\ge9$; on the other hand $h^0(\T\cA(9))\ge\chi(\T\cA(9))=3$ for every
realization of $\Lambda_F$ (here $h^2(\T\cA(9))=h^0(\T\cA(4))=0$, since $\mdr\ge7$ by the preceding
paragraph), so $\mdr\le9$.  (Over $\mathbb F_{10007}$ and $\mathbb F_{10009}$ all $10\,000$
values of $t$ which give the lattice $\Lambda_F$ have $D_0$ zero in degrees $7$ and $8$.)

\emph{No critical member.}  Let $\cA$ be a critical realization and $H\in\cA$.  Then
$D(\cA)=S\,E\oplus D_H(\cA)$, where $E$ is the Euler derivation, with $D_H(\cA)=\{\theta:\theta(\alpha_H)=0\}$, and
$D_H(\cA)\cong D_0(\cA)$ as graded modules; so $D_H(\cA)_6=0\ne D_H(\cA)_7$.  The Ziegler restriction
$D_H(\cA)\to D(\cA^H,m^H)$ has kernel $\alpha_HD_H(\cA)$ \cite{Zie89}, so the multiarrangement
$(\cA^H,m^H)$ on $H$, of total multiplicity $16$, has a nonzero derivation of degree $7$.  Its
derivations of degree $7$ are the kernel of a $16\times16$ matrix (pairs of binary forms of degree $7$,
$16$ conditions), which we write over $\ZZ[t]$ after choosing coordinates on $H$ by two fixed linear
forms; this matrix describes the multirestriction exactly for $t$ outside the zeros of an explicit
polynomial $e_H(t)$ (the two forms must be independent on $H$, and the second one must not vanish at a
point of $\cA^H$).  Let $\Delta_H(t)$ be its determinant.  Then a critical realization is a zero of
$\Delta_He_H$ for every $H$.  For $H\in\{H_1,H_2,H_3,H_4,H_6,H_7,H_8,H_9\}$ the polynomial $\Delta_H$ is
nonzero, and the irreducible factors of the greatest common divisor of the
eight products $\Delta_He_H$ are $t$, $t-1$, $t+1$, $t-2$, $2t-1$, $3t-1$, $3t-2$, $2t^2-2t+1$, all of
which vanish only on $B$.  The same holds for three choices of the two linear forms.  So $\Lambda_F$ has
no critical realization; since it has lines with $8$ points, $\mdr\ge7$ by the first paragraph, hence
$\mdr\ge8$ for all its realizations.
\end{proof}

\begin{proof}[Proof of Theorem~\ref{thm:nu4}]
As for Theorem~\ref{thm:nu}, with $M$ the minimal core of largest number (so that no node $N_j$ of the
path of $\Lambda$ is discarded by the reduction): if none of $N_1,\dots,N_4$ has status \textsf R, then
$(\Lambda,M)$ is a realizable child on $17$ lines of an open node, so by Computational
result~\ref{cr:nu4} it has status \textsf R or $\Lambda=\Lambda_F$, and the latter is excluded by
Lemma~\ref{lem:LF}.
\end{proof}

\paragraph{Checks for defect four.}
(i) For $(\nu,b(M))=(3,42)$ the search with the reduction gives the same conclusion as the search of
Computational result~\ref{cr:nu} (no open node on $17$ lines), and each of its $3\,141$ open nodes
is an open node of the search without the reduction.  For the twelve critical lattices of check (v) of
Section~\ref{sec:nu}, the minimal core with the largest number is followed by the search with the
reduction.

(ii) The second program of Section~\ref{sec:checks}, which does not use the reduction, was run on the
roots $3\,281$ and $13\,396$ of the class $b(M)=43$.  For $13\,396$ the numbers of open nodes and of
nodes of status \textsf R agree at every level, and its only open node on $17$ lines is isomorphic to
$\Lambda_F$.  For $3\,281$ it finds $50$ open nodes instead of $27$: the $27$ open nodes of the search
with the reduction are among them, and each of the other $23$ contains a $12$-line sublattice of weight
$43$ with a larger number than $3\,281$.  Its only open node on $17$ lines there is again $\Lambda_F$,
whose $12$-line sublattices of weight $43$ have the numbers $4,6,30,70,88,155,225,2\,917,3\,281,
10\,416,13\,396$; with the reduction, $\Lambda_F$ is reached only from $13\,396$.

(iii) The six lattices of status \textsf R on $17$ lines are not critical: for each of them the
reduction of its realization modulo three primes has $D_0$ of dimension $0,0,3$ in degrees $7,8,9$, so
$\mdr=9$ (as in the proof of Lemma~\ref{lem:LF}).  No critical lattice of defect $4$ with multiplicities at most $6$ is known, so there is no
positive control at this defect.

(iv) \emph{Statuses of the additional cores.}  The cores of the class $b(M)=43$ with a line with $8$
points occur only at defect four.  For the $4\,166$ of them with a construction with at most one
parameter, the status given by the fast engine was recomputed with the general engine after a random
relabelling of the lines, and the general-engine statuses of all $323$ lattices with $b=43$, multiplicities
at most $6$ and a construction with more parameters ($177$ of them with a line with $8$ points) were
recomputed after two random relabellings: no status changed.  (The other classes of cores are those of
Section~\ref{sec:nu}.)

\section{Candidates with a point of multiplicity seven}\label{sec:sept}

Throughout this section $\Lambda$ is a $7$-candidate of defect $\nu\ge1$ with a point $P$ of
multiplicity $7$.  The seven lines
through $P$ are the \emph{pencil lines}, the other ten the \emph{transversals}; a vertex is \emph{off the
pencil} if it lies on no pencil line.

\begin{lemma}\label{lem:septcount}
\textup{(a)} Every transversal contains at most one vertex off the pencil; so these vertices $q$ lie on
pairwise disjoint sets of transversals, $\sum_qm_q\le10$, and $\sum_q(m_q-1)=\nu+3$.

\textup{(b)} A $7$-candidate has multiplicities at most $8$, and one with a point of multiplicity $8$
has defect $0$; so here all multiplicities are at most $7$.

\textup{(c)} $\nu\le5$.
\end{lemma}
\begin{proof}
(a) A transversal meets the pencil lines in $7$ points and has at most $8$ points.  Counting $b$ on the
pencil lines gives $79+\nu=6+7\cdot10+\sum_q(m_q-1)$.
(b) A line not through a point of multiplicity $m$ has at least $m$ points, so $m\le8$.  If $m=8$, the
nine lines not through the point meet the lines through it in $8$ points each, so there is no vertex off
these lines, and $b=7+8\cdot9=79$.
(c) By (b) $m_q\le7$; with one vertex off the pencil $\sum(m_q-1)\le6$, and with at least two
$\sum(m_q-1)\le10-2$.
\end{proof}

\begin{lemma}\label{lem:pencil}
Let $\cA$ be a realization of $\Lambda$, $p\in\CC^3$ a vector representing $P$, $D_p=\sum p_i\partial_i$
and $G$ the product of the equations of the transversals.  Then $\theta_P=G\,D_p$ is tangent to all
lines, and its class modulo the Euler derivation is a section of $\T\cA(10)$ with a finite zero scheme $Z_P$ of length $\nu+3$, supported on the vertices off the pencil;
at such a vertex $q$, $Z_P$ is the subscheme of length $m_q-1$ of the line $qP$ supported at $q$.
Moreover $D_0(\cA)_7\cong H^0(\cI_{Z_P}(1))$; in particular $\mdr(\cA)\ge7$.
\end{lemma}
\begin{proof}
In an affine chart the field of $\theta_P$ is $X=G\,X_0$, where $X_0$ is the
field of $D_p$: it vanishes only at $P$, where $X_0=-E_P$, and it is tangent to the lines through $P$.
We use the local frames of $T_{\PP^2}(-\log\cA)$: $E_q=x\partial_x+y\partial_y$ and
$\eta_q=\partial_yQ_q\,\partial_x-\partial_xQ_q\,\partial_y$ at a vertex $q$ with local equation $Q_q$ of
degree $m$, for which $\det(E_q,\eta_q)=-mQ_q$; and $\partial_x$, $y\partial_y$ at a point on exactly one
line $\{y=0\}$.

At $P$, $X=-G\,E_P$ with $G(P)\ne0$.  At a vertex $q\ne P$ on a pencil line $\{y=0\}$, $Q_q=y\Pi$,
$G=\Pi\cdot$unit and $X_0$ is a nonvanishing field tangent to $\{y=0\}$; so $\det(E_q,X_0)=y\cdot$unit,
and the coefficient $\det(E_q,X)/\det(E_q,\eta_q)$ of $\eta_q$ in $X$ is a unit.  At a point on exactly one line the coefficient of $y\partial_y$
(on a transversal, to which $X_0$ is transverse) or of $\partial_x$ (on a pencil line) is a unit; off the
arrangement $X$ does not vanish.  So $X$ vanishes, as a logarithmic field, only at the vertices $q$ off
the pencil.  There $G=Q_q\cdot$unit, and $X=aE_q+b\eta_q$ with
\[
 a=\frac{\det(X,\eta_q)}{\det(E_q,\eta_q)}=\text{unit}\cdot X_0(Q_q),\qquad
 b=\frac{\det(E_q,X)}{\det(E_q,\eta_q)}=\text{unit}\cdot\det(E_q,X_0).
\]
The function $\det(E_q,X_0)$ vanishes exactly where $X_0$ is radial from $q$, i.e.\ it is a local equation
of the line $qP$; and the restriction of $X_0(Q_q)$ to this line vanishes to order exactly $m_q-1$ at
$q$, because the line $qP$ is not a line of $\cA$ through $q$.  This gives $Z_P$ near $q$.  Its length
is $\sum_q(m_q-1)=\nu+3=c_2(\T\cA(10))$.  Finally $\theta_P$ gives
$0\to\cO\to\T\cA(10)\to\cI_{Z_P}(4)\to0$; twisting by $\cO(-3)$ and taking global sections gives
$D_0(\cA)_7\cong H^0(\cI_{Z_P}(1))$, since $H^0(\cO(-3))=H^1(\cO(-3))=0$.  Twisting by $\cO(-4)$ gives in
the same way $D_0(\cA)_6\cong H^0(\cI_{Z_P})=0$, because $Z_P\ne\varnothing$ ($\nu+3>0$); so $D_0(\cA)$
vanishes in degrees at most $6$ (a nonzero element of lower degree times a nonzero form would lie in
degree $6$), that is, $\mdr(\cA)\ge7$.
\end{proof}

\begin{corollary}\label{cor:septcrit}
A realization $\cA$ of $\Lambda$ is critical if and only if its vertices off the pencil lie on a line
$\ell$ which, if one of them has multiplicity at least $3$, passes through $P$.  The line $\ell$ is not
a line of $\cA$.  If there is only one vertex off the pencil, every realization is critical.
\end{corollary}
\begin{proof}
$\cA$ is critical iff $Z_P$ lies on a line $\ell$.  At a vertex $q$ of multiplicity $\ge3$ the scheme
$Z_P$ has length $\ge2$ and lies on the line $qP$, so $\ell=qP$.  A pencil line contains no vertex off
the pencil and a transversal at most one, while there are at least two if $\ell\not\ni P$, and
$\ell=qP$ is not a line of $\cA$ otherwise.
\end{proof}

So a critical realization is of one of two kinds.
\begin{itemize}[leftmargin=*]
\item[(A)] $\ell\ni P$.  Then in $\cA^+=\cA\cup\{\ell\}$ any two transversals meet on one of the $8$ lines
through $P$.
\item[(B)] $\ell\not\ni P$.  Then all vertices off the pencil are double points, so $2(\nu+3)\le10$,
i.e.\ $\nu\le2$; they lie on $\ell$, and the other pairs of transversals meet on pencil lines.
\end{itemize}

\paragraph{Transversal structures.}
A \emph{structure of type} A consists of $8$ \emph{base lines} through a point $P$ and $k$
\emph{transversals} not through $P$, any two of which meet on a base line.  In a \emph{structure of
type} B the base lines are $B_1,\dots,B_7$ through $P$ and a line $B_0\not\ni P$; any two transversals
meet on a base line, at most two of them at the same point of $B_0$, and the points $B_0\cap B_i$ lie on
no transversal (this last condition is automatic when all transversals are paired on $B_0$; see the
proof of Theorem~\ref{thm:sept} for the other case).  The combinatorial type of a structure records, for every base line, the partition of
the transversals according to their points on it (the \emph{fibres} are the parts with at least two
elements); it determines the lattice of the $8+k$ lines, and conversely.  The \emph{weight} of a base
line is the sum of $|F|-1$ over its fibres $F$.  Structures are considered up to isomorphisms of
lattices preserving the base lines, the line $B_0$ and the transversals.

A structure with $k$ transversals is obtained from one with $k-1$ by adding a transversal; the
possibilities are described by a choice, for each base line, of a fibre which the new line joins, of a
single transversal with which it forms a new fibre, or of none, such that every old transversal is met
on exactly one base line.  For $k\ge3$ the new line passes through at least two vertices unless all old
transversals lie in one fibre.

\begin{compres}\label{cr:sept}
The realizable structures were classified level by level, starting from all combinatorial types with
$k=6$ ($526$ of type A; $281$ of type B in which the six transversals form three pairs on $B_0$), with
the exact statuses of Section~\ref{sec:real} and, for the children of nodes with a chart (one Galois
orbit, several orbits, or a one-parameter family), the incremental computation of
Section~\ref{sec:search}.  In type B the transversals $7,\dots,10$
were added in the order: one not paired on $B_0$, its partner, one not paired, and a last one, paired or
not.  The numbers of structures whose status is not \textsf E are
\begin{center}\small\upshape
\begin{tabular}{lrrrrr}\toprule
$k$&6&7&8&9&10\\\midrule
type A, fibres of size $\le7$&227&428&135&37&13\\
type A, fibres of size $\le6$&227&427&134&34&10\\
type B&112&851&43&3&1\\\bottomrule
\end{tabular}
\end{center}
Of the $13$ structures of type A with ten transversals, $10$ have one Galois orbit of realizations and
$3$ a one-parameter family.  In the only structure of type B with ten transversals, five pairs of
transversals meet on $B_0$ and two of the lines $B_1,\dots,B_7$ contain no intersection point of
transversals.
\end{compres}

\begin{theorem}\label{thm:sept}
Let $\Lambda$ be a $7$-candidate with a point of multiplicity $7$.
\begin{enumerate}[label=\textup{(\alph*)}]
\item If $\Lambda$ (of defect $\nu\ge1$) has a critical realization, then $\Lambda$ is one of $21$
lattices: $10$ of defect $1$, $7$ of defect $2$ and $4$ of defect $3$.
\item For each of these lattices every realization is critical: $20$ of them have one Galois orbit of
realizations, and the last one (of defect $3$, with a one-parameter family of realizations) has a single
vertex off the pencil.
\end{enumerate}
Consequently $(C_7)$ holds for every $7$-candidate with a point of multiplicity at least $7$, and no
such candidate of defect $\nu\ge4$ has a critical realization.
\end{theorem}
\begin{proof}
Let $\cA$ be a critical realization and $\ell$ the line of Corollary~\ref{cor:septcrit}.  In case (A),
$\cA\cup\{\ell\}$ realizes a structure of type A with ten transversals; its fibres have size at most
$7$ by Lemma~\ref{lem:septcount}(b).  So it is one of the $13$ structures of
Computational result~\ref{cr:sept}, and $\Lambda$ is obtained from it by deleting a base line of weight
$\nu+3$, the other base lines having weight at least $3$ because they have at most $8$ points in
$\Lambda$.  Going through the $13$ structures and their base lines gives $21$ lattices on $17$ lines up
to isomorphism, with the defects stated.  In case (B), $\nu\le2$ and at least $8$ transversals are paired
on $\ell$.  If $\nu=2$, or if $\nu=1$ and the two unpaired transversals do not pass through the points
$\ell\cap B_i$, then $\cA\cup\{\ell\}$ realizes a structure of type B with $B_0=\ell$; taking three of the
pairs first, then a fourth, it appears in the classification, but in the only such structure two pencil
lines contain no intersection point of transversals, so each of them meets the ten transversals in ten
distinct points and, with $P$, has $11>8$ points.  If $\nu=1$ and an unpaired transversal passes through a point $\ell\cap B_i$,
we use a supplementary computation: starting from the $43$ structures of type B with eight transversals
(four pairs), the two unpaired transversals were added in all possible ways, each avoiding the points
$\ell\cap B_i$ or passing through one of them. Exactly two structures with ten transversals arise; in both the two unpaired transversals meet at a
point $\ell\cap B_i$, and one pencil line contains no intersection point of transversals, hence has $11$
points in the same way.  So case (B) does not occur.

(b) The realization spaces of the $21$ lattices were computed with the engines of
Section~\ref{sec:real}: $20$ have one Galois orbit; for the $8$ of them with at least two vertices off
the pencil the collinearity of Corollary~\ref{cor:septcrit} was also checked on the exact realization.
For the last one use Corollary~\ref{cor:septcrit}.  The last assertion follows, with
Lemma~\ref{lem:septcount}(b) for multiplicity $8$ and the remark below for defect $\nu\le0$.
\end{proof}

At defect $\nu\le0$ a realization with $\mdr\ge8$ would make $\T\cA(8)$ semistable with $c_1=0$ and
$c_2=\nu-1<0$, against Bogomolov's inequality \cite{Bog79}; so $(C_7)$ is automatic there (cf.\ \cite{GST16}).

The $21$ lattices are, by defect and by the multiplicities of the vertices off the pencil:
\begin{center}\small
\begin{tabular}{clcc}\toprule
$\nu$&vertices off the pencil&one point of multiplicity $7$&two points of multiplicity $7$\\\midrule
1&$(2,2,2,2)$&3&3\\
1&$(5)$&1&3\\
2&$(2,2,2,2,2)$&1&--\\
2&$(6)$&6&--\\
3&$(3,3,3)$&1&--\\
3&$(7)$&--&3\\\bottomrule
\end{tabular}
\end{center}
For a lattice with two points of multiplicity $7$ the row does not depend on the choice of $P$, except in
the second row: in the three lattices of the first row the vertices off the pencil are $(2,2,2,2)$ for both
points, in the last row they are $(7)$ for both, and in the three lattices of the second row with two
points of multiplicity $7$ the vertices off the pencil of the other point of multiplicity $7$ are
$(2,2,2,2)$.

\paragraph{Checks.}
(i) Since $\cA(G(6,1,3))$ is free with exponents $(1,7,13)$ (those of $D_0$ are $(7,13)$), every arrangement obtained from it by deleting
four lines has $\mdr\le7$.  These deletions give exactly twenty isomorphism classes of critical lattices
with a point of multiplicity $7$, $\nu\ge1$ and at most $8$ points on every line ($10$, $7$ and $3$ of
defect $1$, $2$, $3$), and they are twenty of the $21$ lattices of Computational result~\ref{cr:sept}
(\texttt{check/g613\_sept.py}).
(ii) The remaining lattice, of defect $3$ with vertices $(3,3,3)$ off the pencil, is not such a deletion.
It is realized over $\QQ(\omega)$, $\omega^2+\omega+1=0$, by the lines
{\small
\[
\begin{gathered}
-x-z,\ y+\omega z,\ \omega x+(-1-\omega )y+(1+\omega )z,\ (1+2\omega )x-y+(1+\omega )z,\ (-1+\omega )x-y-z,\\
-2\omega x+y-\omega z,\ x+y+(1+\omega )z,\ x,\ y,\ -\omega x+(1+\omega )y,\\
z,\ -x+\omega z,\ x+y+z,\ y+z,\\
x+(1+\omega )y+(1+\omega )z,\ \omega x-y+\omega z,\ y+(1+\omega )z
\end{gathered}

\]}
an exact computation gives minimal generators of $D_0$ in degrees $7$, $10$, $12$, and an exact
computation over $\QQ(\omega)$ confirms that these lines have the lattice in question, with the three
vertices off the pencil collinear with $P$ (\texttt{check/sept\_nu3\_lattice.py}).

\section{Proof of the main theorem}\label{sec:proof}

\begin{theorem}\label{thm:c7}
Condition $(C_7)$ holds for
\begin{enumerate}[label=\textup{(\alph*)}]
\item every $7$-candidate of defect $\nu\le4$;
\item every $7$-candidate with a point of multiplicity at least $7$, of any defect;
\item every $7$-candidate.
\end{enumerate}
\end{theorem}
\begin{proof}
A $7$-candidate has multiplicities at most $8$, and one with a point of multiplicity $8$ has defect $0$
(Lemma~\ref{lem:septcount}).  At defect $\nu\le0$ the condition is automatic (the remark after
Theorem~\ref{thm:sept}).  A candidate with a point of multiplicity $7$ is covered by
Theorem~\ref{thm:sept}; this proves~(b).  A candidate of defect $\nu\ge1$ with multiplicities at most $6$
belongs to $\mathcal W_\nu$.  For $\nu=1$, $2$, $3$, $(C_7)$ holds by Theorems~\ref{thm:main}
and~\ref{thm:nu}, and for $\nu=4$ by Theorem~\ref{thm:nu4}; this proves~(a).  For $\nu\ge5$ the candidate has no critical
realization by Theorem~\ref{v4:main}, so $(C_7)$ holds for it; this proves~(c).
\end{proof}

\begin{proof}[Proof of Theorems~\ref{thm:main-gst}, \ref{thm:defect} and~\ref{thm:c7-intro}]
Theorem~\ref{thm:defect} is Theorem~\ref{v4:main} and Theorem~\ref{thm:c7-intro} is
Theorem~\ref{thm:c7}(a),(b).  By Theorem~\ref{thm:c7}(c) and \cite[Theorems~A, B and~3.5]{GST16}, the
generic splitting type is combinatorial for arrangements of $17$ and of $18$ lines; for at most $16$ lines
this is \cite{GST16}.
\end{proof}

\begin{proof}[Proof of Corollary~\ref{cor:anyd}]
By \cite[Theorems~A and~B]{GST16}, $(C_s)$ holds for all $s\le6$, and by Theorem~\ref{thm:c7}(c) it holds
for $s=7$.  So \cite[Proposition~8.1]{GST16} applies with $S=7$: its parts (b) and (c) are (a) and (b) of
the corollary.  For the second statement of (a), $2\mdr(\cA)\le d-3$ means $\mdr(\cA)<\lfloor(d-1)/2\rfloor$,
so $k(\cA)=\mdr(\cA)$, and then $k(\cB)=k(\cA)<\lfloor(d-1)/2\rfloor$ gives $\mdr(\cB)=\mdr(\cA)$.
\end{proof}

\section{Computations and verification}\label{sec:verification}

The checks of the computations for defects one, two and three are described in Section~\ref{sec:checks},
and the scope of the verification in Section~\ref{sec:reliability}.  This section lists the computations for the
defect bound, for defect four and for points of multiplicity seven.  All computations for the defect bound are exact; they
are in the directory \texttt{gst1718\_computation/defect\_bound} within the ancillary files.

\paragraph{Local identities.}  Every displayed initial form, weighted initial form and polynomial
identity of Appendices~\ref{sec:residual}--\ref{sec:nonreduced} and of Lemmas~\ref{red5:local},
\ref{tc4:semisimple-trace}, \ref{clm:noninvariant-node} and~\ref{n5:normal-form}--\ref{n5:at-least-two-degree-two}
was recomputed by exact symbolic computation with \textsc{SymPy}
(\texttt{defect\_bound/local/*.py}: the programs written with the text and, independently, the programs
\texttt{audit\_*.py} of an audit, written from the text alone; the table in
\texttt{defect\_bound/README.md} lists which statement each program checks).  This includes the identities $17t(t-1)+2=0$ and $(1-2t)^2=9$ of
Lemma~\ref{lem:L4-smooth-four}, the coefficients $44b^4/25$ of Lemma~\ref{cg:no-physical-five} and
$168A^3B/25$ of Lemma~\ref{cg4:no-triple-four}, and the identity $E^2-D'E+DE'=C\,T_C$ of
Lemma~\ref{n5:normal-form}.

\subsection{The catalogs of critically fixed maps}\label{sec:catalogs}
The program \texttt{catalog\_exact.py} in the directory \texttt{defect\_bound/exact} treats the catalog $\mathcal K_{\mathrm{cf}}$ of
Proposition~\ref{cat:multipliers} and the catalog of Theorem~\ref{ds6:main}.  For four critical points
the interpolation determinant is computed and factored over $\QQ$; for five critical points the saturated
ideal of maximal minors is computed with \textsc{Singular} \cite{Singular} and brought into shape form by the
separating coordinate $u=q+3p$ (\texttt{exact/gen\_shape.py}, outputs \texttt{exact/ar5\_*.out},
\texttt{exact/ds5\_*.out}).  Every component is then treated over its number field with exact rational
arithmetic (\textsc{python-flint} \cite{flint}), as described in the proof of
Proposition~\ref{cat:multipliers}.  The results are:
\begin{center}\small
\begin{tabular}{lrrrl}\toprule
catalog & components & labelled solutions & multipliers & excluded values\\\midrule
$\mathcal K_{\mathrm{cf}}$ & $26$ & $277$ & $867$ & $2,3,4,5,6,7$\\
Theorem~\ref{ds6:main} & $7$ & $1+24+108$ & $425$ & $2,3$\\\bottomrule
\end{tabular}
\end{center}
Here a component is an irreducible factor (four critical points) or a component of the saturated ideal
(five critical points), and a labelled solution is a solution of a labelled normalized system (as in
Section~\ref{sec:catalog}, a map is counted once for each admissible labelled normalization); every component has a one-dimensional interpolation kernel, a map of degree seven
with $P$, $Q$ coprime, and distinct fixed points.  The program \texttt{exact/analyse\_multipliers.py}
evaluates the characteristic polynomials at the excluded values and computes the reciprocal pairs of
Proposition~\ref{cat:multipliers}(b): the gcd is $(3x-7)(4x-7)$, and the three components carrying $7/3$
or $7/4$ have, in our normalizations,
\[
A:\ \frac PQ=\frac{x^5-\frac37x^7}{x^2-\frac37},\qquad
C:\ \frac PQ=\frac{x^4-\frac17x^7}{x^3-\frac17},
\]
with noncritical fixed points the roots of $x^4-\frac43x^2+1$ and of $x^3+1$, and for $B$ a map over
$\QQ(t)$, $t^2-t+1=0$, whose noncritical fixed points are the roots of
$x^4-\frac43(t+1)x^3+2tx^2+\frac43(1-2t)x+t-1$.

\paragraph{The incidence certificate.}  The program \texttt{exact/incidence\_exact.py} proves
Lemma~\ref{clm:incidence}.  It uses the prime $\ell=10009$, for which all the polynomials above split
into distinct linear factors and all fixed points of each map stay distinct.  It enumerates the
$108$ marked pairs and the $6\,696$ matrices; the largest number of distinct intersection points on a
line is $13$ in $888$ cases, $14$ in $3\,912$ and $15$ in $1\,896$.

\paragraph{One free critical point.}  The subresultant computations of
Lemma~\ref{ofc4:three-critical} and Theorem~\ref{n4:main} (the gcd $t+\frac13$ and the four sextics $C(t)$,
with $a_7$ proportional to $C$) are exact computations over $\QQ[t]$ (\texttt{defect\_bound/local/ofc\_n4.py}).  For
Lemma~\ref{n5:uniform}, \texttt{defect\_bound/n5/gen.py} writes the coefficient systems of all six local-degree profiles,
and \textsc{Singular} returns the reduced Gr\"obner basis $\{1\}$ in every case (\texttt{n5/profiles.out}).

\paragraph{Defect four and multiplicity seven.}  The searches of Section~\ref{sec:nu4}, the
classification of Section~\ref{sec:sept} and their checks use the census, the engines and the search
programs of Sections~\ref{sec:census}--\ref{sec:search}; the directory \texttt{gst1718\_computation}
contains them together with the records of the defect-four search
(\texttt{gtree}), the classification of transversal structures (\texttt{sept}), the statuses of the census
with multiplicities up to $7$ (\texttt{scan7}) and the computations of Lemma~\ref{lem:LF}
(\texttt{check/lambdaF}).  About $18.2$ million children were examined at defect four, from
$15\,471$ roots.  The general-engine statuses of the $12$-line lattices with $b=43$ and multiplicities
at most $6$, which include the cores with a line with eight points used only at defect four, were
recomputed after two random relabellings of the lines: $323$ lattices, $177$ of them with a line with
eight points, $646$ runs and no disagreement (\texttt{check/nu4\_recheck}).

\paragraph{Dependencies.}  The results of \cite{GST16} quoted in Section~\ref{sec:reductions} and in
Appendix~\ref{app:A} (the reduction to $(C_7)$, the residual lemma, the vertex
formula, the local initial-form lemma, the Riccati identity and the cubic local tests), Ziegler's
restriction theorem \cite{Zie89}, Abe's bound \cite{Abe}, the census generators of \cite{AS} and the
primary decomposition and Gr\"obner basis routines of \textsc{Singular} are used as inputs.

\section*{Computational data}
{\sloppy The ancillary directories are \texttt{gst1718\_computation} and \texttt{mss\_comparison}.  The
directory \texttt{gst1718\_computation} contains the programs, the census of twelve-line lattices, the
statuses of the cores, the per-root records of the searches for defects one to four, the classification
of Section~\ref{sec:sept}, the computations for Lemma~\ref{lem:LF} and for Section~\ref{sec:18}, the checks
of Sections~\ref{sec:checks}, \ref{sec:nu}, \ref{sec:nu4} and~\ref{sec:sept}, and the computations for the
defect bound (\texttt{defect\_bound}); its README gives the reproduction commands, and \texttt{SHA256SUMS}
the hashes of all files.  The directory \texttt{mss\_comparison} contains the comparison of
Remark~\ref{rem:mss}.  Software: \textsc{Singular}~4.3.2 \cite{Singular}, \textsc{nauty}~2.8.8,
\textsc{python-flint}~0.9 \cite{flint}, \textsc{SymPy}~1.14.\par}

\section*{Acknowledgements}
We thank Felipe Monteiro, Henrique N.~S\'a Earp and Tom\'as S.~R.~Silva for their helpful correspondence
and valuable discussions, and for making the data of their computations publicly available.

\section*{AI assistance}
ChatGPT and Claude were used for exploratory computations, literature searches, programming and interactive feedback on earlier drafts of the manuscript. The authors take full responsibility for all content.

\appendix
\section{The residual quartic and the first defect bounds}\label{app:A}
This appendix and the next prove Theorem~\ref{thm:defect}.  Here we set up the critical foliation and its
residual quartic, prove the curvilinearity of $Z$ and the trace bound, exclude the doubled conic, give the
first bounds for the nonreduced residuals and show that a reduced residual has $\nu\le5$.  Appendix~\ref{app:B}
completes the proof.

\subsection{The critical foliation}\label{sec:foliation}
Let $\cA$ be a critical realization of a $7$-candidate.  We call the lines of $\cA$ \emph{old} and all
other lines \emph{new}.  A nonzero $\theta\in D_0(\cA)_7$ is unique up to a scalar, and as a section of
$\T\cA(7)$ it has a finite zero scheme $Z$ of length $\nu$ \cite[Lemma~4.1]{GST16} (for multiplicities at most
$6$ this is also Lemma~\ref{lem:gen} and Section~\ref{sec:nu}, whose proofs use only the value of
$b(\cA)$, $r_H\le8$ and $\mdr(\cA)=7$):
\[
 0\longrightarrow\cO\longrightarrow\T\cA(7)\longrightarrow\cI_Z(-2)\longrightarrow0 .
\]
The derivation $\theta$ defines a saturated foliation of degree $7$ on $\PP^2$ for which every old line
is invariant; $Z$ is its \emph{logarithmic zero scheme}.  Its inflection polynomial \cite{Per01} has degree
$3\cdot7=21$ and factors as $F_\theta=Q_\cA R$, where the \emph{residual quartic} $R$ satisfies
\[
 R\in H^0(\cI_Z^2(4)),\qquad\operatorname{Supp}Z\subseteq\operatorname{Sing}R
\]
\cite[Theorem~4.3 and Corollary~4.4]{GST16}.  We call $\cA$ together with $\theta$, $Z$ and $R$ the \emph{critical
seventeen-line setting}.

\paragraph{Local formulas.}
In an affine chart a local vector field $X=A\partial_x+B\partial_y$ defining the foliation has
inflection determinant
\begin{equation}\label{ar:D}
 D=\det(X,\nabla_XX)=A^2B_x+AB(B_y-A_x)-B^2A_y ,
\end{equation}
a local equation of the inflection divisor; replacing $X$ by $uX$ with $u$ a unit multiplies $D$ by $u^3$.
Off the arrangement the ideal of $Z$ is $(A,B)$; we call $X$ the \emph{physical} field.  At a point on
exactly one old line $\{y=0\}$ write $X=f\partial_x+yg\partial_y$; then $I_Z=(f,g)$ and
\begin{equation}\label{ar:smooth}
 P=\frac Dy=f^2g_x+fg^2+fygg_y-fgf_x-yg^2f_y .
\end{equation}
At an old vertex of multiplicity $m$, let $q$ be the product of the $m$ old linear forms through it,
$E=x\partial_x+y\partial_y$, $\eta=q_y\partial_x-q_x\partial_y$ and $X=aE+b\eta$; then $I_Z=(a,b)$ and
\begin{equation}\label{ar:vertex}
 \frac Dq=m\bigl(bX(a)-aX(b)\bigr)-m(m-1)a^2b-kb^3,\qquad \operatorname{Hess}_q(\eta,\eta)=qk,\quad\deg k=2m-4
\end{equation}
\cite[Lemma~4.5]{GST16}.  Both frame fields vanish at the vertex, so at a vertex of
$Z$ the residual satisfies $R\in\mathfrak m_pI_{Z,p}^2$.

\paragraph{Restrictions to a line.}  For a new line $L$ which is not invariant, the normal component of
the foliation along $L$ is a nonzero section of $\cO_L(7)$, the \emph{normal restriction}; every point of
$Z$ on $L$ and every tangency of the foliation with $L$ is a zero of it.  For an invariant line the
tangential component is a nonzero section of $\cO_L(8)$, the \emph{tangential restriction}, which
vanishes at the singular points of the foliation on $L$.

\subsection{Critically fixed maps of degree seven}\label{sec:catalog}
A rational map $f\colon\PP^1\to\PP^1$ is \emph{critically fixed} if every critical point is a fixed
point; these maps are classified combinatorially in \cite{CGNPP15,Hlu19}.  At a fixed point $v$ with
multiplier $\lambda=f'(v)\ne1$ the holomorphic fixed-point index is $1/(1-\lambda)$, and for $f\ne\mathrm{id}$
the indices of all fixed points add up to $1$ \cite{Mil06}.  We denote by $\mathcal K_{\mathrm{cf}}$
the set of critically fixed maps of degree $7$ with at most five critical points, each of local degree at
most $5$.  Such a map has $\sum(m_v-1)=12$, hence at least three critical points.  It is determined by
the positions and local degrees of its critical points: if $f_i=P_i/Q_i$ ($i=1,2$), with binary forms
$P_i,Q_i$ of degree $7$, have the same $r$ critical points with the same local degrees, then the form
$P_1Q_2-P_2Q_1$ of degree $14$ vanishes to order at least $m_v$ at every critical point $v$ (at a finite
critical point $[v:1]$ the form $(x-vz)^{m_v}$ divides $P_i-vQ_i$, and $P_1Q_2-P_2Q_1=(P_1-vQ_1)Q_2-(P_2-vQ_2)Q_1$;
at $v=\infty$, $z^{m_\infty}$ divides $Q_1$ and $Q_2$), hence to total order $12+r\ge15$, so it is zero and
$f_1=f_2$.  Normalizing three critical points to $0,1,\infty$
leaves finitely many solutions of the interpolation systems, in which the critical points are labelled:
$277$ solutions of the labelled normalized systems, one with three critical points, $34$ with four and
$242$ with five (a map is counted once for each admissible labelled normalization, so these are neither
Möbius classes nor necessarily distinct maps).  The map with three critical points has local degrees $5,5,5$ and is
$x^5(2x^2-7x+7)/(7x^2-7x+2)$, with noncritical multipliers $35/8$ (three times) and $14/5$ (twice).

\begin{proposition}\label{cat:multipliers}
Let $f\in\mathcal K_{\mathrm{cf}}$.
\begin{enumerate}[label=\textup{(\alph*)}]
\item No noncritical fixed point of $f$ has multiplier in $\{2,3,4,5,6,7\}$.
\item Let $f,g\in\mathcal K_{\mathrm{cf}}$ (possibly equal) and let $\lambda$, $\mu$ be multipliers of
noncritical fixed points of $f$ and of $g$ with $(1-\lambda)(1-\mu)=1$.  Then $\{\lambda,\mu\}=\{7/3,7/4\}$.
The maps with a noncritical multiplier $7/3$ are, up to Möbius transformations, the map $A$ with
critical points $0,1,-1,\infty$ of local degrees $5,3,3,5$, and the map $B$ with critical points
$0,1,\tfrac{1+\alpha}2,\infty$ of local degrees $4,4,4,4$, $\alpha^2=-3$ (its Galois conjugate, with
$\tfrac{1-\alpha}2$, is Möbius equivalent to it via $x\mapsto1-x$); all their four noncritical
multipliers equal $7/3$.  The map with a noncritical multiplier $7/4$ is, up to Möbius transformations,
the map $C$ with critical points $0,1,\tfrac{-1+\alpha}2,\tfrac{-1-\alpha}2,\infty$ of local degrees
$4,3,3,3,4$, which is defined over $\QQ$; all its three noncritical multipliers equal $7/4$.
\end{enumerate}
\end{proposition}

\begin{proof}[Computer-assisted proof]
Normalize three critical points to $0,1,\infty$ and write $f=P/Q$ with
$P=\sum_{i=m_0}^7a_ix^i$ and $Q=\sum_{j=0}^{7-m_\infty}b_jx^j$; a finite critical point $v$ of local degree
$m_v$ gives the linear conditions $(x-v)^{m_v}\mid P-vQ$.  With three critical points the kernel is one
line over $\QQ$.  With four critical points $0,1,t,\infty$ the conditions form a square matrix; after
removing all powers of the factors $t$ and $t-1$, the determinants of the four local-degree profiles have
irreducible factors of total degree $34$.
With five critical points $0,1,p,q,\infty$ the ideal of maximal minors, saturated by
$pq(p-1)(q-1)(p-q)$, is radical of degree $242$ in total over the six local-degree profiles.  For every
irreducible component, with function field a number field $K$, we computed exactly over $K$ the kernel
(always one-dimensional), checked that $\deg P=7$, that $P$ and $Q$ are coprime and that the fixed points
are distinct, and computed the characteristic polynomial $\chi(x)\in\QQ[x]$ of multiplication by
$W/Q^2$, $W=P'Q-PQ'$, on $K[x]/(F)$ viewed as a $\QQ$-vector space, where $F$ is the polynomial of the noncritical fixed points.  The
roots of $\chi$ are the noncritical multipliers of the maps of the component.  None of the polynomials
$\chi$ vanishes at $2,\dots,7$, which proves~(a).  For~(b), let $\chi_0$ be the squarefree part of the
product of all $\chi$; the condition $(1-\lambda)(1-\mu)=1$ means $\mu=\lambda/(\lambda-1)$, and the gcd of
$\chi_0(x)$ and $(x-1)^{\deg\chi_0}\chi_0\bigl(x/(x-1)\bigr)$ is $(x-\tfrac73)(x-\tfrac74)$.  The components
on which $\chi$ vanishes at $7/3$ or at $7/4$ are the three stated ones, and there $\chi$ is a power of
$x-\tfrac73$, respectively $x-\tfrac74$.  The computation is described in Section~\ref{sec:catalogs}.
\end{proof}

\subsection{The residual quartic}\label{sec:residual}

Throughout Appendices~\ref{app:A} and~\ref{app:B} we are in the critical seventeen-line setting of Section~\ref{sec:foliation}: $\cA$ is a
critical realization of a $7$-candidate, $\mathcal F$ its saturated foliation of degree seven,
$F_\theta=Q_\cA R$, and $Z$ the logarithmic zero scheme, of length $\nu$.  We write $I_p$ for the ideal of
$Z$ at $p$, $n_p=\length Z_p$, and $h(H)=r_H$ for the number of vertices of an old line $H$.  An equation of
$R$ in a local chart is understood up to a unit; this causes no difficulty in ideal-membership statements.

\subsubsection{Local consequences of the residual lemma}\label{sec:lemma46}
We use the following consequences of \cite[Lemma~4.6]{GST16}.  If $p\in Z$ and $\operatorname{mult}_pR=2$,
then $n_p=1$ and $p$ is not an old vertex.  If $\operatorname{mult}_pR=3$, then $Z_p$ is curvilinear; if
moreover $n_p\ge3$, the initial cubic of $R$ at $p$ is a cube, and it is the cube of the tangent line of
$Z_p$.  For the last point write $I_p=(y-\varphi(x),x^{n_p})$ with $\operatorname{ord}\varphi\ge2$; if
the initial cubic were $\ell^3$ with $\ell\ne y$ up to a scalar, then $R(x,\varphi(x))$ would have order
three, whereas $R\in I_p^2$ forces it to be divisible by $x^{2n_p}$.  These statements concern logarithmic
foliations and are stronger than the bare inclusion $R\in I_p^2$.

\subsubsection{A high-multiplicity vertex on every old line}
\begin{lemma}\label{lem:hm}
Every old vertex has multiplicity at most eight, and every old line contains a vertex of multiplicity at
least four.
\end{lemma}
\begin{proof}
A line outside a pencil of $m$ old lines meets those lines in $m$ distinct vertices.  Since the
arrangement is not a pencil and $h(H)\le8$, this gives $m\le8$.

Suppose an old line $H$ contains only double and triple vertices.  The sixteen intersection weights on $H$,
together with $h(H)\le8$, force exactly eight triple vertices.  Normalize the homogeneous logarithmic
derivation by subtracting an Euler multiple so that it annihilates the equation of $H$.  Its restriction
to $H$ is a binary derivation of degree seven for the Ziegler multirestriction, whose eight weights are
all two.  This restriction is nonzero: otherwise the normalized derivation is divisible by the equation of
$H$, giving a non-Euler logarithmic derivation of degree six, contrary to minimality.

Choose an affine coordinate $t$ on $H$ whose point at infinity avoids the eight directions, and write
$p(t)$, $q(t)$ for the two coefficients, each of degree at most seven.  If the directions are
$t_1,\dots,t_8$, the double-divisibility conditions are
\[
 p(t_i)-t_iq(t_i)=p'(t_i)-t_iq'(t_i)=0 .
\]
Let $T(t)=\prod_{i=1}^8(t-t_i)$.  The first conditions give $p-tq=cT$.  Differentiating at each $t_i$ gives
$q(t_i)=-cT'(t_i)$, hence $q=-cT'$ and $p=c(T-tT')$.  Its coefficient of degree eight is $-7c$, so $c=0$.
Then $p=q=0$, a contradiction.
\end{proof}

\subsubsection{A global Riccati budget}
\begin{proposition}\label{prop:riccati}
No old line divides $R$.  Every repeated linear factor of $R$ is a new line which is not invariant under
$\mathcal F$.
\end{proposition}
\begin{proof}
Let $\ell$ be invariant and suppose $\ell^2\mid F_\theta$.  The Riccati identity \cite[Lemma~6.2(a)]{GST16}
provides a point $p_0\in\ell$ such that, in an affine chart having $p_0$ at infinity on $\ell$,
\[
 \ell=(u=0),\qquad X=uf\partial_u+g\partial_v,\qquad f(0,v)=0 .
\]
For clarity, the identity is independent of the degree of the foliation: the coefficient of $u$ in $D$ is
$g_0(f_0g_0'-g_0f_0'-f_0^2)$, where $f_0=f(0,v)$, $g_0=g(0,v)$.  Its vanishing gives either $f_0=0$ or
$(g_0/f_0)'=1$, hence $g_0=(v-v_0)f_0$.  The latter says that the foliation and the pencil through
$(0,v_0)$ agree to first order along $\ell$; making that point infinite gives the asserted chart.  The
nonzero tangential restriction $t_\ell$ has degree eight.  At a point $q\ne p_0$ of $\ell$ where $m$ lines
of an invariant arrangement $\cB$ meet, including $\ell$, write $q_\cB=u\prod_{i=1}^{m-1}(v-t_iu)$ and
$X=aE+b\eta$.  Restriction to $u=0$ gives
\[
 f_0=a_0+(m-1)v^{m-2}b_0=0,\qquad g_0=va_0-v^{m-1}b_0=-mv^{m-1}b_0 .
\]
Thus $\operatorname{ord}_q(t_\ell)\ge m-1$.

If $\ell=H$ is old and divides $R$, apply this to $\cB=\cA$.  The sum of $m-1$ over the vertices on $H$ is
sixteen.  At most seven can occur at $p_0$, by Lemma~\ref{lem:hm}.  The other zeros of $t_H$ therefore have
total degree at least nine, contradicting its degree eight.

If $\ell=L$ is new, invariant and $L^2\mid R$, apply the same argument to $\cB=\cA\cup\{L\}$.  The weights
$m-1$ on $L$ now sum to seventeen, counting the old lines.  At most eight occur at $p_0$, so again the
remaining tangential divisor has degree at least nine.  This proves both assertions.
\end{proof}

\subsubsection{Invariant new lines}
\begin{lemma}\label{lem:invline}
Every invariant new line meets the old arrangement in at least five distinct points.
\end{lemma}
\begin{proof}
Let $I$ be the line and let $r$ be its number of old intersections.  Add $I$ to the arrangement.  Its
Ziegler multirestriction to $I$ has total weight seventeen; its weights are the old multiplicities on $I$,
each at most eight by Lemma~\ref{lem:hm}.  Thus $r\ge3$, and the multirestriction is balanced.  If its
exponents are $d_1\le d_2$, Abe's bound \cite[Theorem~0.1 of the arXiv version]{Abe} gives
\[
 d_1+d_2=17,\qquad d_2-d_1\le r-2 .
\]
The minimal derivation of degree seven restricts nontrivially, so $d_1\le7$.  To see nontriviality, choose
coordinates with $I=(z=0)$ and subtract an Euler multiple so that $\theta(z)=0$.  If both remaining
coefficients vanished on $I$, division by $z$ would give a logarithmic derivation of degree six for the
old arrangement, contradicting minimality.  If $r\le4$, the displayed bound gives $d_1\ge8$.  Thus $r\ge5$.
\end{proof}

\subsubsection{Local support}
\begin{lemma}\label{lem:sole}
Suppose $L$ is a non-invariant line and $R=L^2S$ with $L\nmid S$.  Then no point of $Z$ lies on
$L\setminus V(S)$.
\end{lemma}
\begin{proof}
At such a point $p\in Z$ the local multiplicity of $R$ is two, so $n_p=1$ and $p$ is not an old vertex
(Section~\ref{sec:lemma46}).  If $p$ lies on just one old line, the quadratic part of~\eqref{ar:smooth} is a
nonzero scalar times $yf_1$, where $y=0$ is the old line; its being the square of the equation of $L$ would
force $L$ to be that old line, which is impossible.

It remains to consider a reduced zero off the arrangement.  If $J$ is the linear part of $X$ and
$w=(x,y)^{\mathsf t}$, then
\[
 D_2=\det(Jw,J^2w)=\det J\,\det(w,Jw).
\]
A nonzero square as initial form forces a nontrivial Jordan block with $L$ as its eigenline.  Choose
$L=(y=0)$ and
\[
 J=\begin{pmatrix}\lambda&\beta\\0&\lambda\end{pmatrix},\qquad \lambda\beta\ne0 .
\]
If $B(x,0)=b_jx^j+O(x^{j+1})$ is its first nonzero normal term, then $j\ge2$ and~\eqref{ar:D} gives
\[
 D(x,0)=\lambda^2jb_jx^{j+1}+O(x^{j+2})\ne0 ,
\]
contradicting $L\mid F_\theta$.  Hence the normal restriction vanishes identically, making $L$ invariant,
again a contradiction.
\end{proof}

\begin{lemma}[A local intersection inequality]\label{lem:lci}
Let $I\subset\CC[[x,y]]$ be an $\mathfrak m$-primary complete-intersection ideal.  If
$G=\prod_i\ell_i^{a_i}\in I^2$, where the $\ell_i$ are distinct linear forms, then
\begin{equation}\label{eq:lci}
 2\length\bigl(\CC[[x,y]]/I\bigr)\le\sum_ia_i\length\bigl(\CC[[x,y]]/(I,\ell_i)\bigr).
\end{equation}
Consequently, if $R=\prod_iL_i^{a_i}$ globally, then $2\length Z\le\sum_ia_i\length(Z\cap L_i)$.
\end{lemma}
\begin{proof}
Choose general constant linear combinations $f,g$ of the two generators of $I$.  They still generate $I$,
the curve $f=0$ has no common component with $G$, and its restriction to each $\ell_i$ has the minimal
vanishing order of $I|_{\ell_i}$.  Thus
\[
 i_0(f,g)=\length\bigl(\CC[[x,y]]/I\bigr),\qquad i_0(f,\ell_i)=\length\bigl(\CC[[x,y]]/(I,\ell_i)\bigr).
\]
Modulo $f$, the membership $G\in(f,g)^2$ implies $G\in(g^2)$.  Taking valuations on the branches of $f=0$
and summing with their multiplicities gives $i_0(f,G)\ge2i_0(f,g)$; additivity of intersection
multiplicities gives~\eqref{eq:lci}.  Summation over $Z$ proves the global assertion.
\end{proof}

\begin{lemma}\label{lem:fourlines}
Let $P$ be a product of four distinct linear forms through the origin.  If $P\in I^2$ for an
$\mathfrak m$-primary complete-intersection ideal $I\subset\CC[[x,y]]$, then $\length(\CC[[x,y]]/I)\le4$.
If $I$ contains a smooth parameter, the bound is three.
\end{lemma}
\begin{proof}
If $I=(w,t^n)$ for regular parameters $w,t$, the coefficients of $w^0$ and $w^1$ in $P$ must be divisible
by $t^{2n}$ and $t^n$, respectively.  The curve $w=0$ is tangent to at most one of the four lines.  If it
is tangent to none, $P(0,t)$ has order four, so $n\le2$.  If it is tangent to one, the coefficient of $w$
has order exactly three, so $n\le3$.

Otherwise $I\subset\mathfrak m^2$.  Since $P$ has order four, both minimal generators have order two: if
only one did, the initial quartic of every element of $I^2$ would be a scalar multiple of a square.  Their
quadratic initial forms are coprime, since a common factor would have its square dividing $P$.  Their
local intersection multiplicity is therefore $2\cdot2=4$.
\end{proof}

\subsubsection{Curvilinearity}
\begin{lemma}\label{lem:curvi}
Every nonempty local scheme $Z_p$ is curvilinear.
\end{lemma}
\begin{proof}
The ideal $I_p$ is a local complete intersection generated by the two coefficients of the logarithmic
section in a local frame.  A zero-dimensional ideal in a regular local ring of dimension two which is not
contained in the square of the maximal ideal contains a regular parameter, and its quotient is
curvilinear.  Suppose then that both frame coefficients vanish to order at least two.

Off the old arrangement, the physical field $X=A\partial_x+B\partial_y$ then has order at least two, and
its inflection determinant $D=AX(B)-BX(A)$ has order at least five.  Since the equation of the old
arrangement is a unit, this contradicts $\operatorname{ord}_pR\le4$.

At a smooth point of the old arrangement, put its unique old line at $y=0$ and write
$X=f\partial_x+yg\partial_y$.  Both $f$ and $g$ have order at least two, and~\eqref{ar:smooth} has order at
least five, again a contradiction.

At an old ordinary double point, the physical field has order at least three, so $D$ has order at least
eight, and division by the two old linear factors leaves residual order at least six.

Finally, at an old vertex of multiplicity $m\ge3$ use $X=aE+b\eta$, with $\eta$ homogeneous of degree
$m-1$, and~\eqref{ar:vertex}.  Here $a,b\in\mathfrak m_p^2$.  The term $a^2b$ has order at least six,
$kb^3$ has order at least eight, and the terms involving $\eta$ have order at least seven.  In the remaining
Euler part, the terms of degree four of $bE(a)-aE(b)$ cancel, so multiplication by $a$ gives order at least
seven.  Thus the residual has order at least six.  This exhausts all local types.
\end{proof}

\subsubsection{Reduced residuals and the nonreduced list}
\begin{theorem}\label{thm:reduced6}
If the residual quartic $R$ is reduced, then $\nu\le6$.  Equality forces $R=L_1L_2L_3L_4$ with four lines in
general position and $Z$ equal to their six reduced intersection points.
\end{theorem}
\begin{proof}
We consider the five partitions of four describing the degrees of the irreducible components of $R$.

An irreducible quartic has arithmetic genus three.  If all its singularities have multiplicity two, there
are at most three, each with local $Z$-length at most one.  A triple point uses at least three of the genus
budget, so it is the only singularity.  If its local $Z$-length $n$ were at least three, the local
initial-form statements of Section~\ref{sec:lemma46} would give
\[
 R=y^3+R_4(x,y),\qquad I_p=(y-\varphi(x),x^n),\qquad\operatorname{ord}\varphi\ge2 .
\]
The coefficient of $x^4$ in $R_4$ is nonzero, since otherwise $y$ divides $R$.  Hence $R(x,\varphi(x))$ has
order exactly four, contradicting its divisibility by $x^{2n}$.  Thus a triple point contributes at most
two.  A point of multiplicity four would make the quartic a union of lines, contrary to irreducibility.

For an irreducible cubic times a line, a smooth cubic contributes at most three intersection points, all of
residual multiplicity two.  A singular cubic has one node or cusp.  If the line misses it, the bound is
$1+3=4$.  If the line passes through it and the initial cubic of the product is not a cube, the local
contribution is at most two and at most one further intersection remains.  The remaining case is a
cuspidal cubic with its cuspidal tangent.  In affine coordinates its equation is $R=y(y^2-x^3)$, with no
other singularity.  If $I_p=(y-\varphi(x),x^n)$ and $n\ge3$, the coefficient of $w$ after
$y=w+\varphi(x)$ is $3\varphi(x)^2-x^3$, of order three, and the membership
$R\in(w^2,wx^n,x^{2n})$ gives $n\le3$.

Two irreducible conics are smooth and meet in at most four distinct points.  Every singularity of their
product has multiplicity two, so the bound is four, including tangential intersections.  For a smooth conic
and two distinct lines, if the common point of the lines is off the conic there are at most five double
points.  If it lies on the conic, the initial cubic there contains the two distinct line tangents and is
not a cube; its contribution is at most two, with at most one further intersection on each line, giving at
most four.

Finally, four lines in general position have six nodes.  With exactly three concurrent lines, there is a
non-cube triple point of local contribution at most two and three further nodes, giving five.  With four
concurrent lines, Lemma~\ref{lem:fourlines} gives four.  Only the four lines in general position can give
total length six, and then all six nodes occur with length one.
\end{proof}

\begin{corollary}\label{cor:nonreduced}
A nonreduced residual quartic has one of the forms
\[
 C^2,\quad L^2C,\quad L^4,\quad L^3M,\quad L^2M^2,\quad L^2MN,
\]
where $C$ is a smooth irreducible conic and distinct line symbols denote distinct lines.  Every repeated line
in this list is new and not invariant, and the lines $M,N$ are new (Proposition~\ref{prop:riccati}).
\end{corollary}
\begin{proof}
A repeated component of a quartic has degree one or two, and an irreducible conic is smooth.  Listing the
factorizations of a quartic with a repeated component gives the six forms.  The last assertion is
Proposition~\ref{prop:riccati}.
\end{proof}

\subsection{The trace bound}\label{sec:trace}

\begin{lemma}[Three small dual-curve obstructions]\label{lem:dual}
No irreducible curve in the dual plane contains all but $k$ of the seventeen old dual points in any of the
cases
\[
 (\deg C,k)=(3,\le1),\qquad(2,\le4),\qquad(1,\le5).
\]
The same conclusion holds if the indicated degree is replaced by an upper bound, and a constant Gauss image
cannot contain two distinct old dual points.
\end{lemma}
\begin{proof}
By Lemma~\ref{lem:hm} the old arrangement has no vertex of multiplicity greater than eight, and it is not a
pencil.  Thus a dual line cannot contain twelve old dual points, proving the assertion in degree one.

For an irreducible conic, an old vertex of multiplicity at least four contains at most two of the dual
points on the conic, and hence at least two omitted dual points.  Distinct such old vertices use distinct
pairs of omitted lines, since two lines meet only once.  So there are at most $\binom k2$ such vertices,
each covering at most two included old lines.  Every one of the $17-k$ included lines contains a vertex of
multiplicity at least four (Lemma~\ref{lem:hm}), so
\[
 17-k\le2\binom k2 ,
\]
which is false for $0\le k\le4$.

For an irreducible cubic with at most one omitted point, an old vertex of multiplicity at least four must be
exactly quadruple, with three included lines and the unique omitted line $H_0$ (with no omitted line it
cannot exist).  Every included line therefore meets $H_0$ at such a quadruple vertex.  These vertices
partition the sixteen included lines into triples, contradicting $3\nmid16$.  Lower degrees are covered by
the preceding cases.
\end{proof}

\begin{theorem}[Total trace at most four]\label{thm:trace}
Let $L$ be a line which is not invariant under $\mathcal F$.  Then $\length(Z\cap L)\le4$.
\end{theorem}
\begin{proof}
Let $t_p=\length(Z_p\cap L)$ and $T=\sum_{p\in L}t_p$.  The normal restriction of $\mathcal F$ to $L$ is a
nonzero section of $\cO_L(7)$; write $\beta_p$ for its vanishing order, so that $\sum_p\beta_p=7$.  The
homogeneous Gauss map, which sends a point $w$ to the tangent line of the foliation, is represented by
$w\times\theta(w)$, whose three coordinates have degree eight.  Let $b_p$ be their common vanishing order
after restriction to $L$, and $b=\sum_pb_p$.  After cancellation of this common divisor the restriction to
$L$ is a morphism $L\to\check\PP^2$ of degree $N=8-b$, whose image we call the \emph{Gauss image}; its
values at all points of $L$, including limit values at the base points, belong to the image.  In an affine
chart centered at $p$, the Gauss coordinates are $(-B,A,xB-yA)$ for the local field $A\partial_x+B\partial_y$,
so $b_p$ is exactly the smaller of the orders of the two restricted coefficients.  The following local
comparisons are independent of the chart and of nonvanishing local scalar multiples of the field.

At an old vertex $p$ of multiplicity $m\ge2$, choose coordinates with $L=(y=0)$ and write the local old
product as $Q$, a product of $m$ distinct linear forms, none equal to $y$.  A logarithmic frame is
$E=x\partial_x+y\partial_y$, $\eta=Q_y\partial_x-Q_x\partial_y$.  If the section is $aE+e\,\eta$, its local
zero ideal is $(a,e)$, and the normal restriction is $-e(x,0)Q_x(x,0)$, with $Q_x(x,0)$ a nonzero
constant times $x^{m-1}$.  Hence
\begin{equation*}\tag{V}\label{eq:V}
 \beta_p\ge t_p+m-1,\qquad b_p\ge t_p+1 ,
\end{equation*}
the second inequality because both frame fields vanish at the vertex.

At a smooth old point put the old line $H=(y=0)$ and $L=(x=0)$; these two lines are transverse.  The field
has the form $f\partial_x+yg\partial_y$ and the logarithmic zero ideal is $(f,g)$.  Thus
$t_p=\min\{\operatorname{ord}_Lf,\operatorname{ord}_Lg\}$, and both $\beta_p$ and $b_p$ are at least $t_p$.
If $b_p=t_p$, then necessarily $\operatorname{ord}_Lf=t_p$, whereas $yg$ has order at least $t_p+1$, so the
limiting Gauss line at $p$ is $H$ itself.  Consequently an old line which is \emph{not} in the Gauss image,
and whose intersection with $L$ is a smooth old point $p$, satisfies
\begin{equation*}\tag{S}\label{eq:S}
 b_p-t_p\ge1 .
\end{equation*}
This includes smooth old singularities of the foliation not in $Z$, where $t_p=0$.  At points off the
arrangement the logarithmic frame is an ordinary tangent frame, so again $\beta_p,b_p\ge t_p$.

Suppose now $T\ge5$ and set
\[
 B_0=7-T=\sum_p(\beta_p-t_p),\qquad V=b-T=\sum_p(b_p-t_p).
\]
All summands are nonnegative, hence $0\le B_0\le2$ and $V\ge0$.  The Gauss map has degree
\begin{equation*}\tag{G}\label{eq:G}
 N=8-b=8-T-V\le3-V .
\end{equation*}
Let $k$ be the number of old dual points missing from the Gauss image.  Every regular intersection of $L$
with an old line gives that line as its Gauss value, so only singular intersection points can cause an
omission.  At an old vertex of multiplicity $m$, all $m$ incident old lines may conservatively be counted as
omitted, and~\eqref{eq:V} gives $m\le(\beta_p-t_p)+(b_p-t_p)$.  At a smooth old point, an omitted old line
costs at least one in $b_p-t_p$ by~\eqref{eq:S}.  Summing over the distinct intersections on $L$ gives
\begin{equation*}\tag{O}\label{eq:O}
 k\le B_0+V .
\end{equation*}

If $V=0$, then~\eqref{eq:V} excludes all old vertices on $L$ and~\eqref{eq:S} excludes any missing smooth
old line, so $k=0$; the Gauss image is irreducible of degree at most three and contains all seventeen old
dual points.  A vertex of multiplicity at least four, which exists on every old line by
Lemma~\ref{lem:hm}, gives four collinear points of that curve, which is impossible for an irreducible curve
of degree two or three; and a line as image would make all old lines concurrent.  If $V=1$,
\eqref{eq:G} and~\eqref{eq:O} give image degree at most two and $k\le3$, contradicting
Lemma~\ref{lem:dual}.  If $V=2$, they give degree at most one and $k\le4$, also impossible.  Finally
$V\ge3$ makes the Gauss map constant: its three coordinates cannot have a common divisor of degree greater
than eight, so $N\ge0$ and $V=3$, $N=0$; then~\eqref{eq:O} leaves at least $17-5=12$ distinct old dual
points in the image, contradicting Lemma~\ref{lem:dual}.  Thus $T\le4$.
\end{proof}

\subsection{The doubled conic}\label{sec:conic}
This section proves that $R=C^2$, with $C$ a smooth conic, is impossible when $\nu>0$.  The argument
associates a rational self-map to the conic and classifies its fixed-point indices.

\begin{lemma}\label{lem:conic1}
If $R=C^2$ with $C$ smooth, then $Z$ is reduced, contained in $C$, and disjoint from every old line.
\end{lemma}
\begin{proof}
The support is on $C$ since $\operatorname{Supp}Z\subseteq\operatorname{Sing}R$, and residual multiplicity
two implies reducedness and excludes old vertices (Section~\ref{sec:lemma46}).  At a zero on just one old
line $y=0$, write $X=f\partial_x+yg\partial_y$.  The nonzero quadratic part of~\eqref{ar:smooth} is a scalar
times $yf_1$.  Its being a square forces
\[
 f_1=ay,\qquad g_1=bx+cy,\qquad ab\ne0 .
\]
For arbitrary higher jets,~\eqref{ar:smooth} gives
\[
 [y^2](D/y)=a^2b,\qquad [x^2y](D/y)=0 ,
\]
where $[\,\cdot\,]$ denotes a coefficient.  A smooth conic tangent to $y=0$ has local equation
$y-\kappa x^2+O(xy,y^2)$ with $\kappa\ne0$; a unit times its square, whose $y^2$ coefficient is $a^2b$, has
$x^2y$ coefficient $-2\kappa a^2b\ne0$.  This contradiction excludes the remaining old points.
\end{proof}

\subsubsection{The second-intersection map}
Parametrize $C$ by $[s:t]\mapsto[s^2:st:t^2]$ and write the homogeneous field of degree seven as
$\theta=(a,b,c)$.  With pullbacks to $C$ denoted by a subscript, define the binary forms of degree fifteen
\[
 P=sb_C-ta_C,\qquad Q=sc_C-tb_C .
\]
They are unchanged by adding an Euler multiple.  The tangent line of the foliation at the point of $C$
with parameter $[s:t]$ has coordinates $[tQ:-(tP+sQ):sP]$, and its value at $[u^2:uv:v^2]$ is
$(tu-sv)(Qu-Pv)$.  Thus its second intersection with $C$ defines
\begin{equation}\label{eq:secondmap}
 f\colon C\longrightarrow C,\qquad[s:t]\longmapsto[P:Q]
\end{equation}
after cancellation of the common divisor.  The pair $(P,Q)$ is not identically zero, since the foliation is
saturated, and it has the same common vanishing order as the restricted projective field.  Each point of
$Z$ contributes one to the common divisor, since its linearization is invertible.  If $e$ is the sum of the
common orders at the old singularities on $C$, then
\begin{equation}\label{eq:degf}
 \deg f=15-\nu-e .
\end{equation}
A projective change of coordinates conjugates this construction.  In affine coordinates $C\colon y=x^2$,
with parameter $t=x$, it reads
\begin{equation}\label{eq:faff}
 f(t)=\frac{B(t,t^2)-tA(t,t^2)}{A(t,t^2)} .
\end{equation}
The map is not the identity: that would make $C$ invariant, and an invariant conic $y=x^2$ satisfies
$D|_C=2A|_C^3$; since $C\mid D$, the field would vanish identically along $C$, contradicting saturation.

\begin{lemma}\label{lem:conic2}
Every fixed point of $f$ is simple.  The possible multipliers and the corresponding fixed-point indices are
\[
\begin{array}{c|cccc}
 f'(p)&0&-1&1/2&3/5\\\hline
 (1-f'(p))^{-1}&1&1/2&2&5/2
\end{array}
\]
Every point of $Z$ is fixed with multiplier zero.
\end{lemma}
\begin{proof}
Use~\eqref{eq:faff}, centered at the point in question.  At a point of $Z$, the invertible linearization and
the quadratic inflection form, which is a nonzero square, give
\[
 A=\lambda x+\beta y+O(2),\qquad B=\lambda y+dx^2+\text{(other terms of order at least two)},\qquad\lambda\beta\ne0 .
\]
The coefficient of $t^3$ in $D(t,t^2)$ is $2\lambda^2d$.  Since $C\mid D$, $d=0$.  The numerator
of~\eqref{eq:faff} has order at least three and the denominator order one, proving the assertion at $Z$.
At a regular point of the foliation, on or off the old arrangement, which is fixed, $A(0)\ne0$ and
$B(0)=0$, and~\eqref{ar:D} gives $B_x(0)=0$, so $f'(0)=-1$.

At an old singularity with restricted zero order one which is fixed, the tangent eigenvalue is nonzero:
\[
 J=\begin{pmatrix}\alpha&\beta\\0&\gamma\end{pmatrix},\qquad\alpha\ne0 .
\]
The quadratic inflection form vanishes because $Q_\cA C^2$ has order at least three there.  Hence either
$\gamma=0$ or $J=\alpha I$.  Writing $B=\gamma y+dx^2+\cdots$, the coefficient of $t^3$ in $D(t,t^2)$ is
\[
 \alpha\bigl((\alpha+\gamma)d+\gamma(\gamma-\alpha)\bigr) ,
\]
and its vanishing gives $d=0$ in both cases.  After cancellation, $f'(0)=(\gamma-\alpha)/\alpha\in\{-1,0\}$.

Consider now old singularities with restricted zero order at least two.  At an old smooth or double point,
the linearization is nonzero because the logarithmic section does not vanish on the arrangement.  Its kernel
must be the tangent to $C$ and is an old line; put this line at $y=0$.  Then $B=yg$ with $g(0)\ne0$, while
$A(t,t^2)=O(t^2)$.  The numerator in~\eqref{eq:faff} has leading term $g(0)t^2$, so the map after
cancellation does not fix the point.

At an old vertex of multiplicity $m\ge3$, use~\eqref{ar:vertex}.  If $a(0)\ne0$, the linearization is scalar
and the restricted order is one, already treated.  Otherwise $b(0)=b_0\ne0$.  Multiplying $X$ by the unit
$b_0/b$ leaves $f$ unchanged and lets us assume $b=b_0$ constant in the coefficient calculations below.

For $m=3$, a fixed value requires the quadratic part $X_2$ to be radial or zero on the tangent of the conic.
Since $\det(E,X_2)=-3b_0q$, that tangent is an old line.  Normalize the old directions to $y$, $x$, $x-y$,
keeping $C\colon y=x^2$.  The quadratic field can be written
\[
 X_2=x(ux+vy)\partial_x+y(wx+(u+v-w)y)\partial_y,\qquad u-w\ne0 .
\]
Its reduced inflection quadratic equals
\[
 -(u-w)\bigl(uwx^2+2u(u+v-w)xy+v(u+v-w)y^2\bigr) .
\]
It must be a nonzero multiple of $y^2$.  If $u\ne0$, the first two coefficients force $w=0$ and $v=-u$,
making the whole expression zero.  If $u=0$, then $w\ne0$, and invariance of the old lines gives
$B(t,t^2)=wt^3+O(t^4)$, whereas $A(t,t^2)=O(t^3)$.  The numerator also has leading term $wt^3$, so the point
is not fixed.

For $m\ge4$, the quadratic residual form is $-m(m-2)b_0a_1^2$.  Thus $a_1=\alpha y$ with $\alpha\ne0$.  Set
$h=[x^2]a$ and $K=[y^2](D/q)$.  The residual is a unit times $(y-x^2)^2$, which forces
\begin{equation}\label{eq:conic12}
 [x^3](D/q)=0,\qquad[x^2y](D/q)=-2K,\qquad[x^4](D/q)=K .
\end{equation}
Terms of positive degree of the unit cannot affect these three coefficients.

If $m=4$, the $x^3$ coefficient is $4b_0^2\alpha\eta_y(1,0)$; its vanishing makes $y=0$ an old line.
Normalize $q=x^3y+{}$terms divisible by $y^2$.  Equation~\eqref{ar:vertex} yields
\[
 [y^2](D/q)=-8\alpha^2b_0,\quad[x^2y](D/q)=-12\alpha b_0(h+b_0),\quad[x^4](D/q)=4b_0(h+b_0)(3b_0-h).
\]
Equations~\eqref{eq:conic12} give $h+b_0=-4\alpha/3$ and $b_0/\alpha=1/24$.  Along the conic the denominator
and the numerator of~\eqref{eq:faff} are $-\frac\alpha3t^3+O(t^4)$ and $-4b_0t^4+O(t^5)$.  The multiplier is
therefore $12b_0/\alpha=1/2$.

If $m=5$, set $q_0=q(1,0)$.  The corresponding coefficients are
\[
 [y^2](D/q)=-15\alpha^2b_0,\qquad[x^2y](D/q)=-25\alpha b_0h,\qquad[x^4](D/q)=-10b_0h^2-25\alpha b_0^2q_0 .
\]
They force $h=-6\alpha/5$ and $b_0q_0=3\alpha/125$.  The denominator and the numerator are then
$-\alpha t^3/5+O(t^4)$ and $-5b_0q_0t^4+O(t^5)$, and the multiplier is $25b_0q_0/\alpha=3/5$.

For $m\ge6$, the Hamiltonian terms enter too late to affect these coefficients, and the equations become
\[
 h=-\frac{2(m-2)}{2m-5}\,\alpha,\qquad(m-3)h^2=(m-2)\alpha^2 .
\]
They would require $4(m-3)(m-2)=(2m-5)^2$, whose right side is larger by one, a contradiction.  The
computations involve only $a_1$, $h$ and $b_0$: omitted higher jets carry an excessive power of $y$ or an
excessive total degree for the indicated coefficients.  This exhausts the old singularities.  All
multipliers found differ from one, so every fixed point is simple and has the displayed index.
\end{proof}

\begin{theorem}\label{thm:conic}
In the critical seventeen-line setting with $\nu>0$, the residual quartic is not $C^2$ for a smooth conic $C$.
\end{theorem}
\begin{proof}
Suppose first that $f$ is not constant.  For a rational self-map of degree at least one other than the
identity, with simple fixed points, the residue theorem gives
\begin{equation}\label{eq:indexsum}
 \sum_{f(p)=p}\frac1{1-f'(p)}=1 ;
\end{equation}
indeed, choose infinity not fixed and apply the residue theorem to $dz/(z-f(z))$, whose residue at infinity
is $-1$.  All summands are positive by Lemma~\ref{lem:conic2}, and the $\nu$ points of $Z$ contribute one
each.  For $\nu\ge2$ this is impossible.  For $\nu=1$ there are $\deg f+1\ge2$ fixed points, all simple, so
an additional positive summand again gives a contradiction.

It remains that $f$ is constant, with its unique fixed point $p_0\in Z$ and $\nu=1$.
Equation~\eqref{eq:degf} gives $e=14$.  Every point of $C\cap\cA$ is a singularity of the foliation:
otherwise the tangent of the foliation there is its old line $H$, whose second intersection with $C$ would
be $f(p)=p_0$, contradicting $p_0\notin\cA$ (Lemma~\ref{lem:conic1}).  None of these old singularities is
fixed.

Let $r_p$ be the restricted zero order at such an old point.  At a smooth old point,
$I_p(C,\cA)\le2\le2r_p$.  At an old vertex, write $X=aE+b\eta$ as above.  Since $p$ is not fixed,
$b(p)\ne0$, because $b(p)=0$ would give $a(p)\ne0$, a scalar linearization and a fixed point.  In conic
coordinates
\[
 \det(E,X)=-mbq,\qquad\det(E,X)|_C=t\bigl(B(t,t^2)-tA(t,t^2)\bigr) .
\]
Since the map after cancellation does not fix $p$, its numerator has order exactly $r_p$.  Thus
$I_p(C,\cA)=r_p+1\le2r_p$.  B\'ezout now gives the contradiction
\[
 34=C\cdot\cA=\sum_{p\in C\cap\cA}I_p(C,\cA)\le2\sum_pr_p=2e=28 .\qedhere
\]
\end{proof}

\subsection{First bounds for the nonreduced residuals}\label{sec:firstbounds}
By Theorem~\ref{thm:trace}, every non-invariant line $L$ satisfies $T_L:=\length(Z\cap L)\le4$.  By
Corollary~\ref{cor:nonreduced}, the repeated lines are new and not invariant, and Lemma~\ref{lem:sole}
locates the zeros on a doubled line at its intersections with the other components.

\begin{proposition}\label{prop:nonred91}
If the residual quartic is $L^2C$ with $C$ a smooth irreducible conic, or $L^2M^2$ with $L\ne M$, or
$L^2MN$ with $L,M,N$ distinct, then $\nu\le5$.  If $L$ meets $C$ transversely, then $\nu\le4$.
\end{proposition}
\begin{proof}
For $L^2C$ the zero scheme is supported on $L\cap C$.  At a transverse crossing the cubic tangent cone of
the residual is not a cube, so the local initial-form statements of Section~\ref{sec:lemma46} give length at
most two.  This proves the bound four in the transverse case.

If $L$ is tangent to $C$, there is only one support point $p$.  Put $n=\length Z_p$ and
$t=\length(Z_p\cap L)$.  We show $n\le\max\{4,t+1\}$.  For $n\ge3$ the cubic tangent cone forces the
curvilinear tangent to be $L$.  Choose local analytic coordinates $(u,v)$ at $p$ with $L=(u=0)$ and
$C=(u+v^2=0)$ up to a unit, and $I=(u-\varphi(v),v^n)$; this is legitimate here, since only memberships
in ideals and lengths are used, and they do not depend on the choice of local coordinates.  Membership $u^2(u+v^2)\in I^2$ implies
\[
 v^{2n}\mid\varphi^2(\varphi+v^2),\qquad v^n\mid\varphi(3\varphi+2v^2) .
\]
If $r=\operatorname{ord}\varphi\ge n$ (including $\varphi\equiv0$), then $t=n$, and the bound $n=t\le4$ comes
from the trace.  Otherwise $t=r$.  If $r=2$, the two conditions imply $n\le4$; if $r\ge3$, the first gives
$n\le r+1$.  Since $t\le4$, $n\le5$ follows.

For $L^2M^2$, the zero scheme is supported at $p=L\cap M$.  By Lemma~\ref{lem:curvi} its tangent coincides
with at most one of the two lines, so one of the traces $t_L,t_M$ is one.  The other is at most four, and
Lemma~\ref{lem:lci} gives $n\le t_L+t_M\le5$.

For $L^2MN$ with the three lines not concurrent, the possible support points are $L\cap M$, $L\cap N$ and
$M\cap N$.  The first two have non-cube cubic tangent cones and hence local length at most two; the third
has residual multiplicity two and local length at most one.  Thus $\nu\le5$.

If the three lines are concurrent at $p$, then $Z$ is supported at $p$; put $t_J=\length(Z_p\cap J)$.  If
the curvilinear tangent is $L$, then $t_M=t_N=1$ and $2n\le2t_L+2$, so $n\le5$.  If it is none of the three
lines, all traces are one and $n\le2$.  If the tangent is $M$ (the case $N$ is identical), then
$t_L=t_N=1$, whence $2n\le3+t_M$.  When $M$ is not invariant, $t_M\le4$.  When $M$ is invariant,
Lemma~\ref{lem:invline} gives at least five distinct old intersections on $M$.  The tangential field on $M$,
of degree eight, vanishes at these five points, and the trace of $Z$ at $p$ adds at least $t_M$ to its
order: if $p$ is one of the old points, the order there is at least $1+t_M$ \cite[Lemma~6.4(b)]{GST16}.  Thus $5+t_M\le8$, i.e.\
$t_M\le3$.  In either case $n\le3$.
\end{proof}

\begin{proposition}\label{prop:L3M7}
If $R=L^3M$ with $L\ne M$, then $\nu\le7$.  If $\nu=7$, then $p_0=L\cap M$ belongs to $Z$, the local scheme
$Z_{p_0}$ has length three and tangent $M$, and
\[
 \length(Z_{p_0}\cap L)=1,\qquad\sum_{p\ne p_0}\length(Z_p\cap L)=3,\qquad\sum_{p\ne p_0}\length Z_p=4 .
\]
\end{proposition}
\begin{proof}
The support of $Z$ is contained in $L=\operatorname{Sing}R$.  Put $t_p=\length(Z_p\cap L)$ and
$T=\sum_pt_p\le4$.  Away from $p_0$ the local residual is a unit times the cube of the equation of $L$, so
Lemma~\ref{lem:lci} gives $2n_p\le3t_p$.

If $p_0\notin Z$, summation gives $\nu\le6$.  If $p_0\in Z$ and its curvilinear tangent is not $M$, then
$\length(Z_{p_0}\cap M)=1$, Lemma~\ref{lem:lci} gives $2n_{p_0}\le3t_{p_0}+1$, and therefore
$2\nu\le3T+1\le13$, again $\nu\le6$.

It remains to consider a local scheme tangent to $M$.  Then $t_{p_0}=1$.  Use affine linear coordinates with
$L=(x=0)$ and $M=(y=0)$, so that the residual is a unit times $x^3y$.  Curvilinearity gives
$I_{p_0}=(y-\varphi(x),x^n)$ with $\operatorname{ord}\varphi\ge2$.  In the coordinate $w=y-\varphi(x)$, the
coefficient of $w$ in $x^3(w+\varphi(x))$ is $x^3$ times a unit, and membership in $(w^2,wx^n,x^{2n})$ forces
$n\le3$.  The remaining trace is at most $T-1\le3$, and its total length at most
$\lfloor3(T-1)/2\rfloor\le4$.  Hence $\nu\le3+4=7$, and equality forces the stated local conditions.
\end{proof}

\begin{proposition}\label{prop:L4off}
If $R=L^4$, then $Z$ is disjoint from the complement of the old arrangement.
\end{proposition}
\begin{proof}
The line $L$ is new and not invariant (Corollary~\ref{cor:nonreduced}).  Let $p\in Z$ lie off the old
arrangement.  In affine coordinates centered at $p$ the physical field $X=A\partial_x+B\partial_y$ has
isolated zero ideal $(A,B)$, and~\eqref{ar:D} reads $D=U\ell^4$ with $U(0)\ne0$, where $\ell$ is a linear
equation of $L$.  The zero scheme is curvilinear by Lemma~\ref{lem:curvi}.

If its length is one, the linearization is invertible.  Vanishing of the quadratic inflection form forces
it to be scalar, $J=\lambda I$ with $\lambda\ne0$.  Put $L=(x=0)$.  If $A(0,y)=a_ky^k+\cdots$ were nonzero,
then $k\ge2$ and $B(0,y)=\lambda y+\cdots$, and the coefficient of $y^{k+1}$ in $D(0,y)$ is
$-ka_k\lambda^2\ne0$, contradicting $D|_L=0$.  Thus $L$ would be invariant.

It remains to consider length $n\ge2$, when the linearization has rank one.  If it is semisimple, choose
linear coordinates with $A=x+O(2)$, $B=O(2)$, the kernel being $x=0$.  Write
$A_2=a_{20}x^2+a_{11}xy+a_{02}y^2$.  Vanishing of $D_3$ forces $B_2=b_{11}xy$, and with
$B_3=b_{30}x^3+b_{21}x^2y+b_{12}xy^2+b_{03}y^3$,
\[
 D_4=x\bigl((a_{02}b_{11}-b_{03})y^3+(-a_{20}b_{11}+b_{11}^2+b_{21})x^2y+2b_{30}x^3\bigr) .
\]
If it is nilpotent, use $A=O(2)$, $B=x+O(2)$.  Here $D_3=-x^2\partial_yA_2=0$, so $A_2=a_{20}x^2$, and with
$B_2=b_{20}x^2+b_{11}xy+b_{02}y^2$, $A_3=a_{30}x^3+a_{21}x^2y+a_{12}xy^2+a_{03}y^3$,
\[
 D_4=-x^2\bigl(3a_{03}y^2+(2a_{12}-2a_{20}b_{02})xy+(a_{20}^2-a_{20}b_{11}+a_{21})x^2\bigr) .
\]
In either case the kernel divides the nonzero initial form of order four, which is a fourth power;
therefore $L=(x=0)$.

Put $a(y)=A(0,y)$, $b(y)=B(0,y)$, $c(y)=A_x(0,y)$ and $d(y)=B_x(0,y)$.  Non-invariance means $a\not\equiv0$.
The identity $D(0,y)=0$ is
\begin{equation}\label{eq:lineid}
 a^2d-acb+abb'-b^2a'=0 .
\end{equation}
In the semisimple case $c$ is a unit, $d(0)=0$ and $k=\operatorname{ord}a\ge2$.  If $b\not\equiv0$, put
$w=b/a$; then~\eqref{eq:lineid} gives $d=cw-aww'$.  Holomorphy of $d$ and $d(0)=0$ force
$\operatorname{ord}b\ge k+1$: if $\operatorname{ord}b<k$, the pole of $cw$ cannot be canceled by $aww'$, of
strictly larger order, and equality would give $d(0)\ne0$.  Thus $w(0)=0$.  The smooth curve $A=0$ is
$x=\phi(y)$ with $\operatorname{ord}\phi=k$ and $\phi=-a/c+O(a^2)$, and Taylor expansion gives
\[
 B(\phi(y),y)=b-da/c+O(a^2)=a^2ww'/c+O(a^2) ,
\]
so $n=\operatorname{ord}B|_{A=0}\ge2k$.  If $b\equiv0$, then~\eqref{eq:lineid} gives $d\equiv0$ and the same
bound follows.  On the other hand, restricting $D=Ux^4$ to $A=0$ gives $-B^2A_y=U\phi^4$; since
$A_y(\phi(y),y)=-A_x(\phi(y),y)\phi'(y)$ has order $k-1$, this gives $2n+k-1=4k$, so $n=(3k+1)/2<2k$, a
contradiction.

In the nilpotent case $d$ is a unit and $c(0)=0$.  The function $b$ cannot vanish identically, since
then~\eqref{eq:lineid} would give $a^2d=0$.  Put $k=\operatorname{ord}b\ge2$, $h=\operatorname{ord}a\ge2$
and $r=a/b$.  The identity becomes $c=rd-br'/r$.  If $h<k$, the pole of $rd$ cannot be canceled by $br'/r$,
whose order is $k-1\ge1$; if $h=k$, then $c(0)=r(0)d(0)\ne0$.  Consequently $h>k$ and
$r'/r=(h-k)/y+O(1)$.  The smooth curve $B=0$ has equation $x=\psi(y)=-b/d+O(b^2)$, with
$\operatorname{ord}\psi=k$, and
\[
 A(\psi(y),y)=a-cb/d+O(b^2)=\frac{b^2}d\,\frac{r'}r+O(b^2) .
\]
The leading term has exact order $2k-1$, so $n=2k-1$.  But restricting $D=Ux^4$ to $B=0$ gives
$A^2B_x=U\psi^4$, forcing $2n=4k$.  This contradiction proves the proposition.
\end{proof}

\begin{lemma}[The strict vertex inequality]\label{lem:vertexineq}
Let $p\in Z$ be an old vertex.  Then $R\in\mathfrak m_pI_p^2$.  If moreover near $p$ the residual is a unit
times $\ell^k$, where $\ell$ is a linear equation of a line $L$, then, with $t_p=\length(Z_p\cap L)$,
\[
 2n_p+1\le k\,t_p .
\]
\end{lemma}
\begin{proof}
In the frame $X=aE+b\eta$ both frame fields vanish at the vertex and $I_p=(a,b)\subset\mathfrak m_p$, so
$X(a),X(b)\in\mathfrak m_pI_p$, and~\eqref{ar:vertex} gives $R\in\mathfrak m_pI_p^2$, including for
multiplicity $m=2$.  By Lemma~\ref{lem:curvi} write $I_p=(s,t^{n_p})$ in regular parameters, and choose a
general smooth generator $g=s+ct^{n_p}$, so that the order of $\ell$ on the smooth curve $g=0$ is exactly
$t_p$: when the restriction of $\ell$ to $s=0$ has order less than $n_p$, that order is unchanged, and
otherwise a general $c$ gives order exactly $n_p$.  Restricting $R\in\mathfrak m_pI_p^2$ to $g=0$ gives order
at least $2n_p+1$, whereas $\ell^k$ has order exactly $kt_p$.
\end{proof}

\subsubsection{The cubic local tests}\label{sec:cubictests}
We use the local tests of \cite[Lemmas~6.12, 6.13 and Appendix~A, (C1)--(C3)]{GST16}.  Their application here
is local: at a point of $L$ other than $L\cap M$, the residual $L^3M$ is a unit times $\ell^3$.  The tests use
divisibility of the physical inflection determinant by $\ell^3$ and the nonzero leading coefficient after
division by the local equation of the old arrangement.  The calculation of \cite{GST16} already allows an arbitrary
unit contributed by old lines not through the point, and multiplication by the additional unit $M$ changes
neither the vanishing equations nor the nonzero leading coefficient.  In the test at a smooth old point the
induction of \cite[Lemma~6.13(a)]{GST16} runs through order seven (it uses the coefficients of $xy^j$ and
$x^2y^j$, $2\le j\le7$, of $P=D/y$; for the determinant $D$ itself the exponents of $y$ are one larger) and
forces normal vanishing order at least eight, which contradicts the normal degree seven here.  Explicitly, at a zero where $R$ is a unit times $\ell^3$ and $L$ is not
invariant:
\begin{enumerate}[label=\textup{(\roman*)},leftmargin=*]
\item there is no zero off the old arrangement or at a smooth old point;
\item at an old vertex of multiplicity $m$, with $\beta=\operatorname{ord}(b|_L)$ in the frame
$X=aE+b\eta$, the normal vanishing order is $(m-1)+\beta$; a nonreduced zero has $\beta\ge8$ if $m=2$,
$\beta\ge6$ if $m=3$, and $\beta\ge3$ if $m=4$, with length exactly three if $m=4$ and $\beta=3$.
\end{enumerate}

\begin{theorem}\label{thm:L3M6}
If $R=L^3M$ with $L\ne M$, then $\nu\le6$.
\end{theorem}
\begin{proof}
Suppose $\nu=7$.  By Proposition~\ref{prop:L3M7}, $Z_{p_0}$ has length three, is tangent to $M$ and has trace
one on $L$; outside $p_0$ the trace on $L$ is three and the total length four.  By the cubic local tests,
every support point outside $p_0$ is an old vertex, where Lemma~\ref{lem:vertexineq} gives
$3t_p\ge2n_p+1$, i.e.\ $n_p\le\lfloor(3t_p-1)/2\rfloor$.  The only way to have total trace three and total
length four is therefore a single old vertex $q\ne p_0$ with $t_q=3$ and $n_q=4$: the partitions $2+1$ and
$1+1+1$ give total length at most three.

Let $m$ be the old multiplicity at $q$ and $\beta=\operatorname{ord}(b|_L)\ge t_q=3$.  The normal section on
$L$ has degree seven; its order at $q$ is $(m-1)+\beta$, and its order at $p_0$ is at least one.  For $m=2$
the local test gives $\beta\ge8$, and for $m=3$ it gives $\beta\ge6$, both exceeding the budget.  For $m=4$
the budget forces $\beta=3$, and the local test gives $n_q=3$, contradicting $n_q=4$.  Finally $m\ge5$ gives
normal cost at least $4+3+1=8$.  So $\nu\ne7$, and $\nu\le6$ by Proposition~\ref{prop:L3M7}.
\end{proof}

\subsection{The reduced residual with six nodes}\label{sec:sixnodes}

\begin{theorem}\label{thm:sixnodes}
A reduced residual quartic satisfies $\nu\le5$.
\end{theorem}

By Theorem~\ref{thm:reduced6}, $\nu\le6$, with equality only when
\begin{equation}\label{eq:sixnodes}
 R=L_1L_2L_3L_4,\qquad Z=\{L_i\cap L_j:1\le i<j\le4\},
\end{equation}
the four lines being in general position and the six points reduced.  The lines $L_i$ are new by
Proposition~\ref{prop:riccati}.  We assume~\eqref{eq:sixnodes} until the end of the proof.

\begin{lemma}\label{lem:sixoff}
All six points of~\eqref{eq:sixnodes} lie off the old arrangement.  Their foliation singularities are
nondegenerate, with the two residual lines as distinct eigendirections.
\end{lemma}
\begin{proof}
At an old vertex in $Z$, $R\in\mathfrak m_pI_p^2$ (Lemma~\ref{lem:vertexineq}), so its multiplicity would be at
least three, whereas a node of~\eqref{eq:sixnodes} has multiplicity two.  At a smooth old point put its old
line at $y=0$ and write $X=f\partial_x+yg\partial_y$; the local ideal is $(f,g)$, and its reducedness makes
$f_1,g_1$ independent.  The quadratic part of~\eqref{ar:smooth} is a nonzero multiple of $yf_1$, and for an
ordinary node of two residual lines this would force one of them to be the old line $y=0$, which is
impossible.  Off the arrangement the coefficients of $X$ generate $I_Z$, so their linearization $J$ is
invertible at a reduced zero, and the identity
\begin{equation}\label{eq:in2D}
 \operatorname{in}_2D=\det(J)\det(w,Jw),\qquad w=(x,y)^{\mathsf t},
\end{equation}
shows that the two residual tangent directions are the two distinct eigendirections of $J$.
\end{proof}

\begin{lemma}\label{lem:tangency}
Every non-invariant residual line has normal tangency orders $(3,2,2)$ at its three residual nodes and no
other tangencies.  Its node of order three lies on an invariant residual line.  More generally, at a
nondegenerate singularity on a non-invariant straight inflection component, with tangent eigenvalue
$\lambda$, transverse eigenvalue $\mu$ and normal vanishing order $n\ge2$,
\begin{equation}\label{eq:resonance}
 \mu=-(n-1)\lambda .
\end{equation}
\end{lemma}
\begin{proof}
First we exclude a simple tangency away from the three nodes.  Put the residual line at $y=0$ and write
$X=f\partial_x+g\partial_y$.  The identity $D=f^2g_x+fgg_y-fgf_x-g^2f_y$ vanishes on this line.  If $g(x,0)$
vanishes to order $n\ge1$ at a regular point with $f(0)\ne0$, the first term has order $n-1$ and the others
order at least $n$, a contradiction.  So every tangency is singular.

A singularity off the arrangement is in $Z$, so it gives no additional point.  At a smooth old point the
initial inflection form is a nonzero product of the old line and the residual line, and~\eqref{eq:in2D} makes
the residual line an eigendirection, excluding a simple normal zero.  At an old vertex of multiplicity at
least three, the linearization preserves at least three directions and is scalar; again its normal component
on the residual line has no linear term.  At an old double point, put its old lines on the axes:
$X=xa\partial_x+yb\partial_y$.  A simple normal zero would make the diagonal linearization non-scalar.  Since
the inflection polynomial has order three there, its two eigenvalues cannot both be nonzero; after exchanging
the axes, $a(0)=\lambda\ne0$ and $b(0)=0$.  Direct expansion gives
\[
 D/(xy)=xa^2b_x+yabb_y-xaba_x-yb^2a_y+ab(b-a),\qquad\operatorname{in}_1\bigl(D/(xy)\bigr)=-\lambda^2b_y(0)y ,
\]
so the cubic initial form is zero or a multiple of $xy^2$, contrary to the required product of three distinct
lines.

The normal section on a non-invariant line has degree seven.  Its three residual nodes are eigendirection
singularities, contributing at least two each.  The only possibilities are $(3,2,2)$, or $(2,2,2)$ with an
additional simple zero, which has just been excluded.

At a nondegenerate singularity as in the last assertion, in adapted coordinates $f(x,0)=\lambda x+O(x^2)$ and
$g(x,0)=bx^n+O(x^{n+1})$ with $b\ne0$, and $\partial_yg(0)=\mu$; the coefficient of $x^{n+1}$ in $D(x,0)=0$
gives~\eqref{eq:resonance}.  If both residual lines at a node were non-invariant, applying this in both
directions would give $(n_i-1)(n_j-1)=1$, hence $n_i=n_j=2$.  A node of order three therefore has an
invariant partner.
\end{proof}

\begin{lemma}\label{lem:sixmult}
Every old vertex has multiplicity at most five.
\end{lemma}
\begin{proof}
For an old vertex $p$ of multiplicity $m$, represented by $p\in\CC^3$, the angular polynomial
$W_p(w)=\det(p,w,\theta(w))$ has degree eight.  Each old line through $p$ divides $W_p$, so their product
$Q_p$ divides it.  At the six residual nodes the foliation is singular and $Q_p$ is nonzero, so $W_p/Q_p$, of
degree $8-m$, vanishes at those six points.  No nonzero form of degree at most two vanishes at all six nodes:
its restriction to any of the four residual lines would have three distinct zeros, forcing that line to be
a component.  If $m\ge6$ we obtain $W_p=0$; the tangent line of every regular leaf would then pass through
$p$, making the foliation a pencil of lines, whose inflection polynomial vanishes identically, contrary to
$F_\theta\ne0$.
\end{proof}

\begin{lemma}\label{lem:invresidual}
An invariant residual line $I$ meets the old arrangement in exactly five distinct points.  Its tangential zero
divisor consists of these five points and its three residual nodes, all with multiplicity one.  There is a
rational map $f_I\colon I\to I$ of degree seven such that
\begin{enumerate}[label=\textup{(\arabic*)}]
\item its eight fixed points are the five old intersections and the three residual nodes;
\item at an old intersection incident with $m_i$ old lines its local degree is $m_i$, and these are all its
critical points with $m_i>1$; moreover $1\le m_i\le5$ and $\sum_{i=1}^5m_i=17$;
\item at a fixed point with tangent eigenvalue $\lambda$ along $I$ and normal eigenvalue $\mu$,
\begin{equation}\label{eq:indexratio}
 \frac\mu\lambda=\frac1{1-f_I'} ;
\end{equation}
in particular, at a residual node whose other residual line is not invariant and has tangency order $n$
there, the multiplier is $n$.
\end{enumerate}
\end{lemma}
\begin{proof}
The tangential section on $I$ is nonzero of degree eight.  It vanishes at the $r$ old intersections and at
the three residual nodes, so $r\le5$, and Lemma~\ref{lem:invline} gives equality; the eight required zeros
exhaust the divisor and are simple.

Put $I=(z=0)$ and arrange $\theta(z)=0$ by subtracting an Euler multiple.  Set $A(x,y)=\theta(x)|_{z=0}$ and
$B(x,y)=\theta(y)|_{z=0}$.  In the chart $(u,t)=(z/x,y/x)$ the normal linear coefficient is $-uA(1,t)$ and
the tangential restriction is $B(1,t)-tA(1,t)$; the latter has the eight simple zeros just described.  All
these singularities are nondegenerate.  This is Lemma~\ref{lem:sixoff} at residual nodes.  At an old point
with at least two old lines, the additional invariant line $I$ makes the linearization scalar, and the simple
tangential zero makes that scalar nonzero.  At an intersection with one old line, the quadratic initial
inflection form is the product of that line and $I$, so~\eqref{eq:in2D} forces invertibility.  Thus $A,B$ have
no common zero, since a common zero would be a tangential zero with vanishing normal eigenvalue, and
$f_I=[A:B]$ is a map of degree seven whose fixed points are exactly the zeros of $xB-yA$.

The $m_i$ old lines through a point of $I$ become parallel lines $\ell(x,y)=c_jz$ in the affine complement of
$I$.  Their invariance and $\theta(z)=0$ imply that $\theta(\ell)$ is divisible by $\prod_j(\ell-c_jz)$; on
$I$, the target coordinate $\ell(A,B)$ is therefore divisible by $\ell^{m_i}$, so the local degree is at
least $m_i$.  Since
\[
 \sum_{i=1}^5(m_i-1)=17-5=12=2\deg f_I-2 ,
\]
Riemann--Hurwitz gives equality and no further critical points; $m_i\le5$ is Lemma~\ref{lem:sixmult}.
Formula~\eqref{eq:indexratio} follows from the chart formula, and at a node whose other residual line has
tangency order $n$, \eqref{eq:resonance} applied to that line reads $\lambda=-(n-1)\mu$, so $f_I'=n$.
\end{proof}

\begin{proof}[Proof of Theorem~\ref{thm:sixnodes}]
Suppose that~\eqref{eq:sixnodes} holds.  Suppose first that one residual line is not invariant.  By
Lemma~\ref{lem:tangency}, its node of order three has an invariant partner $I$, and by
Lemma~\ref{lem:invresidual} the map $f_I$, which belongs to $\mathcal K_{\mathrm{cf}}$, has multiplier three
at that node.  This contradicts Proposition~\ref{cat:multipliers}(a).  Hence all four residual lines are
invariant.

Let $f_i\in\mathcal K_{\mathrm{cf}}$ be the map of $L_i$.  At $L_i\cap L_j$ let $\lambda_i,\lambda_j$ be the
multipliers of $f_i,f_j$.  The two tangent directions have nonzero eigenvalues, and their normal-to-tangent
ratios are reciprocal, so~\eqref{eq:indexratio} gives
\begin{equation}\label{eq:reciprocal}
 (1-\lambda_i)(1-\lambda_j)=1 .
\end{equation}
By Proposition~\ref{cat:multipliers}(b), $\{\lambda_i,\lambda_j\}=\{7/3,7/4\}$, and every map with a
noncritical multiplier $7/3$ (resp.\ $7/4$) has all its noncritical multipliers equal to $7/3$ (resp.\ $7/4$).
So each line $L_i$ receives one of the two values at all three of its nodes, and two lines meeting at a node
receive different values.  This is a two-colouring of the complete graph on four vertices, which is
impossible already on a triangle.
\end{proof}

\section{The global defect bound}\label{app:B}
We treat the remaining residual factorizations: the fourth power, the cubic--line residual and the other
nonreduced residuals, and then the reduced residuals at defect five.  Theorem~\ref{v4:main} at the end
collects the bounds.

\subsection{A fourth-power residual}\label{sec:L4}
In this section $R=L^4$, where $L$ is a new line which is not invariant
(Corollary~\ref{cor:nonreduced}).  We show that $Z$ is empty.

\subsubsection{Smooth old points}
\begin{lemma}[Smooth old points for a fourth-power residual]\label{lem:L4-smooth-even}
Let $X=f\partial_x+yg\partial_y$ be a holomorphic vector-field germ at
the origin, with $f,g\in\mathfrak m=(x,y)$ and
$I=(f,g)$ curvilinear of finite colength $n>0$. Suppose that
$L=(x=0)$ is not invariant and that
\[
 P=f^2g_x+fg^2+fyg g_y-fg f_x-yg^2f_y=Ux^4,
 \qquad U(0)\ne0.
\]
Then $f\in\mathfrak m^2$, $g_x(0)\ne0$, and, writing
$t=\operatorname{length}\mathcal O/(I,x)$, one has
\[
 n=2t,\qquad \operatorname{ord}_y f(0,y)\ge2t.
\]
In particular, every such local logarithmic zero has even length.
\end{lemma}
\begin{proof}
First $n\ne1$. Otherwise the linear parts
$f_1=ax+by$ and $g_1=cx+dy$ are independent, and the quadratic
part of $P$ is
\[
 f_1^2c-f_1g_1a=-(ad-bc)yf_1\ne0,
\]
contrary to $\operatorname{ord}P=4$. Hence $n\ge2$.

Suppose first that $f_x(0)\ne0$. Set
\[
 a(y)=f(0,y),\quad b(y)=g(0,y),\quad
 c(y)=f_x(0,y),\quad d(y)=g_x(0,y).
\]
Non-invariance gives $a\not\equiv0$. Put $k=\operatorname{ord}a\ge1$.
If $b\not\equiv0$, the meromorphic quotient $w=b/a$ satisfies,
by restriction of $P=Ux^4$ to $L$,
\[
 d=cw-aw^2-yaww'.
\]
If $\ell=\operatorname{ord}b<k$, then $cw$ has order $\ell-k<0$,
whereas both remaining terms have order at least $2\ell-k>\ell-k$.
This contradicts holomorphy of $d$. Thus $w$ is holomorphic. If
$b\equiv0$, the same conclusions below follow with $w=0$ and $d=0$.

The smooth curve $f=0$ has equation $x=\phi(y)$, where
$\operatorname{ord}\phi=k$ and
$\phi=-a/c+O(a^2)$. Taylor expansion now gives
\[
 g(\phi(y),y)=b-da/c+O(a^2)
 =\frac{a^2}{c}w(w+yw')+O(a^2).
\]
Here $O(a^2)$ stands for terms of order at least $2k$; consequently $n=\operatorname{ord}g|_{f=0}\ge2k$
(also when $w(0)=0$).
On the other hand, restricting $P=Ux^4$ to this curve gives
$-yg^2f_y=U\phi^4$. Since
$f_y|_{f=0}=-f_x|_{f=0}\phi'$ has order $k-1$, we obtain
$2n+k=4k$, or $2n=3k$. This contradicts $n\ge2k$.

If $f$ is smooth but $f_x(0)=0$, then $f_y(0)\ne0$ and its zero
curve is $y=\phi(x)$, with order $k\ge2$, unless $\phi\equiv0$.
The latter case makes $P$ identically zero on $f=0$, contradicting
$P=Ux^4$. In the former case, restriction to $f=0$ gives
$k+2n=4$, impossible for $n\ge2$ and $k\ge2$.
Thus $f\in\mathfrak m^2$.

Curvilinearity forces $g$ to be smooth. If $g_x(0)=0$, its zero
curve is $y=\psi(x)$ of order $k\ge2$, or is identically $y=0$.
The latter case gives $P|_{g=0}=0$, again impossible. Otherwise
$g_x|_{g=0}=-g_y|_{g=0}\psi'$ has order $k-1$; the restriction
$P|_{g=0}=f^2g_x=Ux^4$ gives $2n+k-1=4$, impossible for $n\ge2$.
Therefore $g_x(0)\ne0$.

Write $g=0$ as $x=\psi(y)$. This curve cannot equal $L$: otherwise
$P|_{g=0}=f^2g_x$ would force $f|_L=0$, contradicting non-invariance.
Put $k=\operatorname{ord}\psi\ge1$. The same restricted identity now
gives $2n=4k$, so $n=2k$. In particular
$\operatorname{ord}b=k$ and $d$ is a unit.

It remains to bound the normal vanishing order. Put $r=a/b$; restriction
of $P$ to $L$ gives
\[
 c=rd+b\left(1-y\frac{r'}r\right).
\]
The meromorphic function $r'/r$ has at most a simple pole at the
origin. If $h=\operatorname{ord}a<k$, the pole of $rd$ cannot be
canceled by the second term, which has order at least $k$.
If $h=k$, the right side has nonzero constant term, contrary to
$c(0)=0$. Hence $h>k$.

Differentiating the exact formula for $P$ and restricting to $g=0$ yields
\[
 P_x|_{g=0}=fg_x(f_x+yg_y)+f^2g_{xx}.
\]
The left side has exact order $3k$, because
$P_x=4Ux^3+U_xx^4$ and $x=\psi(y)$ has order $k$.
The term $f^2g_{xx}$ has order at least $4k$, whereas
$\operatorname{ord}f|_{g=0}=2k$ and $g_x$ is a unit. Therefore
$(f_x+yg_y)|_{g=0}$ has order exactly $k$. The identity
$g_y|_{g=0}=-g_x|_{g=0}\psi'$ shows that $yg_y|_{g=0}$ has order $k$;
hence $f_x|_{g=0}$ has order at least $k$. Taylor expansion from
$x=0$ to $x=\psi(y)$ now gives $\operatorname{ord}c\ge k$.
The second term in the displayed equation for $c$ also has order
at least $k$, so $\operatorname{ord}r\ge k$ and $h\ge2k$.
Finally $t=\min\{\operatorname{ord}a,\operatorname{ord}b\}=k$,
which proves both assertions.
\end{proof}

\subsubsection{The weighted initial form at a smooth old point}
\label{L4w:section}

The smooth-point calculation admits a further sharpening.

\begin{lemma}\label{lem:L4-smooth-four}
Every point $p\in Z$ lying on exactly one old line satisfies
\[
 \operatorname{length}Z_p=4,
 \qquad \operatorname{length}(Z_p\cap L)=2.
\]
The normal restriction to $L$ vanishes there to order at least five.
\end{lemma}
\begin{proof}
Put the old line at $y=0$ and $L$ at $x=0$, and write
$X=f\partial_x+yg\partial_y$. The exact residual is
\begin{equation}\label{L4w:residual}
 P=f^2g_x+fg^2+fyg g_y-fg f_x-yg^2f_y=Ux^4,
 \qquad U(0)\ne0.
\end{equation}
By Lemma~\ref{lem:L4-smooth-even}, $g_x(0)\ne0$,
$f\in\mathfrak m_p^2$, and, writing
$n=\operatorname{length}Z_p$ and
$t=\operatorname{length}(Z_p\cap L)$, one has $n=2t$ and
$\operatorname{ord}(f|_L)\ge2t$.

The curve $g=0$ has equation $x=\psi(y)$ with
$\operatorname{ord}\psi=t$, and $f|_{g=0}$ has order $2t$.
The exact differentiated identity
\[
 P_x|_{g=0}=fg_x(f_x+yg_y)+f^2g_{xx}
\]
has left-hand side of order exactly $3t$, while its last term has
order at least $4t$. Since $yg_y|_{g=0}=-yg_x\psi'$ has order $t$,
it follows that $f_x|_{g=0}$ has order at least $t$.
Taylor expansion from $g=0$ to $L$ therefore gives
\[
 \operatorname{ord}(f_x|_L)\ge t.
\]

Give $x$ weight $t$ and $y$ weight one. Constant rescalings of $x$
and of the field normalize the weighted initial form of $g$ to
$G=X-Y$, where $X=x$ and $Y=y^t$. The two preceding order estimates
show that the weighted initial form of $f$ has weight at least $2t$.
It has weight exactly $2t$, since its restriction to $g=0$ has
that order. Write it as
\[
 F=\alpha X^2+\beta XY+\gamma Y^2,
 \qquad S:=F(1,1)\ne0.
\]
Taking weight $4t$ in~\eqref{L4w:residual} gives the polynomial identity
\begin{equation}\label{L4w:identity}
 F^2+FG\bigl((1-2\alpha)X-(1+t+\beta)Y\bigr)
 -tYG^2(\beta X+2\gamma Y)=uX^4,
 \qquad u\ne0.
\end{equation}

For clarity, the necessary coefficients of this identity can be solved
without elimination software. Set $Y=1$, $z=X-1$ and
$B=2\alpha+\beta$, so that $F=\alpha z^2+Bz+S$.
The coefficients of $z^0,z^1,z^2,z^4$ are respectively
\[
 S^2,\qquad S(B-t),\qquad S(1-2t),\qquad
 \alpha(1-\alpha).
\]
Comparing with $u(z+1)^4$ and using $S\ne0$ gives
\begin{equation}\label{L4w:necessary}
 u=S^2,
 \qquad S=\frac{1-2t}{6},
 \qquad B=\frac{2-t}{3},
 \qquad \gamma=\alpha-\frac12,
 \qquad \alpha(1-\alpha)=S^2.
\end{equation}
The coefficient of $Y^4$ in~\eqref{L4w:identity} is
\[
 \gamma(\gamma+\beta+1-t)=0.
\]
If $\gamma\ne0$, equations~\eqref{L4w:necessary} force
$\alpha=(7-8t)/6$. The last equality in~\eqref{L4w:necessary}
then becomes
\[
 17t(t-1)+2=0,
\]
which is impossible for a positive integer $t$.
Thus $\gamma=0$ and $\alpha=1/2$. It follows that
$(1-2t)^2=9$, whose only positive-integer solution is $t=2$.
Then $n=4$, $\beta=-1$, and
\[
 F=\frac12X^2-XY.
\]
In particular the weight-four term of $f|_L$ vanishes, so its order,
which is the normal vanishing order, is at least five.
\end{proof}

\begin{corollary}\label{L4w:no-smooth}
If $R=L^4$, its logarithmic zero scheme has no point lying on exactly
one old arrangement line.
\end{corollary}
\begin{proof}
Suppose such a point exists and retain the normalization in the proof
of Lemma~\ref{lem:L4-smooth-four}. Put
\[
 a(y)=f(0,y),\quad b(y)=g(0,y),\quad
 c(y)=f_x(0,y),\quad d(y)=g_x(0,y).
\]
The weighted initial forms $G=X-Y$ and $F=X^2/2-XY$ with $t=2$ give
\[
 b=-y^2+O(y^3),\qquad c=-y^2+O(y^3),\qquad d=1+O(y).
\]
If $a$ vanished identically, $L$ would be invariant. Thus
$h=\operatorname{ord}a$ is finite, and the same weighted calculation
gives $h\ge5$. Set $r=a/b$, whose order is $h-2\ge3$.
The exact identity $P|_L=0$, divided by $ab$, reads
\[
 c=rd+b\left(1-y\frac{r'}r\right).
\]
Its right-hand side has coefficient $h-3$ at $y^2$, because $rd$
has order at least three and $yr'/r=h-2+O(y)$. Its left-hand side
has coefficient $-1$. Hence $h-3=-1$, contradicting $h\ge5$.
\end{proof}

\subsubsection{Old vertices}
\label{l4v3:section}

By Proposition~\ref{prop:L4off}, $Z$ has no point off the old arrangement, and by
Corollary~\ref{L4w:no-smooth} it has no point on exactly one old line.  Thus every support point of $Z$
is an old vertex.  The following local calculations exclude those vertices as well.

\paragraph{The local residual at an ordinary double point}

Take the two old lines as $x=0,y=0$ and write
\[
 X=xf\partial_x+yg\partial_y,\qquad I_{Z,p}=(f,g).
\]
The residual after division by the two old factors is exactly
\begin{equation}\label{l4v3:double-residual}
 P=xf^2g_x+fg(g-f)+yfg g_y-xfg f_x-yg^2f_y.
\end{equation}
Since the remaining old factors are units at $p$, the equation
$R=L^4$ says that $P$ is a unit times the fourth power of the linear
equation of $L$.

\begin{lemma}\label{l4v3:no-reduced-double}
A reduced logarithmic zero at an old ordinary double point cannot
occur when $R=L^4$.
\end{lemma}
\begin{proof}
Write $f_1=ax+by$, $g_1=cx+dy$. Reducedness gives $ad-bc\ne0$.
The cubic part of~\eqref{l4v3:double-residual} is
\begin{equation}\label{l4v3:general-cubic}
 P_3=-(acx^2+2adxy+bdy^2)
       \bigl((a-c)x+(b-d)y\bigr).
\end{equation}
It must vanish identically. The second factor cannot vanish
identically, since that would make $f_1=g_1$. Therefore
$ac=ad=bd=0$. As $ad-bc\ne0$, it follows that $b,c\ne0$ and
$a=d=0$. Thus
\[
 f_1=by,\qquad g_1=cx.
\]
This is the reduced-double-point case which the cubic initial form
alone does not exclude.

The line $L$ is distinct from the two old lines, so constant
rescalings of $x,y$ put it at $x-y=0$. Write
\[
 f_2=Ax^2+Bxy+Cy^2,\qquad
 g_2=Dx^2+Exy+Fy^2.
\]
The quartic part of~\eqref{l4v3:double-residual} is
\begin{align}\label{l4v3:reduced-quartic}
 P_4={}&Ac^2x^4-2Abcx^3y
 +(Db^2-Bbc+Ebc-Cc^2)x^2y^2\notag\\
 &\quad+2Fbcxy^3-Fb^2y^4.
\end{align}
Higher jets of $f,g$ cannot contribute to this degree.
If $P_4=u(x-y)^4$ with $u\ne0$, the coefficients of $x^4,x^3y$
give
\[
 Ac^2=u,\qquad -2Abc=-4u,
 \qquad \frac bc=2.
\]
The coefficients of $y^4,xy^3$ give
\[
 -Fb^2=u,\qquad 2Fbc=-4u,
 \qquad \frac cb=2.
\]
Multiplying the two ratios gives $1=4$, a contradiction over
$\mathbb C$.
\end{proof}

\begin{lemma}\label{l4v3:no-nonreduced-double}
A nonreduced logarithmic zero at an old ordinary double point cannot
occur when $R=L^4$.
\end{lemma}
\begin{proof}
Curvilinearity implies that $f_1=A\ell$, $g_1=B\ell$ for a linear
form $\ell$, where $(A,B)\ne(0,0)$. Substitution in
\eqref{l4v3:double-residual} gives
\[
 P_3=AB(B-A)\ell^3.
\]
Thus either $A=B\ne0$ or exactly one of $A,B$ vanishes.

If $A=B\ne0$, put $a=(f+g)/2$, $b=(f-g)/2$ and use the
Euler--Hamiltonian frame for $q=xy$. Here $a$ is smooth and
$b\in\mathfrak m_p^2$. Its quartic residual is $-4a_1^2b_2$.
It follows that $a_1$ is proportional to the equation of $L$ and
$b_2$ to its square. On the smooth curve $a=0$, let $r$ be the
vanishing order of that linear equation and let
$n=\operatorname{length}Z_p$. The restriction of the exact vertex
formula~\eqref{l4v3:vertex-formula} below (which holds for every $m\ge2$, \cite[Lemma~4.5]{GST16}; here
$k=-2$) is
\[
 b^2\bigl(2\eta(a)-kb\bigr).
\]
The function $\eta(a)$ has exact order one, since $L$ is new,
whereas $b$ has order at least two. Comparison with a unit times
the fourth power gives $4r=2n+1$, impossible. If the curve $a=0$
were exactly $L$, its residual restriction would be nonzero by the
same lowest-order comparison, again a contradiction.

If, say, $A\ne0$, $B=0$, the quartic part is divisible by $\ell$;
its being a nonzero fourth power therefore forces $\ell\sim L$.
Normalize $L=(x-y=0)$, and write $g_2=ux^2+vxy+wy^2$.
Direct expansion gives
\[
 P_4=A^2(x-y)y\bigl(-(u+v)x^2-2wxy+wy^2\bigr).
\]
It has the old factor $y$, so cannot be a nonzero multiple of
$(x-y)^4$. The other case is symmetric.
\end{proof}

\paragraph{Vertices of higher multiplicity}

At an old vertex of multiplicity $m\ge3$, let $q$ be the homogeneous
product of the old local equations, and write
\[
 E=x\partial_x+y\partial_y,
 \qquad \eta=q_y\partial_x-q_x\partial_y,
 \qquad X=aE+b\eta.
\]
The logarithmic ideal is $(a,b)$ and the exact vertex formula is
\begin{equation}\label{l4v3:vertex-formula}
 P=m\bigl(bX(a)-aX(b)\bigr)-m(m-1)a^2b-kb^3,
 \qquad \operatorname{Hess}_q(\eta,\eta)=qk,
 \quad\deg k=2m-4.
\end{equation}

\begin{lemma}\label{l4v3:no-higher-vertex}
No old vertex of multiplicity at least three belongs to $Z$ when
$R=L^4$.
\end{lemma}
\begin{proof}
The cubic part of~\eqref{l4v3:vertex-formula} is
$-m(m-1)a_1^2b_1$. If $Z_p$ were reduced, $a_1,b_1$ would be
independent and this polynomial would be nonzero. Thus suppose
$n=\operatorname{length}Z_p\ge2$. Curvilinearity says that at least
one of $a,b$ is smooth. Vanishing of the cubic gives two cases.

If $b$ is smooth and $a\in\mathfrak m_p^2$, the Euler terms have
order at least five, the terms $b^2\eta(a)$ and $ab\eta(b)$ have
order at least $m+2$, and $kb^3$ has order at least $2m-1$.
This contradicts residual order four.

Therefore $a$ is smooth and $b\in\mathfrak m_p^2$. Its quartic
part is
\[
 P_4=-m^2a_1^2b_2.
\]
Put $L=(x=0)$. We obtain $a_1\sim x$ and $b_2\sim x^2$.
In particular $\beta=\operatorname{ord}b(0,y)\ge3$; this order
is finite, because $L$ is non-invariant. The normal restriction to
$L$ is $b(0,y)q_y(0,y)$, and $q(0,1)\ne0$, so its order is
$\beta+m-1$. Its global degree is seven. Consequently
\begin{equation}\label{l4v3:normal}
 3\le\beta\le8-m,\qquad m\le5.
\end{equation}

Write the curve $a=0$ as $x=\varphi(y)$ and set
$r=\operatorname{ord}\varphi\ge2$. The case $\varphi\equiv0$
is impossible: on $L$ the residual is
$b^2(m\eta(a)-kb)$, whose first term has strictly lowest order
and is nonzero. Thus $r$ is finite. Since $b_2\sim x^2$,
$b|_{a=0}$ has order $n\ge3$. On this curve $\eta(a)$ has exact
order $m-1$, while $kb$ has larger order. Restricting the residual
therefore gives
\begin{equation}\label{l4v3:restriction-order}
 4r=2n+m-1.
\end{equation}
It follows that $m$ is odd, so only $m=3$ or $m=5$ remains.

Because $b_2\sim x^2$, one also has
\begin{equation}\label{l4v3:beta}
 \operatorname{ord}b_x(0,y)\ge2,
 \qquad \beta\ge\min\{n,r+2\}.
\end{equation}
The latter inequality follows by Taylor expansion from $x=\varphi(y)$
to $x=0$; its first correction has order at least $r+2$, and the
second has order at least $2r\ge r+2$.

For $m=5$, equation~\eqref{l4v3:restriction-order} gives
$n=2r-2\ge3$, hence $r\ge3$. Formula~\eqref{l4v3:beta} gives
$\beta\ge4$, contradicting~\eqref{l4v3:normal}.

For $m=3$, the same equation gives $n=2r-1$. Combining
\eqref{l4v3:normal} and~\eqref{l4v3:beta}, the only possibilities
for $(r,\beta)$ are
\begin{equation}\label{l4v3:four-pairs}
 (2,3),\quad(2,4),\quad(2,5),\quad(3,5).
\end{equation}
We now use the exact residual on $L$, not just its order on the
curve $a=0$. Write
\[
 a(0,y)=A y^r+O(y^{r+1}),\qquad
 b(0,y)=B y^\beta+O(y^{\beta+1}),\qquad AB\ne0.
\]
The Euler terms in~\eqref{l4v3:vertex-formula}, restricted to $L$,
are
\[
 3a^2b\left(y\frac{a'}a-y\frac{b'}b-2\right).
\]
Their order is $2r+\beta$ and their leading coefficient is
$3A^2B(r-\beta-2)\ne0$. On the other hand,
\[
 \operatorname{ord}(\eta(a)|_L)=2,\qquad
 \operatorname{ord}(\eta(b)|_L)\ge\min\{4,\beta+1\}.
\]
The second inequality uses~\eqref{l4v3:beta} and the fact that both
coefficients of $\eta|_L$ have order at least two. Hence the three
remaining contributions have order at least
\[
 2\beta+2,\qquad
 r+\beta+\min\{4,\beta+1\},\qquad 3\beta+2.
\]
For every pair in~\eqref{l4v3:four-pairs}, all three numbers are
strictly greater than $2r+\beta$; explicitly their lower bounds are
$(8,9,11)$, $(10,10,14)$, $(12,11,17)$ and $(12,12,17)$, whereas
the Euler orders are $7,8,9,11$. The nonzero Euler leading term
therefore cannot cancel. This contradicts $P|_L=0$.
\end{proof}

\begin{theorem}\label{l4v3:no-positive-fourth-power}
In the critical seventeen-line setting, a fourth-power residual
$R=L^4$ forces $Z=\varnothing$, equivalently $\nu=0$.
In particular this factorization is impossible at positive defect.
\end{theorem}
\begin{proof}
The support of $Z$ consists of old vertices, as noted at the beginning of
Section~\ref{l4v3:section}.  Lemmas~\ref{l4v3:no-reduced-double} and
\ref{l4v3:no-nonreduced-double} exclude every ordinary double
vertex, and Lemma~\ref{l4v3:no-higher-vertex} excludes every vertex
of higher multiplicity. Hence the support is empty.
\end{proof}

\begin{remark}
The displayed initial forms of this section were rechecked by exact symbolic computation
(Section~\ref{sec:verification}).  The proof does not depend on any incidence enumeration.
\end{remark}

\subsection{A cubic--line residual}\label{sec:L3M}
In this section $R=L^3M$, where $L,M$ are distinct new lines and $L$ is not invariant
(Corollary~\ref{cor:nonreduced}); put $p_0=L\cap M$.  We prove $\nu\le3$
(Theorem~\ref{cg4:cubic-three}).  The first three subsections study a zero at $p_0$, according to
its position with respect to the old arrangement.

\subsubsection{The crossing point off the arrangement}
\label{ae:section}

The following local calculation concerns the physical inflection determinant
of a vector field. It is independent of any classification of arrangement
lattices.

\begin{lemma}\label{ae:cubic-crossing}
Let $X$ be a holomorphic vector field with an isolated curvilinear zero at
$p$ in a smooth complex surface. Suppose, in affine linear coordinates,
that its inflection determinant is
\[
 D=U\ell^3m,
\]
where $U$ is a unit, $\ell,m$ are independent linear forms vanishing at
$p$, and $L=(\ell=0)$ is non-invariant. Then the local zero scheme has
length at most five.  If its length is at least two and its tangent is $L$, then the linearization
is nilpotent, the length is five, and the trace on $L$ has length three.
\end{lemma}
\begin{proof}
Write $X=A\partial_x+B\partial_y$ and $I=(A,B)$. Thus
\begin{equation}\label{ae:inflection}
 D=A^2B_x+AB B_y-AB A_x-B^2A_y.
\end{equation}
Suppose first that the zero has length one. Its linearization $J$ is
invertible, and the quadratic part of $D$ is
$\det(J)\det(w,Jw)$, where $w=(x,y)^{\mathsf t}$. Since $D$ has order
four, $J=\lambda\operatorname{Id}$ with $\lambda\ne0$. Put $L=(x=0)$.
If $A(0,y)=a_ky^k+\cdots$ is nonzero, then $k\ge2$, while
$B(0,y)=\lambda y+\cdots$. The leading term of $D(0,y)$ is
$-k a_k\lambda^2y^{k+1}$, contradicting $D|_L=0$. Otherwise $L$ is
invariant, also a contradiction. Consequently the length is at least
two. Curvilinearity now implies that $J$ has rank one.

The kernel line of $J$ divides the fourth-order initial form of $D$.
Here is a direct check, which also covers the nilpotent case. If $J$ is
semisimple, rescale $X$ and choose coordinates with
$A=x+O(2)$, $B=O(2)$. Vanishing of the cubic part of $D$ forces
$B_2=b_{11}xy$. If
\[
 A_2=a_{20}x^2+a_{11}xy+a_{02}y^2,\qquad
 B_3=b_{30}x^3+b_{21}x^2y+b_{12}xy^2+b_{03}y^3,
\]
then
\[
 D_4=x\bigl((a_{02}b_{11}-b_{03})y^3
 +(-a_{20}b_{11}+b_{11}^2+b_{21})x^2y+2b_{30}x^3\bigr).
\]
If $J$ is nilpotent, choose $A=O(2)$, $B=x+O(2)$. Vanishing of the
cubic part forces $A_2=a_{20}x^2$, and, writing
\[
 B_2=b_{20}x^2+b_{11}xy+b_{02}y^2,\qquad
 A_3=a_{30}x^3+a_{21}x^2y+a_{12}xy^2+a_{03}y^3,
\]
one obtains
\[
 D_4=-x^2\bigl(3a_{03}y^2+(2a_{12}-2a_{20}b_{02})xy
 +(a_{20}^2-a_{20}b_{11}+a_{21})x^2\bigr).
\]
Thus the kernel, which is also the curvilinear tangent, is either $L$
or $M=(m=0)$.

If the tangent is $M$, use $L=(x=0)$ and $M=(y=0)$ and write
$I=(y-\varphi(x),x^n)$ with $\operatorname{ord}\varphi\ge2$.
The determinant~\eqref{ae:inflection} belongs to $I^2$. In the coordinate
$w=y-\varphi(x)$, the coefficient of $w$ in
\[
 U(x,w+\varphi(x))x^3(w+\varphi(x))
\]
is
$x^3\bigl(U(x,\varphi(x))+\varphi(x)U_y(x,\varphi(x))\bigr)$.
The factor in parentheses is a unit. Membership in
$(w^2,wx^n,x^{2n})$ therefore forces $n\le3$.

It remains that the tangent is $L$. Put $L=(x=0)$; the equation of
$M$ then has nonzero $y$ coefficient. Set
\[
 a=A(0,y),\quad b=B(0,y),\quad
 c=A_x(0,y),\quad d=B_x(0,y).
\]
The identity $D|_L=0$ is
\begin{equation}\label{ae:line-identity}
 a^2d-acb+abb'-b^2a'=0.
\end{equation}
Non-invariance says $a\not\equiv0$.

For semisimple $J$, normalize $A=x+O(2)$, $B=O(2)$. Thus $c$ is a
unit, $d(0)=0$, and $k=\operatorname{ord}a\ge2$. If $b\not\equiv0$,
put $v=b/a$. Equation~\eqref{ae:line-identity} gives
$d=cv-avv'$. If $\operatorname{ord}b<k$, the term $cv$ has a pole
which cannot be canceled by $avv'$, whose valuation is strictly larger.
If $\operatorname{ord}b=k$, the equation gives $d(0)\ne0$.
Hence $\operatorname{ord}b\ge k+1$ and $v$ vanishes at the origin.
On the smooth curve $A=0$, write
$x=\phi(y)=-a/c+O(a^2)$, of order $k$. Taylor expansion gives
\[
 B(\phi(y),y)=b-da/c+O(a^2)
 =a^2vv'/c+O(a^2).
\]
Thus its order, equal to $n=\operatorname{length}\mathcal O/I$, is at
least $2k$. If $b\equiv0$, equation~\eqref{ae:line-identity} gives
$d\equiv0$ and the same estimate follows.
On the other hand, restriction of $D=U\ell^3m$ to $A=0$ gives
\[
 -B^2A_y=U\phi^3m(\phi(y),y).
\]
The orders on the two sides are $2n+k-1$ and $3k+1$, respectively:
$A_y|_{A=0}=-A_x|_{A=0}\phi'$ has order $k-1$, and transversality
of $M$ and $L$ gives $\operatorname{ord}m(\phi(y),y)=1$.
Consequently $n=k+1<2k$, a contradiction.

For nilpotent $J$, choose $A=O(2)$, $B=x+O(2)$. Here $d$ is a
unit and $c(0)=0$. The function $b$ cannot vanish identically, since
\eqref{ae:line-identity} would then give $a^2d=0$. Put
$k=\operatorname{ord}b\ge2$, $h=\operatorname{ord}a\ge2$ and $r=a/b$.
The same equation gives
\[
 c=rd-b\frac{r'}r.
\]
If $h<k$, the first term has a pole, while the second has order $k-1$;
if $h=k$, its constant term contradicts $c(0)=0$. Thus $h>k$ and
$r'/r=(h-k)/y+O(1)$.
On the smooth curve $B=0$, write $x=\psi(y)=-b/d+O(b^2)$, of order
$k$. Then
\[
 A(\psi(y),y)=a-cb/d+O(b^2)
 =\frac{b^2}{d}\frac{r'}r+O(b^2).
\]
Its leading term has exact order $2k-1$, so $n=2k-1$. Restricting the
inflection identity to $B=0$ gives
\[
 A^2B_x=U\psi^3m(\psi(y),y),
\]
whence $2n=3k+1$. Thus $k=3$ and $n=5$. Finally the trace on $L$
has length $\min\{h,k\}=3$. This completes the proof.
\end{proof}

\subsubsection{The crossing point on a smooth old line}
\begin{proposition}[No smooth old zero at a triple--simple crossing]
\label{prop:L3M-no-smooth-crossing}
Let $X=f\partial_x+yg\partial_y$ be a holomorphic vector-field germ,
where $f,g\in\mathfrak m=(x,y)$ and $I=(f,g)$ is curvilinear of finite
colength $n>0$. Suppose that $L=(x=0)$ is not invariant and that
\[
 P=f^2g_x+fg^2+fyg g_y-fg f_x-yg^2f_y
   =Ux^3(ax+by),\qquad U(0)ab\ne0.
\]
Then no such germ exists.
\end{proposition}
\begin{proof}
The assumption $ab\ne0$ says that the second residual line is distinct
from both $L$ and the old line $(y=0)$. Non-invariance means
$f(0,y)\not\equiv0$. Multiplying $X$ by a holomorphic unit multiplies
$P$ by its cube, so it preserves every hypothesis.

First $n\ne1$: independent linear parts $f_1=Ax+By$ and
$g_1=Cx+Dy$ give
$P_2=-(AD-BC)yf_1\ne0$, whereas the right side has order four.
Thus $n\ge2$.

Suppose that $f$ is smooth. If $f_x(0)=0$, its zero curve is
$y=\phi(x)$ of order $k\ge2$, or is the old line $y=0$.
The latter is impossible by restriction of the identity to $f=0$.
In the former case, the same restriction gives
$k+2n=4$, again impossible. Hence $f_x(0)\ne0$.
After division by a unit write $f=x-\phi(y)$, where
$k=\operatorname{ord}\phi$ is finite by non-invariance. Restriction
to $f=0$ gives
\[
 y\phi' g(\phi(y),y)^2
   =U(\phi(y),y)\phi(y)^3(a\phi(y)+by).
\]
If $k\ge2$, comparison of orders gives $2n=2k+1$, impossible.
If $k=1$ and this curve is not tangent to $ax+by=0$, it gives
$2n=3$, also impossible. Scale $y$ so that the second residual line
is $x-y=0$. The identity $\phi\equiv y$ is impossible, since restriction
to $f=0$ would then give $yg(\phi(y),y)^2\equiv0$, contrary to finite
colength. Necessarily
\[
 \phi=y+O(y^s),\qquad s=2n-2\ge2.
\]
Write temporarily $\phi=y+q y^2+O(y^3)$. Since $n\ge2$,
$g_1=c(x-y)$ for some constant $c$. If
$g_2=A x^2+Bxy+C y^2$, direct expansion gives
\[
 P_3=(x-y)\bigl((A+c^2)x^2-2(A+c^2)xy
       +(-B-C+c^2-cq)y^2\bigr).
\]
Thus $A=-c^2$ and $C=-B+c^2-cq$. In particular, substituting
$x=\phi(y)$ cancels the degree-two term of $g$, so $n\ge3$.
Consequently $s\ge4$, $q=0$, and $f=x-y+O(\mathfrak m^4)$.
Write
\[
 g=c(x-y)-c^2x^2+Bxy+(c^2-B)y^2
       +A_3x^3+B_3x^2y+C_3xy^2+D_3y^3+O(\mathfrak m^4).
\]
The quartic term is $(x-y)Q_3$, where the coefficients of
$x^3,x^2y,xy^2$ in $Q_3$ are respectively
\[
 2A_3-2c^3,\qquad
 B_3-3A_3+3Bc+c^3,\qquad
 -2B_3-6Bc+4c^3.
\]
But $P_4=u x^3(x-y)$, $u\ne0$. The last two displayed coefficients
must vanish, giving $B_3=-3Bc+2c^3$ and $A_3=c^3$.
The first coefficient is then zero, contradicting $u\ne0$.
This excludes every smooth $f$.

It follows that $f\in\mathfrak m^2$. Curvilinearity forces $g$ to
be smooth. If $g_x(0)=0$, its zero curve is $y=\psi(x)$ of order
$k\ge2$, or is identically $y=0$. On $g=0$ the identity is
$f^2g_x=Ux^3(ax+by)$. The identically zero case is impossible;
otherwise comparison of orders gives $2n+k-1=4$, also impossible.
Hence $g_x(0)\ne0$. Divide by a unit to normalize
$g=x-\psi(y)$, where $\psi\not\equiv0$: if $\psi\equiv0$,
restriction to $g=0$ would force $f|_L=0$.

If $\operatorname{ord}\psi=1$, scale $y$ so that $g_1=x-y$.
Put $f_2=A x^2+Bxy+C y^2$. The coefficients of $P_4$ in
descending powers of $x$ are
\[
 \begin{split}
 &A(1-A),\\
 &A(2A-B-3),\\
 &3AB+2A-B-C,\\
 &2AC+B^2+BC+B+C,\\
 &C(B+C).
 \end{split}
\]
Since $x^3$ divides $P_4$, the last three expressions vanish.
If $C=0$, the fourth gives $B(B+1)=0$; the third then gives
$A=0$ when $B=0$, and $A=1$ when $B=-1$.
If $C\ne0$, the fifth gives $B=-C$, the fourth gives $2AC=0$,
and hence $A=0$. In every case $A(1-A)=0$, contradicting the
nonzero $x^4$ coefficient of $Ux^3(ax+by)$.

It remains to consider $k=\operatorname{ord}\psi\ge2$.
Restriction to $g=0$ gives
\[
 f(\psi(y),y)^2=U(\psi(y),y)\psi(y)^3(a\psi(y)+by),
 \qquad 2n=3k+1.
\]
Thus $k$ is odd and $k\ge3$. Differentiating the exact formula for
$P$ and then restricting to $g=0$ gives
\[
 P_x|_{g=0}=f(f_x-y\psi')|_{g=0}.
\]
The left side has exact order $2k+1$, because its leading contribution
is $3bU\psi^2y$. Therefore, putting $\ell=(k+1)/2<k$, one has
\[
 \operatorname{ord}(f_x-y\psi')|_{g=0}=\ell,
 \qquad \operatorname{ord}f_x(0,y)=\ell\ge2.
\]
The second assertion follows since $y\psi'$ has order $k$ and
Taylor expansion from $x=\psi(y)$ to $x=0$ changes $f_x$ by
a term of order at least $k$.

The normal order can be determined exactly.
Set $\alpha=f(0,y)$, $\beta=g(0,y)=-\psi(y)$,
$c=f_x(0,y)$, and $r=\alpha/\beta$. Restriction of $P$ to $L$ gives
the exact meromorphic identity
\[
 c=r+\beta\left(1-y\frac{r'}r\right).
\]
If $h=\operatorname{ord}\alpha\le k$, the nonpositive-order term
$r$ cannot be canceled by the second term, of order at least $k$;
this contradicts $c(0)=0$. Hence $h>k$. Since
$\operatorname{ord}c=\ell<k$, the same identity gives
$\operatorname{ord}r=\ell$ and
\[
 h=k+\ell=n,\qquad
 \operatorname{length}\mathcal O/(I,x)=k.
\]
In particular $h\ge5$ and $\operatorname{ord}c\ge2$, so
$f_2=A x^2$. Now $g_1=x$ and direct expansion gives
$P_4=A(1-A)x^4$, which has zero $x^3y$ coefficient.
This contradicts $bU(0)\ne0$ and completes the proof.
\end{proof}

\subsubsection{The crossing point at an old double point}
\label{dc:section}

\begin{proposition}\label{dc:node-crossing}
Let $p=L\cap M$ be an ordinary double point of the old arrangement in
the critical seventeen-line setting. Suppose $R=L^3M$, with $L,M$ new
and distinct and $L$ non-invariant. If the logarithmic zero scheme is
nonreduced at $p$, then its tangent is $L$. By the trace bound
$\operatorname{length}(Z\cap L)\le4$ of Theorem~\ref{thm:trace}, its local length and trace
are necessarily
\[
 \operatorname{length}Z_p=3,\qquad
 \operatorname{length}(Z_p\cap L)=2.
\]
In particular, a local length-six zero of trace four is impossible.
No assertion is made here about a reduced zero at $p$.
\end{proposition}
\begin{proof}
Take the two old lines to be $x=0,y=0$, normalize $L=(x-y=0)$, and
write $X=xf\partial_x+yg\partial_y$. Its exact physical inflection
residual is
\begin{equation}\label{dc:physical}
 P=xf^2g_x+fg(g-f)+yfg g_y-xfg f_x-yg^2f_y.
\end{equation}
It equals a unit times $(x-y)^3m(x,y)$, where $m=0$ defines $M$.
The logarithmic zero ideal is $(f,g)$. Curvilinearity and
nonreducedness mean that the linear parts have the form
$f_1=A K$, $g_1=B K$, with $K$ the tangent equation and at least one
of $A,B$ nonzero. The cubic part of~\eqref{dc:physical} is
$AB(B-A)K^3$. Therefore either $A=B\ne0$, or exactly one is zero.

The latter cases are impossible. For example, put
$K=\alpha x+\beta y$, $f_1=A K$ and
$g_2=ux^2+vxy+wy^2$. Direct expansion gives
\[
 P_4=A^2K y\bigl((\beta u-\alpha v)x^2
                       -2\alpha wxy-\beta wy^2\bigr),
\]
which has the old factor $y$. The case with $g_1\ne0,f_1=0$ has
the old factor $x$. Neither can equal the nonzero fourth-order form
of $L^3M$, which has no old factor.

Thus the Euler and Hamiltonian coefficients
$a=(f+g)/2$, $b=(f-g)/2$ satisfy $a_1=A K$ and $b_1=0$.
Their quartic residual is exactly $-4a_1^2b_2$. Consequently $K^2$
divides $L^3M$, forcing $K=L$. After a constant rescaling of the field
we may assume $a_1=x-y$.

Use coordinates $u=x-y$, $v=y$. The logarithmic frame is
\[
 E=u\partial_u+v\partial_v,\qquad
 H=(u+2v)\partial_u-v\partial_v,
\qquad X=aE+bH.
\]
Write $M=(u+\kappa v=0)$, where $\kappa\ne0$ because $M\ne L$.
The preceding initial-form calculation gives
\begin{equation}\label{dc:initial}
 a_1=u,\qquad b_2=c u(u+\kappa v),\qquad c\ne0.
\end{equation}
In this frame the exact residual identity is
\begin{equation}\label{dc:frame}
 P=2\bigl(bX(a)-aX(b)-a^2b+b^3\bigr).
\end{equation}

Let $a=0$ have equation $u=\psi(v)$, and put $k=\operatorname{ord}\psi$.
The curve cannot be $L$ itself: on $a=0$, formula~\eqref{dc:frame}
is $2b^2(H(a)+b)$, and $H(a)$ has exact order one there, whereas
$b$ has order $n=\operatorname{length}Z_p\ge3$. The latter inequality
follows from~\eqref{dc:initial}, since restricting $b_2$ to a smooth
curve tangent to $u=0$ has order at least three.
Thus $k$ is finite and $k\ge2$. Comparing orders on $a=0$ gives
\[
 3k+1=2n+1,\qquad 3k=2n.
\]
The trace on $L$ is $t=\min(k,n)$; hence $n>k$ and $t=k$.
The trace bound $t\le4$ leaves only $(n,t)=(3,2)$ or $(6,4)$.

We exclude the latter by restricting the exact identity to $L$.
Suppose $n=6,t=4$. Divide analytically by the smooth function $a$:
\begin{equation}\label{dc:division}
 b=a h+B(v),\qquad
 \operatorname{ord}B=6.
\end{equation}
Here $B(v)=b(\psi(v),v)$, and $h$ is holomorphic. Equation
\eqref{dc:initial} implies $h_1=c(u+\kappa v)$.
On $L$, therefore, $a_0=a(0,v)$ has order four, $h_0$ has exact
order one, and $b_0=a_0h_0+B$ has exact order five.
Moreover
\[
 b a_u-a b_u=B a_u-a^2h_u
\]
has exact order six on $L$, since $a_u$ is a unit there and the
second term has order at least eight.

Restricting~\eqref{dc:frame} to $u=0$ now gives
\begin{equation}\label{dc:line}
 \frac{P|_L}{2}=
 2v b_0(b a_u-a b_u)|_L
 +v(a_0-b_0)(b_0a_0'-a_0b_0')
 -b_0(a_0^2-b_0^2).
\end{equation}
The first term has exact order $1+5+6=12$ and nonzero leading
coefficient. Each of the other two has order at least thirteen:
$a_0,b_0,a_0',b_0'$ have orders $4,5,3,4$, respectively.
Thus $P|_L$ is nonzero, contradicting $L\mid R$. This excludes
$(n,t)=(6,4)$ and proves the proposition.
\end{proof}

\subsubsection{Points away from the crossing; defect six}
\label{cg:section}

Write $n_p=\operatorname{length}Z_p$, $t_p=\operatorname{length}(Z_p\cap L)$ and
$T=\sum_p t_p$; thus $T\le4$ by Theorem~\ref{thm:trace}.  Since
$\operatorname{Supp}Z\subseteq\operatorname{Sing}R=L$, every point of $Z$ lies on $L$.
By the cubic local tests of Section~\ref{sec:cubictests}, a zero away from $p_0$ lies at an old vertex, not at a smooth old point or off
the arrangement.  At such a vertex the residual is a unit times $\ell^3$, $\ell$ an equation of $L$,
so Lemma~\ref{lem:vertexineq} gives the \emph{strict vertex inequality}
\begin{equation}\label{cg:strict}
 2n_p+1\le3t_p,\qquad\text{i.e.}\qquad n_p\le\Bigl\lfloor\frac{3t_p-1}{2}\Bigr\rfloor .
\end{equation}

\paragraph{The normal budget.}  The normal restriction of the foliation to $L$ is a nonzero section of
$\cO_L(7)$ (Section~\ref{sec:foliation}).  Its vanishing order at a point $p\in L$ is
called the \emph{normal cost} of $p$; the costs add up to at most $7$.  A point of $Z$ has cost at least
$t_p\ge1$.  At an old vertex $p$ of multiplicity $m$, in the frame $X=aE+b\eta$ of~\eqref{ar:vertex}, the
normal restriction is $b|_L$ times the restriction of $\eta(\ell)$, where $\ell$ is a local linear equation of
$L$ (the term $aE(\ell)=a\ell$ vanishes on $L$); since $L$ is not an old line, $\eta(\ell)|_L$ is a nonzero
constant times $x^{m-1}$ in a coordinate $x$ on $L$ centred at $p$, so the cost is $(m-1)+\operatorname{ord}(b|_L)$;
when $p\in Z$ this is at least $(m-1)+t_p$, since $b\in I_{Z,p}$.

The following arguments retain reduced old double zeros until their full cubic initial form has been
examined.

\begin{lemma}[No reduced zero away from the crossing]
\label{cg:no-reduced-off}
If $p\in\operatorname{Supp}Z$ and $p\ne p_0$, then $Z_p$ is not reduced.
\end{lemma}
\begin{proof}
At an old vertex of multiplicity $m\ge3$, use the logarithmic frame
$X=aE+b\eta$. If $Z_p$ were reduced, $a_1,b_1$ would be independent
linear forms. The cubic initial form of the residual is
$-m(m-1)a_1^2b_1$. It cannot be a nonzero cube, whereas the local
residual is a unit times the cube of the equation of $L$.

At an old double vertex use the old axes and write
$X=xf\partial_x+yg\partial_y$, with
$f_1=ax+by$, $g_1=cx+dy$ and $ad-bc\ne0$. The complete cubic
residual, including its derivative terms, is
\[
 P_3=-(f_1-g_1)(acx^2+2adxy+bdy^2).
\]
Rescale the axes so that $L=(x-y=0)$. If $P_3=u(x-y)^3$, $u\ne0$,
the quadratic factor must be a nonzero multiple of $(x-y)^2$.
Its coefficients give $ac=bd=\lambda$ and $ad=-\lambda$, with
$\lambda\ne0$. Hence $d=-c$, $b=-a$, making $f_1,g_1$
proportional and contradicting $ad-bc\ne0$.
\end{proof}

\begin{lemma}[No old double zero on a pure-cubic residual]
\label{cg:no-node-off}
An old ordinary double vertex cannot support a logarithmic zero at a
point of $L\setminus M$.
\end{lemma}
\begin{proof}
The reduced case is Lemma~\ref{cg:no-reduced-off}. Suppose $n\ge2$.
Put $u=x-y$, $v=y$, so the old product is $v(u+v)$ and
\[
 E=u\partial_u+v\partial_v,\qquad
 H=(u+2v)\partial_u-v\partial_v,\qquad X=aE+bH.
\]
Let $K$ define the curvilinear tangent. Write $a_1=A K$, $b_1=B K$.
The cubic residual is $2B(B^2-A^2)K^3$. Since it is a nonzero
multiple of $u^3$, we have $K=u$, $B\ne0$ and $B^2\ne A^2$.
The exact residual divided by two is
\[
 W=bX(a)-aX(b)-a^2b+b^3.
\]
On $L=(u=0)$ it becomes
\begin{equation}\label{cg:node-line}
 W|_L=2vb(ba_u-ab_u)+v(a-b)(ba_v-ab_v)-b(a^2-b^2),
\end{equation}
and it must vanish identically.

If $A\ne0$, divide by the smooth function $a$ to write
$b=ah+B_0(v)$, where $\operatorname{ord}B_0=n$ and
$h(0)=B/A\ne0,\pm1$. Write $a=0$ as $u=\psi(v)$ and let
$k=\operatorname{ord}\psi$. This curve cannot be exactly $L$:
on $a=0$, $W=b^2(H(a)+b)$ and $H(a)$ has exact order one,
while $b$ has order $n\ge2$. Thus $k$ is finite. Restricting
$W=Uu^3$ to this curve gives $3k=2n+1$.
Consequently $n>k=t\ge2$, where $t$ is the trace on $L$.
On $L$, both $a,b$ have order $t$. The identity
\[
 ba_u-ab_u=B_0a_u-a^2h_u
\]
has exact order $n<2t$, so the first summand of
\eqref{cg:node-line} has exact order $1+t+n$.
The second summand has order at least $3t+1$, because
\[
 ba_v-ab_v=B_0a_v-a^2h_v-aB_0'
\]
has order at least $2t$. The last summand has exact order $3t$.
But $1+t+n=(5t+1)/2<3t$, contradicting the identity.

If $A=0$, then $b$ is smooth and $a=bh+A_0(v)$, with
$h(0)=0$ and $\operatorname{ord}A_0=n$. Let $b=0$ have equation
$u=\psi(v)$. The curve cannot be $L$, since that would make $L$
invariant. If $k=\operatorname{ord}\psi$, restriction gives
$3k=2n+k$, because $W|_{b=0}=-a^2E(b)$ and
$E(b)|_{b=0}$ has exact order $k\ge2$. Hence $n=k=t$.
On $L$, both $a,b$ have order $t$, and
$ba_u-ab_u=b^2h_u-A_0b_u$ has exact order $t$.
The first summand of~\eqref{cg:node-line} therefore has exact
order $2t+1$. The leading terms of $bA_0'-A_0b_v$ cancel, so
$ba_v-ab_v=b^2h_v+bA_0'-A_0b_v$ has order at least $2t$.
The second summand has order at least $3t+1$ and the last at least
$3t$. Since $2t+1<3t$, this is again impossible.
\end{proof}

\begin{lemma}[The order identity at the remaining old vertices]
\label{cg:vertex-order}
Let $p\in L\setminus M$ be an old vertex of multiplicity $m\ge3$
supporting a nonreduced zero. Let $a$ be the Euler coefficient in the
frame $X=aE+b\eta$ of~\eqref{ar:vertex} (it is smooth, see the proof),
and let $k$ be the contact order of $a=0$ with $L$.
Then
\[
 3k=2n_p+m-1,\qquad t_p=\min\{k,n_p\}.
\]
\end{lemma}
\begin{proof}
The cubic residual is $-m(m-1)a_1^2b_1$ and is nonzero, so both
$a,b$ are smooth with common tangent $L$. On $a=0$ the residual is
\[
 b^2(m\eta(a)-\kappa b),\qquad\deg\kappa=2m-4.
\]
Because $L$ is new, the leading homogeneous term of $\eta(a)$
has exact order $m-1$ along that curve. The term $\kappa b$
has strictly larger order. The curve cannot coincide with $L$,
for then the left side would be nonzero while the residual vanishes.
Comparison with a unit times the cube of the equation of $L$ gives
the identity. The formula for the trace follows by restricting
the two generators of the curvilinear ideal to $L$.
\end{proof}

\begin{lemma}[The simple-factor tangent costs at least four normal zeros]
\label{cg:off-cross-cost}
Suppose $p_0$ is off the old arrangement and its physical zero scheme
is nonreduced with tangent $M$. Then the normal restriction to $L$
vanishes to order at least four at $p_0$.
\end{lemma}
\begin{proof}
The linearization has rank one and kernel $M$. In the nilpotent
case, choose $A=O(2)$ and $B=x+O(2)$ with $M=(x=0)$.
Vanishing of the cubic inflection term forces $A_2=a_{20}x^2$,
and its quartic term is divisible by $x^2$. This contradicts the
simple factor $M$ in the initial form of $D=U L^3M$.

In the semisimple case choose $A=x+O(2)$, $B=O(2)$, again with
$M=(x=0)$. Cubic vanishing forces $B_2=b_{11}xy$.
Writing $A_2=a_{20}x^2+a_{11}xy+a_{02}y^2$ and
$B_3=b_{30}x^3+b_{21}x^2y+b_{12}xy^2+b_{03}y^3$, the quartic
initial form is
\[
 D_4=x\bigl((a_{02}b_{11}-b_{03})y^3+
 (-a_{20}b_{11}+b_{11}^2+b_{21})x^2y+2b_{30}x^3\bigr).
\]
It has no $x^2y^2$ term. Write $L=(\alpha x+\beta y=0)$,
where $\beta\ne0$. Comparing with $u x(\alpha x+\beta y)^3$
forces $\alpha=0$. Thus $L=(y=0)$, the nonzero eigendirection.
The target initial form has no $x^4$ coefficient either, so
$b_{30}=0$. Since $B_2=b_{11}xy$, the normal coefficient
$B(x,0)$ has order at least four.
\end{proof}

\begin{lemma}[The all-support crossing resonance at higher vertices]
\label{cg:no-six-at-high-crossing}
An old vertex of multiplicity at least three cannot support a
length-six zero with trace four on $L$ at $p_0$.
\end{lemma}
\begin{proof}
By Lemma~\ref{cg:crossing-traces} below, such a zero has tangent $L$. Write $a_1=A\ell$, $b_1=B\ell$.
Since the residual has order four, its cubic term gives $A^2B=0$.
If $A=0$, then $b$ is smooth, $a\in\mathfrak m^2$, and every
residual term has order at least five, a contradiction. Thus $a$
is smooth, $b\in\mathfrak m^2$, and the quartic term is
$-m^2a_1^2b_2$. Along $a=0$, the same exact restriction as in
Lemma~\ref{cg:vertex-order} gives
$3\cdot4+1=2\cdot6+m-1$, forcing $m=2$.
\end{proof}

\begin{lemma}[Numerical alternatives at the crossing]
\label{cg:crossing-traces}
Suppose $Z$ is nonreduced at $p_0$. If its tangent is $M$, then
$t_{p_0}=1$ and $n_{p_0}\le3$. Otherwise its trace on $M$ is one,
and $2n_{p_0}\le3t_{p_0}+1$. In particular, a crossing zero of
length six has tangent $L$ and trace four on $L$.
\end{lemma}
\begin{proof}
For a curvilinear scheme, any line different from its tangent has
trace one. If its tangent is $M$, use $L=(x=0)$, $M=(y=0)$ and
$I=(y-\varphi(x),x^n)$, with $\operatorname{ord}\varphi\ge2$.
Set $w=y-\varphi(x)$. The coefficient of $w$ in
$Ux^3(w+\varphi(x))$ is $x^3$ times a unit. Since the residual
belongs to $I^2=(w^2,wx^n,x^{2n})$, it follows that $n\le3$.
For every other tangent, the local complete-intersection inequality
gives $2n\le3t_L+t_M=3t_L+1$. At length six and $t_L\le4$,
the only remaining possibility is tangent $L$ with $t_L=4$.
\end{proof}

\begin{theorem}[Exclusion of the cubic--line residual at defect six]
\label{cg:no-nu-six}
In the critical seventeen-line setting, $R=L^3M$ implies $\nu\le5$.
\end{theorem}
\begin{proof}
By Theorem~\ref{thm:L3M6}, $\nu\le6$; suppose $\nu=6$. At every support point away from $p_0$, the strict
vertex inequality~\eqref{cg:strict} holds.  By
Lemmas~\ref{cg:no-reduced-off}--\ref{cg:no-node-off} such a zero is nonreduced and lies at an old
vertex of multiplicity at least three; since $n_p\ge2$, \eqref{cg:strict} gives $t_p\ge2$.

If $p_0\notin Z$, there is either one support point, of length at
most five, or two points of traces two and two, each of length at
most two. Neither gives total length six.

Suppose $p_0\in Z$. There is at most one support point $q\ne p_0$.
If there is none, all six units of length are at $p_0$.
Lemma~\ref{cg:crossing-traces} forces tangent $L$ and trace four.
The off-arrangement crossing lemma
\ref{ae:cubic-crossing}, the smooth-old crossing exclusion
\ref{prop:L3M-no-smooth-crossing}, the old-double exclusion
\ref{dc:node-crossing}, and Lemma~\ref{cg:no-six-at-high-crossing}
rule out every location of $p_0$.

It remains that $q$ exists. If $t_q=2$, its length is two, leaving
length four at $p_0$ with trace at most two. This is impossible:
if its tangent is $M$, its length is at most three; otherwise
$2n_{p_0}\le3t_{p_0}+1\le7$.
Thus $t_q=3$, $t_{p_0}=1$ and $n_q=3$ or $4$ (since $n_q\ge t_q=3$).
Let $m\ge3$ be the old multiplicity at $q$. The normal budget
gives $(m-1)+t_q+1\le7$, where the last summand is the cost of $p_0\in Z$; so $m\le4$.
Lemma~\ref{cg:vertex-order} excludes $n_q=4$, because then its
contact order is three and $9=8+m-1$ would force $m=2$.
For $n_q=3$ it forces $m=4$ and contact order three. Consequently
the normal cost at $q$ is at least six, leaving cost at most one
at $p_0$. Its length is three, its trace on $L$ is one, and its
tangent must be $M$; otherwise the local intersection inequality
would give $2n_{p_0}\le4$.

If $p_0$ is an old vertex, its normal cost is at least two.
If it is a smooth old point, Proposition~\ref{prop:L3M-no-smooth-crossing}
excludes it. If it is off the old arrangement,
Lemma~\ref{cg:off-cross-cost} gives normal cost at least four.
Every case contradicts the remaining budget one.
\end{proof}

\subsubsection{Defect five}

\begin{lemma}[Exclusion of the nilpotent length-five crossing]
\label{cg:no-physical-five}
A nonreduced off-arrangement zero at $p_0$ has tangent $M$, local
length at most three, and normal cost at least four on $L$.
\end{lemma}
\begin{proof}
The initial-form calculation in Lemma~\ref{ae:cubic-crossing}
puts the kernel of the rank-one linearization on $L$ or $M$.
Its proof excludes the semisimple case with kernel $L$; in the
nilpotent case with kernel $L$ it gives exactly
$\operatorname{length}Z_{p_0}=5$ and $t_{p_0}=3$.
We exclude this last possibility by its weighted initial form.

Use $L=(x=0)$ and write $M=(\alpha x+\beta y=0)$, $\beta\ne0$,
with $A=O(2)$ and $B=x+O(2)$. The cubic inflection term forces
$A_2=a_{20}x^2$. Since $\operatorname{ord}B(0,y)=3$, one has
$b_{02}=0$. Comparison of the quartic inflection initial form with
$U_0x^3(\alpha x+\beta y)$ gives $a_{03}=0$ and $c:=a_{12}\ne0$.
Thus $A_x(0,y)$ has exact order two. If
$h=\operatorname{ord}A(0,y)$ and $r=A(0,y)/B(0,y)$, the line
identity is
\[
 A_x(0,y)=r B_x(0,y)-B(0,y)r'/r.
\]
The earlier nilpotent calculation gives $h>3$. If $h=4$, the
first term on the right has order one and the second order two,
contrary to the left side. Hence $h\ge5$.

Assign weights $\mathrm{wt}(x)=3$, $\mathrm{wt}(y)=1$.
The lowest weighted terms are
\[
 A_{[5]}=y^2(cx+ay^3),\qquad B_{[3]}=x+by^3,
 \qquad b\ne0,\quad c\ne0;
\]
here $a=0$ is allowed. Every omitted term has larger weight.
The weight-ten inflection term, after division by $y$ and writing
$z=y^3$, is
\[
 -2cx^3+(-bc-5a)x^2z+
 (ca+cb^2-c^2b-7ab)xz^2+a(a-2b^2-bc)z^3.
\]
It must be a nonzero multiple of $x^3$.
The first unwanted coefficient gives $a=-bc/5$; the next then
becomes $6bc(2b-c)/5$, so $c=2b$ and $a=-2b^2/5$.
The last coefficient becomes $44b^4/25\ne0$, a contradiction.
Thus the tangent is $M$. Lemma~\ref{cg:crossing-traces} gives
length at most three and Lemma~\ref{cg:off-cross-cost} gives
normal cost at least four.
\end{proof}

\begin{lemma}[A triple crossing zero uses four normal zeros]
\label{cg:length-three-cost}
If $n_{p_0}=3$ and $t_{p_0}\le2$, its normal cost on $L$ is
at least four.
\end{lemma}
\begin{proof}
The off-arrangement case is Lemma~\ref{cg:no-physical-five},
and the smooth-old case is excluded by
Proposition~\ref{prop:L3M-no-smooth-crossing}.
At an old double point, Proposition~\ref{dc:node-crossing}
gives tangent $L$, trace two, and a smooth coefficient $a$ with
$b\in\mathfrak m^2$. On the curve $a=0$ the coefficient $b$ has
order three. Its contact with $L$ is two, so Taylor expansion gives
$\operatorname{ord}(b|_L)\ge\min\{3,2+1\}=3$.
The normal cost is therefore at least four.
At an old vertex of multiplicity $m\ge3$, the quartic initial form
forces $a$ smooth, $b\in\mathfrak m^2$, and tangent $L$.
Since the trace is at most two, it is two and the contact order
of $a=0$ with $L$ is two. Restriction to that curve gives
$3\cdot2+1=2\cdot3+m-1$, forcing $m=2$, a contradiction.
\end{proof}

\begin{lemma}[The single old triple zero of length five]
\label{cg:no-pure-cube-five}
An old vertex cannot support a pure-cubic zero of length five and
trace four in the critical setting.
\end{lemma}
\begin{proof}
The ordinary-double case is Lemma~\ref{cg:no-node-off}.
For $m\ge3$, Lemma~\ref{cg:vertex-order} gives
$12=10+m-1$, hence $m=3$. Use coordinates $L=(u=0)$, $v$
along $L$, and write $b=ah+B_0(v)$, with $\operatorname{ord}B_0=5$;
$h$ is a unit, since $a$ and $b$ are smooth with the same tangent (proof of
Lemma~\ref{cg:vertex-order}), so that $b_1=h(0)a_1\ne0$. Both $a|_L,b|_L$ have order four.
The determinant $ba_u-ab_u=B_0a_u-a^2h_u$ has exact order five,
whereas $ba_v-ab_v$ has order at least eight.
The coefficient $\eta(u)|_L$ has exact order two and nonzero
leading coefficient, because $L$ is new.
In the residual restricted to $L$, the transverse Hamiltonian term
$m b\eta(u)(ba_u-ab_u)$ therefore has exact order eleven.
The Euler term $-m(m-1)a^2b$ has order twelve, the remaining
Euler derivative term has order at least thirteen, the tangential
Hamiltonian derivative term at least fourteen, and $\kappa b^3$
at least fourteen. A unique lowest-order term cannot vanish.
\end{proof}

\begin{theorem}[The cubic--line residual has defect at most four]
\label{cg:cubic-four}
In the critical seventeen-line setting, $R=L^3M$ implies $\nu\le4$.
\end{theorem}
\begin{proof}
Theorem~\ref{cg:no-nu-six} reduces the question to $\nu=5$.
If $p_0\notin Z$, two support points have traces two and two and
total length at most four. A single support point of length five
has trace four and is excluded by Lemma~\ref{cg:no-pure-cube-five}.

Suppose all five units of length are at $p_0$.
Lemma~\ref{cg:no-physical-five},
Proposition~\ref{prop:L3M-no-smooth-crossing}, and
Proposition~\ref{dc:node-crossing} exclude the off-arrangement,
smooth-old and old-double locations. At a higher old vertex the
quartic initial form gives $a$ smooth and $b\in\mathfrak m^2$,
with tangent $L$. Since $t=\min\{k,5\}\le4$, the contact order $k$ of $a=0$ with $L$ equals $t$.
The exact residual order is $3t+1=10+m-1$, forcing $t=4,m=4$.
Taylor expansion along $a=0$ gives
$\operatorname{ord}(b|_L)\ge\min\{5,4+1\}=5$,
so its normal cost is at least $3+5=8$, impossible.

It remains that $p_0$ and one point $q\ne p_0$ support $Z$.
If $t_q=2$, then $n_q=2$, leaving $n_{p_0}=3$ with trace at
most two. The normal budget gives $(m_q-1)+t_q+1\le7$, so $m_q\le5$.
Lemma~\ref{cg:vertex-order} gives $3k=m_q+3$ with $k\ge2$,
so necessarily $m_q=3,k=2$. The point $q$ costs at least four
normal zeros; Lemma~\ref{cg:length-three-cost} gives at least
four more at $p_0$, exceeding seven.
If $t_q=3$, its length is three or four. Length four would force
$m_q=2$ by Lemma~\ref{cg:vertex-order} and is impossible.
Thus $n_q=3$. The normal budget
$(m_q-1)+t_q+1\le7$, the last summand being the cost of $p_0\in Z$, gives $m_q\le4$, and
Lemma~\ref{cg:vertex-order} gives $3k=m_q+5$ with $k\ge3$.
Hence $m_q=4$, using six normal zeros and leaving
$n_{p_0}=2,t_{p_0}=1$ with budget at most one. An old vertex
costs at least two, a smooth old point is excluded, and an
off-arrangement point costs at least four by
Lemma~\ref{cg:no-physical-five}. Every case is impossible.
\end{proof}

\subsubsection{Defect four}

\begin{lemma}[The length-four triple-vertex crossing]
\label{cg4:no-triple-four}
An old triple vertex at $p_0=L\cap M$ cannot support a logarithmic
zero of length four and trace three on $L$.
\end{lemma}
\begin{proof}
Use affine linear coordinates $L=(u=0)$ and $v$ along $L$.
Let $q$ be the homogeneous cubic product of the old line equations,
and write $X=aE+b\eta$, with
$E=u\partial_u+v\partial_v$, $\eta=q_v\partial_u-q_u\partial_v$.
The linear logarithmic coefficients are proportional to $u$.
Since the residual has order four, the cubic initial form forces
$a$ smooth and $b\in\mathfrak m^2$; the alternative with $b$
smooth would make its order at least five. Scale the field so that
$a_1=u$. The quartic residual is $-9a_1^2b_2$, so
\[
 b_2=c u(\alpha u+\beta v),\qquad c\beta\ne0.
\]
The trace is three and the local length four, hence the curve
$a=0$ has contact exactly three with $L$. Assign weights
$\mathrm{wt}(u)=3$, $\mathrm{wt}(v)=1$. The lowest terms are
\[
 a_{[3]}=u+Av^3,\qquad
 b_{[4]}=v(Bu+Cv^3),\qquad A\ne0,\quad B\ne0.
\]
Here $B=c\beta$. In fact a pure $v^3$ term in $b$ would give
local length three upon restriction to $a=0$, so it vanishes.
Every omitted coefficient term has strictly larger weight.

Since $L$ is new, $\eta(u)=\lambda v^2$ in its lowest weight,
where $\lambda=3q(0,1)\ne0$. The tangent component of $b\eta$
and the term $\kappa b^3$ have weights too large to contribute
to the weight-ten residual. Substituting the displayed lowest
terms into the exact logarithmic residual, dividing by $3v$ and
putting $z=v^3$, gives
\[
\begin{aligned}
 -3Bu^3&+(-4AB-5C)u^2z\\
 &+(-A^2B-AB^2\lambda-8AC+BC\lambda)uz^2\\
 &+(-3A^2C-ABC\lambda+C^2\lambda)z^3.
\end{aligned}
\]
It must be a nonzero multiple of $u^3$, because the leading
weighted residual is a unit times $u^3v$.
The first unwanted coefficient gives $C=-4AB/5$; the next is then
$9AB(3A-B\lambda)/5$, so $B\lambda=3A$.
The last coefficient becomes $168A^3B/25\ne0$, a contradiction.
\end{proof}

\begin{theorem}[The cubic--line residual has defect at most three]
\label{cg4:cubic-three}
In the critical seventeen-line setting, $R=L^3M$ implies $\nu\le3$.
\end{theorem}
\begin{proof}
Theorem~\ref{cg:cubic-four} reduces the question to $\nu=4$.
Every off-crossing zero is a nonreduced old vertex of multiplicity
at least three and has trace at least two, by
Lemmas~\ref{cg:no-reduced-off}--\ref{cg:vertex-order}.

Suppose first that $p_0\notin Z$. Two support points must both
have trace and length two. At either, the vertex order identity
gives $3k=m+3$ with $k\ge2$, hence $m\ge3$, and its normal cost
is at least four. The two costs exceed seven.
For a single support point of length four, trace three would
give contact three and $9=8+m-1$, forcing $m=2$.
Trace four instead gives contact $k\ge4$ and
$3k=8+m-1$, so $m\ge5$ and normal cost at least eight.
Both cases are impossible.

Suppose all four units of length are at $p_0$.
Lemma~\ref{cg:no-physical-five} excludes the off-arrangement
location, Proposition~\ref{prop:L3M-no-smooth-crossing} the
smooth-old location, and Proposition~\ref{dc:node-crossing}
the old-double location. At a higher old vertex the quartic
initial form forces $a$ smooth, $b\in\mathfrak m^2$, and
tangent $L$. Let $k$ be the contact of $a=0$ with $L$.
The exact restriction identity is
\[
 3k+1=8+m-1,\qquad t=\min\{k,4\}.
\]
If $t=2$, it would give $m=0$. If $t=3$, it gives $k=3,m=3$,
excluded by Lemma~\ref{cg4:no-triple-four}.
If $t=4$, then $k\ge4$, $m\ge6$, and the normal cost is
at least nine. This excludes the single-support case.

Finally suppose $p_0,q$ support $Z$. If $t_q=2$, then $n_q=2$
and $n_{p_0}=2$ with trace at most two. The normal cost at $q$
is at least four. At $p_0$, an off-arrangement zero costs at
least four by Lemma~\ref{cg:no-physical-five}; a smooth-old
zero does not exist; an old double zero cannot have length two
by Proposition~\ref{dc:node-crossing}. At a higher old vertex
of multiplicity $m\ge3$ the quartic initial form forces, as in the
single-support case, $a$ smooth, $b\in\mathfrak m^2$ and tangent $L$,
hence trace two, and its normal cost is at least $(m-1)+2\ge4$.
Every possibility exceeds
the remaining budget three.

If $t_q=3$, then necessarily $n_q=3$ and $n_{p_0}=1$.
The vertex order identity gives $3k=m_q+5$ with $k\ge3$,
so $m_q\ge4$ and the normal cost at $q$ is at least six.
A reduced crossing zero cannot be off the old arrangement
by the first case in the proof of Lemma~\ref{ae:cubic-crossing}, or smooth-old by
Proposition~\ref{prop:L3M-no-smooth-crossing}.
At an old vertex of multiplicity at least three its independent
logarithmic linear coefficients give a nonzero cubic initial
residual, whereas the crossing has order four. Thus only an
old double vertex remains; it costs at least two normal zeros,
again exceeding seven. No case is possible.
\end{proof}

\subsection{The other nonreduced residuals}\label{sec:nonreduced}
By Corollary~\ref{cor:nonreduced} and Theorem~\ref{thm:conic}, the remaining nonreduced residuals are $L^2C$,
$L^2M^2$ and $L^2MN$, with $C$ a smooth conic and distinct line symbols denoting distinct new lines, $L$ not
invariant; for them $\nu\le5$, and $\nu\le4$ if $R=L^2C$ with $L$ transverse to $C$
(Proposition~\ref{prop:nonred91}).

\subsubsection{Equality cases at defect five}
\label{df5:section}

We exclude the equality cases $\nu=5$ for $R=L^2C$ and $R=L^2M^2$ by a local argument.

\begin{lemma}\label{df5:off-inflection-trace}
Let $X=A\partial_x+B\partial_y$ be a holomorphic vector-field germ
with finite curvilinear zero scheme of length $n\ge2$. Suppose the
line $L=(x=0)$ is non-invariant and is a component of the inflection
divisor. If the trace of the zero scheme on $L$ has length $t\ge2$,
then
\[
 n\ge2t-1.
\]
\end{lemma}
\begin{proof}
Curvilinearity and $n\ge2$ make the linearization rank one. Its
kernel is $L$, since the trace has length at least two. A linear
coordinate change preserving $L$, followed by a constant field
rescaling, puts the linearization in semisimple form
$A=x+O(2), B=O(2)$, or nilpotent form $A=O(2),B=x+O(2)$.
Set
\[
 a=A(0,y),\quad b=B(0,y),\quad
 c=A_x(0,y),\quad d=B_x(0,y).
\]
The identity that the inflection determinant vanishes on $L$ is
\begin{equation}\label{df5:off-identity}
 a^2d-acb+abb'-b^2a'=0.
\end{equation}
Non-invariance says $a\not\equiv0$.

In the semisimple case $c$ is a unit and $d(0)=0$. If $b\equiv0$,
\eqref{df5:off-identity} gives $d\equiv0$; on the curve $A=0$,
$B$ then has order at least $2\operatorname{ord}a$, proving the
desired stronger bound $n\ge2t$. Otherwise put $w=b/a$. The exact
identity is
\[
 d=cw-aww'.
\]
A pole in $w$ cannot cancel: if $k=\operatorname{ord}a$ and
$h=\operatorname{ord}b<k$, the two terms on the right have orders
$h-k$ and $2h-k-1$, and $h\ge2$. Thus $w$ is holomorphic.
Also $w(0)=0$, since otherwise $d(0)=c(0)w(0)\ne0$.
The trace length is consequently $t=k$. Writing the smooth curve
$A=0$ as $x=\phi(y)=-a/c+O(a^2)$ gives
\[
 B(\phi(y),y)=b-da/c+O(a^2)
             =a^2ww'/c+O(a^2),
\]
so $n\ge2k=2t$.

In the nilpotent case $d$ is a unit and $c(0)=0$. The function
$b$ cannot vanish identically, by~\eqref{df5:off-identity} and
non-invariance. Put $k=\operatorname{ord}b$, $h=\operatorname{ord}a$
and $r=a/b$. Equation~\eqref{df5:off-identity} becomes
\[
 c=rd-b\frac{r'}r.
\]
If $h<k$, the pole of $rd$ cannot be canceled by the second term,
which is holomorphic of order $k-1\ge1$. If $h=k$, its constant term contradicts $c(0)=0$.
Therefore $h>k$ and $t=k$. On the smooth curve $B=0$, whose equation
is $x=\psi(y)=-b/d+O(b^2)$, one obtains
\[
 A(\psi(y),y)=a-cb/d+O(b^2)
             =\frac{b^2}{d}\frac{r'}r+O(b^2).
\]
Since $r'/r=(h-k)/y+O(1)$, this has exact order $2k-1$.
Thus $n=2t-1$ in this case.
\end{proof}

\begin{proposition}\label{df5:double-conic-four}
If $R=L^2C$ with $C$ a smooth irreducible conic, then $\nu\le4$.
\end{proposition}
\begin{proof}
The transverse case has this bound by Proposition~\ref{prop:nonred91}.  In the tangential
case the only possible zero is the tangency point $p$, and the local bound $n\le\max\{4,t_L+1\}$ of the
proof of Proposition~\ref{prop:nonred91}, with $t_L\le4$, shows that $\nu=5$ would force
\[
 \operatorname{length}Z_p=5,
 \qquad \operatorname{length}(Z_p\cap L)=4.
\]

Suppose first that $p$ is an old vertex. Choose a general smooth
generator of the curvilinear ideal $I_{Z,p}$ whose zero curve has
contact order four with $L$. Its contact order with $C$ is two,
since $C$ has contact exactly two with $L$. Thus the restriction of
$L^2C$ to this curve has order ten. But the vertex identity gives
$R\in\mathfrak m_p I_{Z,p}^2$, so the same restriction has order
at least $2\cdot5+1=11$, a contradiction.

If $p$ lies off the old arrangement, its logarithmic ideal is the
physical zero ideal. Lemma~\ref{df5:off-inflection-trace} gives
$5\ge2\cdot4-1=7$, again impossible.

Finally suppose that exactly one old line passes through $p$.
Take that line as $y=0$ and $L$ as $x=0$, and write
$X=f\partial_x+yg\partial_y$, with logarithmic ideal $(f,g)$.
The conic has local affine equation
\[
 C=x+ux^2+vxy+wy^2,\qquad w\ne0.
\]
The logarithmic residual $P=D/y$ is a unit times $x^2C$.
Since the zero scheme has trace four on $L$, every smooth member
of $(f,g)$ has tangent $L$.

If $f$ is smooth, its zero curve is $x=\phi(y)$ with order exactly
four.  Indeed the trace on $L$ is $\min\{\operatorname{ord}\phi,\operatorname{ord}g(0,y)\}=4$; if
$\operatorname{ord}\phi>4$, then $g(\phi(y),y)\equiv g(0,y)$ modulo $\phi$ would have order four, whereas it
has order equal to the length five.  The same argument, with the roles of $f$ and $g$ exchanged, gives contact
order exactly four in the case $g$ smooth below and in
Proposition~\ref{df5:two-double-lines-four} and Corollary~\ref{oc:three-concurrent-lines}.  Restricting the exact residual
formula to this curve gives
\[
 -y g^2 f_y=U\phi^2 C(\phi,y).
\]
The left side has order $2\cdot5+4=14$, whereas the right side has
order $2\cdot4+2=10$, a contradiction.
Otherwise curvilinearity forces $g$ to be smooth and
$f\in\mathfrak m_p^2$. The exact formula
\[
 P=f^2g_x+fg^2+fyg g_y-fg f_x-yg^2f_y
\]
then has order at least four, while $x^2C$ has order three.
This final contradiction excludes $\nu=5$.
\end{proof}

\begin{proposition}\label{df5:two-double-lines-four}
If $R=L^2M^2$ with $L\ne M$, then $\nu\le4$.
\end{proposition}
\begin{proof}
The zero scheme is supported at $p=L\cap M$. If its length is five,
the local intersection and trace bounds force traces $(4,1)$ on
the two lines, in some order. Denote the line of trace four by $L$.
At an old vertex, restriction to a general smooth generator of the
curvilinear ideal gives residual order $2\cdot4+2\cdot1=10$,
contradicting $R\in\mathfrak m_p I_{Z,p}^2$ and order at least
eleven. Off the old arrangement, Lemma~\ref{df5:off-inflection-trace}
again contradicts length five and trace four.

At a smooth old point put its old line at $y=0$, $L$ at $x=0$,
and write $M=ax+by$; here $a\ne0$ because $M$ is not the old line $y=0$, and $b\ne0$ because
$M\ne L$. The logarithmic residual is
\[
 P=Ux^2(ax+by)^2,\qquad U(0)\ne0.
\]
If $f$ is smooth, its zero curve has equation $x=\phi(y)$ of
order four, and the restriction $P|_{f=0}=-yg^2f_y$ has order
fourteen. The displayed right side has order ten, impossible.

Otherwise $g_x(0)\ne0$ and $f\in\mathfrak m_p^2$.
After a constant field rescaling $g_1=x$. Its zero curve has
contact order four with $L$. Since $f|_{g=0}$ has order five,
the coefficient of $y^2$ in $f_2$ is zero. Write
$f_2=Ax^2+Bxy$. Direct expansion of the exact logarithmic
residual gives
\[
 P_4=A(1-A)x^4-ABx^3y.
\]
Its $x^2y^2$ coefficient is zero, whereas the corresponding
coefficient of $U(0)x^2(ax+by)^2$ is $U(0)b^2\ne0$.
Thus the equality case is impossible.
\end{proof}

\subsubsection{Transverse crossings with a doubled line}
\label{oc:section}

\begin{proposition}\label{oc:reduced-only}
Suppose that the residual germ at a point $p\in Z$ is a unit times
$\ell^2 c$, where $\ell$ is the equation of a new non-invariant
straight line $L$, and $c=0$ is a smooth curve transverse to $L$.
Then $Z_p$ is reduced.
\end{proposition}
\begin{proof}
Suppose $n=\operatorname{length}Z_p\ge2$. Curvilinearity makes
the linearization of the two local logarithmic coefficients rank one.
The residual has nonzero cubic initial form $u\ell^2m$, with
$u\ne0$ and $m$ independent of $\ell$; let $M=\{m=0\}$ be the tangent
line of $c=0$ at $p$.

At an old vertex of multiplicity $s\ge3$, the cubic residual
is $-s(s-1)a_1^2b_1$, with $a_1,b_1$ proportional. It is therefore
a cube or zero. At an old ordinary double point, proportional
linear parts $f_1=A h,g_1=B h$ give cubic residual
$AB(B-A)h^3$, again a cube or zero. Neither can equal
$u\ell^2m$. Thus $p$ is not an old vertex.

Suppose $p$ is a smooth point of one old line. Put that old line at
$y=0$ and $L$ at $x=0$. The physical field and exact residual are
\[
 X=f\partial_x+yg\partial_y,\qquad
 P=f^2g_x+fg^2+fyg g_y-fg f_x-yg^2f_y.
\]
If $f\in\mathfrak m_p^2$ and $g$ is smooth, every term in $P$
has order at least four, a contradiction. Otherwise $f$ is smooth
and $g_1=cf_1$ for a constant $c$. More explicitly, writing
$f=f_1+f_2+O(3)$ and $g=cf_1+g_2+O(3)$, expansion gives
\[
 P_3=f_1\bigl(f_1\partial_xg_2+(cf_2-g_2)\partial_xf_1
             -cf_1\partial_xf_2+c^2f_1^2\bigr).
\]
Thus its cubic part is divisible by $f_1$.
Consequently $f_1$ is proportional to $x$ or to $m$.

If $f_1\sim x$, division of the field by a unit normalizes
$f=x-\phi(y)$. Write
\[
 \phi=qy^2+O(y^3),\qquad
 g=cx+Ax^2+Bxy+Cy^2+O(3).
\]
Direct expansion gives
\[
 P_3=(A+c^2)x^3-(C+cq)xy^2.
\]
It has zero $x^2y$ coefficient, whereas that coefficient of
$u x^2m$ is nonzero because $M\ne L$.

If $f_1\sim m$ and $M$ is not the old line, constant rescaling puts
$M=(x-y=0)$.
Normalize $f=x-y-qy^2+O(y^3)$ and write
$g=c(x-y)+Ax^2+Bxy+Cy^2+O(3)$. The exact cubic is
\[
 P_3=(x-y)\bigl((A+c^2)x^2-2(A+c^2)xy
               +(-B-C+c^2-cq)y^2\bigr).
\]
It cannot be $u x^2(x-y)$: the $x^2$ coefficient in the parenthesis
would be $u$, whereas its $xy$ coefficient would then be $-2u$
instead of zero. This excludes smooth old points. The argument uses
only the tangent of $c=0$, so curvature of that curve is harmless.
If instead $M$ is the old line $y=0$, the remaining possibility
is $f_1\sim y$. Normalize
$f=y-qx^2+O(3)$, $g=cy+Ax^2+Bxy+Cy^2+O(3)$.
The cubic residual is
\[
 P_3=2(A+cq)xy^2+(B+c^2)y^3,
\]
which has zero $x^2y$ coefficient and cannot equal $u x^2y$.
This also excludes tangency of $c=0$ to the old line.

It remains to treat a point off the old arrangement. The logarithmic
ideal is the physical zero ideal. Put the rank-one linearization in
semisimple or nilpotent form.

In the semisimple case write
$X=A\partial_x+B\partial_y$ with $A=x+O(2)$, $B=O(2)$.
The kernel line is $x=0$, and, writing
$B_2=b_{20}x^2+b_{11}xy+b_{02}y^2$, its cubic inflection form is
\[
 D_3=x(b_{20}x^2-b_{02}y^2).
\]
Its factorization as a nonzero doubled line times a transverse
simple line forces the simple line to be the kernel $x=0$ and
the doubled line to be $y=0$: if the kernel were the doubled
line, the missing $xy$ coefficient in the displayed quadratic
would be a contradiction; if the kernel is the simple line,
that quadratic is a square only when its repeated factor is $y$.
Thus $L$ is the nonzero eigendirection of the linearization.

But an inflection line in this eigendirection must be invariant.
Indeed, if $B(x,0)$ is not identically zero, let its first term be
$b_kx^k$, $k\ge2$. On $y=0$ the exact determinant
\[
 D=A^2B_x+AB(B_y-A_x)-B^2A_y
\]
has first term $(k-1)b_kx^{k+1}\ne0$. Here
$A(x,0)=x+O(x^2)$ and $B_y(x,0),A_y(x,0)=O(x)$, so all omitted
contributions have strictly higher order. This contradicts $D|_L=0$.
Hence $L$ would be invariant, contrary to the hypothesis.

In the nilpotent case write $A=O(2)$, $B=x+O(2)$ and
$A_2=a_{20}x^2+a_{11}xy+a_{02}y^2$. Then
\[
 D_3=-x^2(a_{11}x+2a_{02}y).
\]
The doubled line is the kernel $L=(x=0)$ and $a_{02}\ne0$.
The physical ideal consequently has length two and trace two on
$L$: its smooth generator $B$ has zero curve $x=O(y^2)$, and
$A$ restricts there to $a_{02}y^2+O(y^3)$.
Lemma~\ref{df5:off-inflection-trace} would give $2\ge3$, a final
contradiction.
\end{proof}

\begin{corollary}\label{oc:transverse-double-conic}
If $R=L^2C$ and the smooth irreducible conic $C$ meets $L$
transversely, then $\nu\le2$.
\end{corollary}
\begin{proof}
The zero scheme is supported at the two transverse intersections,
and Proposition~\ref{oc:reduced-only} makes each contribution at
most one, including when the tangent of $C$ is an old line.
\end{proof}

\begin{corollary}\label{oc:three-nonconcurrent-lines}
If $R=L^2MN$, where $L,M,N$ are distinct and not concurrent,
then $\nu\le3$.
\end{corollary}
\begin{proof}
The support is contained in $L\cap M$, $L\cap N$ and $M\cap N$.
At the first two points Proposition~\ref{oc:reduced-only} gives
local length at most one. At the last point the residual is an
ordinary node and also has local length at most one.
\end{proof}

\begin{corollary}\label{oc:three-concurrent-lines}
If $R=L^2MN$, where the three distinct lines are concurrent,
then $\nu\le4$.
\end{corollary}
\begin{proof}
By the proof of Proposition~\ref{prop:nonred91}, $\nu\le5$, and the equality case is a
single curvilinear scheme of length five at the common point,
with trace four on $L$ and trace one on $M,N$.
At an old vertex, restriction to a general smooth generator of
its ideal gives residual order ten, contradicting the strict
vertex condition $R\in\mathfrak m I_Z^2$, which requires at
least eleven. Off the old arrangement,
Lemma~\ref{df5:off-inflection-trace} contradicts length five
and trace four.

At a smooth old point use its old line $y=0$, $L=(x=0)$,
$M=ax+by$, $N=cx+dy$, with $abcd\ne0$. If $f$ is smooth,
restriction to $f=0$ has order fourteen, whereas
$x^2MN$ has order ten, exactly as in
Proposition~\ref{df5:two-double-lines-four}. Otherwise $g_1=x$
after a constant rescaling, and $f_2=Ax^2+Bxy$. The quartic
residual $A(1-A)x^4-ABx^3y$ has zero $x^2y^2$ coefficient,
whereas the corresponding coefficient of $U(0)x^2MN$ is
$U(0)bd\ne0$. This excludes equality.
\end{proof}

\subsection{A reduced residual at defect five}\label{sec:reduced}
By Theorems~\ref{thm:reduced6} and~\ref{thm:sixnodes}, a reduced residual has $\nu\le5$.  In this section
we exclude $\nu=5$.

\subsubsection{The three configurations}
\label{red5:section}

The following local observations use the foliation equations, in addition
to the ideal-square condition. This distinction is necessary: an
arbitrary curvilinear length-two scheme can satisfy the ideal-square
condition at a tacnode, but not the residual identities below.

\begin{lemma}\label{red5:local}
At a point $p\in Z$:
\begin{enumerate}
\item if the residual has multiplicity two, then
$\operatorname{length}Z_p=1$ and $p$ is not an old vertex;
\item if the residual has multiplicity three and
$\operatorname{length}Z_p\ge3$, its cubic initial form is a cube.
\end{enumerate}
\end{lemma}
\begin{proof}
At residual multiplicity at most three, $R\in I_Z^2$ excludes
$I_Z\subset\mathfrak m^2$. Thus the two-generated finite local ideal
contains a smooth parameter and is curvilinear; a nonreduced zero has
rank-one logarithmic linear parts.
Off the old arrangement, a nonreduced curvilinear zero makes the
physical linearization rank one, so the quadratic inflection form
$\det(J)\det(w,Jw)$ vanishes. At a smooth old point, write
$X=f\partial_x+yg\partial_y$. Its quadratic residual is
$-\det(df_1,dg_1)yf_1$, which also vanishes at a nonreduced zero.
At an old vertex $R\in\mathfrak m I_Z^2\subset\mathfrak m^3$.
These facts prove the first assertion.

For the second assertion, first work off the arrangement. In rank-one
semisimple coordinates $A=x+O(2)$, $B=B_2+O(3)$,
\[
 D_3=x(b_{20}x^2-b_{02}y^2).
\]
Length at least three means that $B$ restricted to the smooth curve
$A=0$ has no quadratic term; hence $b_{02}=0$.
In nilpotent coordinates $A=A_2+O(3)$, $B=x+O(2)$,
\[
 D_3=-x^2(a_{11}x+2a_{02}y),
\]
and length at least three similarly forces $a_{02}=0$.
Thus the nonzero cubic initial form is a cube in either case.

At a smooth old point, if $f\in\mathfrak m^2$ and $g$ is smooth,
every residual term has order at least four. Otherwise $f$ is smooth.
Write $g_1=cf_1$ and put $k_2=g_2-cf_2$. Length at least three makes
$k_2$ vanish on the line $f_1=0$, so $k_2=f_1h_1$ for a linear form
$h_1$. The exact cubic residual is
\[
 P_3=f_1\bigl(f_1\partial_xk_2-k_2\partial_xf_1+c^2f_1^2\bigr)
     =(\partial_xh_1+c^2)f_1^3.
\]
Finally, at an old ordinary double point proportional logarithmic
linear parts $f_1=Ah$, $g_1=Bh$ give residual cubic
$AB(B-A)h^3$. At an old vertex of multiplicity $m\ge3$, the
radial--Hamiltonian frame $X=aE+b\eta$ gives residual cubic
$-m(m-1)a_1^2b_1$; the proportionality of $a_1,b_1$ again makes
this a cube or zero. This proves the second assertion.
\end{proof}

\begin{theorem}\label{red5:classification}
Suppose the critical residual quartic is reduced and $\nu=5$.
Then exactly one of the following necessary configurations occurs:
\begin{enumerate}
\item a smooth conic and two lines, with five distinct ordinary nodes,
and $Z$ the reduced scheme of all five;
\item four lines in general position, with $Z$ five of their six nodes;
\item four lines with exactly three concurrent, with local zero length
two at the triple point and length one at each of the other three nodes.
\end{enumerate}
\end{theorem}
\begin{proof}
The reduced quartic degree partitions are
$(4),(3,1),(2,2),(2,1,1),(1,1,1,1)$.
The support of $Z$ lies in its singular set.

An irreducible quartic has total delta-invariant at most three.
Its double singularities each contribute at most one by
Lemma~\ref{red5:local}. A triple point uses the entire genus budget,
so there is no other singularity. If its local zero length were at
least three, the same lemma would make the cubic initial form $y^3$.
In affine linear coordinates its equation is $y^3+R_4(x,y)$, and
irreducibility forces a nonzero $x^4$ term. The curvilinear supporting
curve has tangent $y=0$, so restriction of this equation to it has
order exactly four. This contradicts membership in $I_Z^2$ at length
at least three. A multiplicity-four point would make the entire
quartic a product of lines. Thus the irreducible case contributes at
most three.

For an irreducible cubic and a line, if the line avoids the possible
unique singular point of the cubic, its at most three intersection
points contribute at most one each, and that cubic singularity at
most one. If the line passes through the cubic singularity, the
common point has multiplicity three and there is at most one further
intersection. A noncube cubic initial form contributes at most two.
The cube case is a cusp with its tangent line: putting that line at
$y=0$, the product is $y^3+yC_3(x,y)$ with a nonzero $x^3$ term in
$C_3$. On a curvilinear support $y=\phi(x)$, $\operatorname{ord}\phi\ge2$,
the linear coefficient in the transverse generator has order exactly
three. Membership in $I_Z^2$ therefore bounds the local length by
three. The tangent already exhausts its three intersections with the
cubic. In every case the total is at most four.

Two irreducible conics are smooth. Their at most four intersection
points have residual multiplicity two, including at tangencies, and
so contribute at most four in total.

For a conic and two lines meeting outside the conic, there are at
most five singular points, each of multiplicity two. Length five
therefore forces five distinct ordinary nodes: both conic--line
intersections on each line must be distinct and transverse.
If the two lines meet on the conic, the common point has a noncube
cubic initial form and contributes at most two; each line has at
most one further conic intersection, so the total is at most four.

For four lines in general position there are six nodes, each of
length at most one, giving the second alternative. If exactly three
are concurrent, their noncube triple tangent cone contributes at
most two, and the other three nodes at most one each; length five
forces the third alternative. If all four are concurrent, the
ideal-square bound for four distinct concurrent lines gives total
length at most four. These exhaust the intersection patterns.
\end{proof}

\subsubsection{Four lines, three of them concurrent}
\label{tc4:section}

\begin{lemma}[A semisimple kernel inflection line]
\label{tc4:semisimple-trace}
Let a physical holomorphic vector field have an isolated curvilinear
zero of length $N\ge2$ and a rank-one semisimple linearization.
If its kernel line $K$ is a non-invariant straight inflection component,
and its physical zero scheme has trace $t\ge2$ on $K$, then $N\ge2t$.
\end{lemma}
\begin{proof}
Put $K=(x=0)$ and normalize the linear part to
$A=x+O(2)$, $B=O(2)$. Write
$a=A|_K$, $b=B|_K$, $c=A_x|_K$, and $d=B_x|_K$.
Non-invariance means $a\not\equiv0$, and restriction of the
inflection determinant gives
\[
 a^2d-acb+abb'-b^2a'=0.
\]
If $b\equiv0$, this gives $d\equiv0$ and restriction to the smooth
curve $A=0$ gives $N\ge2\operatorname{ord}a=2t$.
Otherwise put $w=b/a$. The identity reads $d=cw-aww'$.
A pole in $w$ cannot cancel: writing $k=\operatorname{ord}a$,
$h=\operatorname{ord}b\ge2$, the orders of the two terms are
$h-k$ and at least $2h-k-1>h-k$. Hence $w$ is holomorphic.
Since $c$ is a unit and $d(0)=0$, one has $w(0)=0$, so $t=k$.
The curve $A=0$ is $x=-a/c+O(a^2)$; on it,
\[
 B=b-da/c+O(a^2)=a^2ww'/c+O(a^2).
\]
This gives $N\ge2k$.
\end{proof}

\begin{lemma}[The invariant line through the length-two triple zero]
\label{tc4:kernel-invariant}
Suppose that the residual germ at $p$ is a unit times the product of
three distinct linear forms and that $\operatorname{length}Z_p=2$. Then exactly one of these three
branches is an invariant line $K$. At $p$ the physical linearization
is rank-one semisimple, with zero tangential eigenvalue along $K$ and
nonzero transverse eigenvalue.
\end{lemma}
\begin{proof}
At an old vertex of multiplicity at least three, the rank-one
logarithmic linear parts make the cubic residual a cube or zero.
At an old ordinary double point the same conclusion follows from
the proportional logarithmic linear parts. Neither is an ordinary
triple tangent cone, so $p$ is not an old vertex.

Off the old arrangement, put the rank-one physical linearization
in semisimple or nilpotent form. In the latter form,
\[
 D_3=-x^2(a_{11}x+2a_{02}y),
\]
which has a repeated factor. In semisimple form,
\[
 A=x+O(2),\quad B=O(2),\qquad
 D_3=x(b_{20}x^2-b_{02}y^2).
\]
Distinctness forces $b_{20}b_{02}\ne0$. The kernel line $K=(x=0)$
is a residual branch, and the physical zero scheme has length and
trace both two there. Lemma~\ref{tc4:semisimple-trace} makes
non-invariance impossible.

At a smooth point of one old line put that line at $y=0$ and write
$X=f\partial_x+yg\partial_y$. If $f\in\mathfrak m^2$, the residual
has order at least four, contrary to the ordinary cubic initial form.
Thus $f$ is smooth and $g_1=cf_1$. The cubic residual is divisible by
$f_1$, so its tangent is one of the three residual lines. It cannot
be the old tangent: with $f_1=y$, the cubic residual is divisible by
$y^2$. Put the new tangent $K$ at $x=0$. The physical linearization
then has nonzero normal eigenvalue and zero tangent eigenvalue.
The physical zero ideal is $(f,yg)$, of colength three: on the smooth
curve $f=0$, the factors $y$ and $g$ have orders one and two.
If $K$ were non-invariant, its physical trace would be at least two,
contradicting Lemma~\ref{tc4:semisimple-trace}.

In either local position the other two residual tangents are not
eigendirections. Indeed the displayed semisimple cubic has neither
of its other factors equal to the nonzero eigendirection. In the
smooth-old coordinates its analogue is
$x(\alpha x^2-\beta y^2)$ with $\alpha\beta\ne0$.
An invariant smooth branch must have an eigendirection as tangent.
Consequently neither of the other two residual lines is invariant.
\end{proof}

\begin{lemma}[Critically fixed maps have no parabolic fixed point]
\label{tc4:critical-fixed}
Let $F:\mathbb P^1\to\mathbb P^1$ have degree $d\ge2$, and suppose
that every critical point is fixed. Every fixed point outside the
critical set has multiplier of modulus greater than one.
\end{lemma}
\begin{proof}
Write $P=\operatorname{Crit}(F)$. It has at least two points by
Riemann--Hurwitz, since one point has ramification at most $d-1$.
If it has exactly two points, both are totally ramified. Conjugating
them to zero and infinity gives $F(z)=c z^d$, and every other fixed
point has multiplier $d$.

Suppose $|P|\ge3$. Put
$X=\mathbb P^1\setminus P$ and
$X'=\mathbb P^1\setminus F^{-1}(P)$. The map $F:X'\to X$ is an
unbranched covering and $X'\subsetneq X$. Indeed all critical values
belong to $P$, and
\[
 |F^{-1}(P)|=d|P|-(2d-2)>|P|.
\]
Both surfaces are hyperbolic. Their Poincar\'e metrics satisfy
$\rho_{X'}(z)=\rho_X(F(z))|F'(z)|$, while strict Schwarz--Pick \cite{Mil06}
monotonicity for the proper inclusion gives $\rho_{X'}>\rho_X$.
For clarity, strictness follows by lifting the inclusion and universal
disk coverings to a disk self-map: equality at one point would make
that lift an automorphism and make the inclusion surjective.
At a fixed point outside $P$, the metric identity therefore gives
$|F'|=\rho_{X'}/\rho_X>1$.
\end{proof}

\begin{theorem}[Exclusion of the triple-concurrent defect-five branch]
\label{tc4:excluded}
The residual cannot consist of four distinct lines with exactly three
concurrent, with local zero length two at the triple point and reduced
zeros at the other three residual nodes.
\end{theorem}
\begin{proof}
Write the concurrent lines as $K,L,M$, their common point as $p$,
and the fourth line as $N$. Every residual component is new by
Proposition~\ref{prop:riccati}; each concurrent line contains
logarithmic trace at least two, and $N$ contains trace three.
Lemma~\ref{tc4:kernel-invariant} makes
$K$ invariant, with trace of $Z$ on $K$ equal to three: two at $p$
and one at $q=K\cap N$.

Let $h$ be the number of distinct old intersection points on $K$.
Its nonzero tangential field is a section of $\mathcal O_K(8)$
(Section~\ref{sec:foliation}).
At each old intersection it vanishes to order at least one plus
the local logarithmic trace; away from the old arrangement its zero
order is at least that trace. Hence $h+3\le8$, so $h\le5$.

There is also $h\ge5$. Add the invariant line $K$ and choose a
homogeneous representative $\delta$ of degree seven with
$\delta(K)=0$, by subtracting an Euler multiple. Its restriction
to $K$ is a nonzero logarithmic derivation of the Ziegler
multirestriction of the eighteen-line arrangement. Nonzeroness follows
from saturation: an identically zero restriction would make the
projective field vanish along $K$. The multiplicities $m_1,\ldots,m_h$
sum to seventeen and are at most eight: the critical arrangement is
not a pencil, and an old pencil of at least nine lines would give an
outside old line at least nine vertices. The
balanced rank-two multiplicity gap bound \cite{Abe} gives
$|d_2-d_1|\le h-2$. Since $d_1\le7$ and $d_1+d_2=17$, it gives
$h\ge5$ (this is Lemma~\ref{lem:invline}). Thus $h=5$, and every tangential vanishing bound on $K$ is an
equality.

Choose homogeneous coordinates $[s:t:z]$ with $K=(z=0)$, and write
$\delta|_K=(A(s,t),B(s,t),0)$. Define
$F=[A:B]$. We check that the degree-seven forms $A,B$ are coprime.
A common zero would be a singularity of the projective field on $K$
with zero normal eigenvalue. All such singularities are among the
already exhausted tangential zeros: the five old intersections,
$p$ if it is not old, and $q$.
At $p$, Lemma~\ref{tc4:kernel-invariant} gives a nonzero normal
eigenvalue. The reduced residual node $q$ is off the old arrangement
and has invertible physical linearization. Indeed an old vertex is
excluded by the strict vertex residual order, and at a smooth old
point the quadratic logarithmic residual contains the old tangent,
whereas neither residual line is old.
At any other old intersection, tangential equality gives a nonzero
tangent eigenvalue. If at least two old lines pass there, the three
invariant tangents including $K$ force a nonzero scalar linearization.
If only one old line passes there, the residual is locally a simple
copy of $K$: the nonzero quadratic inflection term, equal to a unit
times the old line and $K$, forces both eigenvalues nonzero.
Thus there is no common zero, and $F$ is a morphism of degree seven.

Every old intersection is a fixed point of $F$, with local degree at
least $m_i$. To see the multiplicity directly, in linear coordinates
centered there the $m_i$ old lines have equations $t-a_jz=0$.
Since $\delta(z)=0$, their invariance makes $\delta(t)$ divisible
by their product. Restricting to $z=0$ gives divisibility by $t^{m_i}$.
Coprimality then gives the asserted local degree. But
\[
 \sum_{i=1}^5(m_i-1)=17-5=12=2\cdot7-2.
\]
Riemann--Hurwitz implies that every critical point of $F$ is fixed.

Finally $p$ is a parabolic fixed point of $F$. In an affine chart of
$K$ centered at $p$, the normalized representative gives the physical
field
\[
 \dot z=-zA(1,t)+O(z^2),\qquad
 \dot t=B(1,t)-tA(1,t)+O(z).
\]
If $\mu\ne0$ and $\lambda=0$ are respectively its normal and tangent
eigenvalues at $p$, then
\[
 F'(p)=1-\frac{\lambda}{\mu}=1.
\]
This contradicts Lemma~\ref{tc4:critical-fixed} and excludes the branch.
\end{proof}

\subsubsection{A conic and two lines}
\label{clm:section}
Assume that $R=CLM$, where $C$ is a smooth conic, $L,M$ are
 distinct lines transverse to $C$, and their five pairwise intersections
 are distinct and form the reduced scheme $Z$. We exclude this
 configuration by forcing both lines to be invariant and applying an
 exact finite incidence certificate.

Each of $L,M$ contains three reduced points of $Z$, so both are new by
 Proposition~\ref{prop:riccati}.
 Their common node lies off the old arrangement. Indeed, an old vertex
 is excluded by $R\in\mathfrak m I_Z^2$, which would force residual
 multiplicity at least three. At a smooth old point, the quadratic
 residual form is a nonzero scalar times $y f_1$ in old-line
 coordinates $y=0$; a crossing of the two straight components $L,M$
 would make one of them the old line. Thus their common physical
 singularity is nondegenerate, with $L,M$ as its eigendirections.

\begin{lemma}\label{clm:noninvariant-node}
Let $J$ be one of $L,M$. At a point of $C\cap J$ outside the old
 arrangement, if $J$ is non-invariant, its normal vanishing order is
 at least three.
\end{lemma}
\begin{proof}
The physical linearization is invertible and has the distinct tangent
 directions of $C,J$ as eigendirections. Projective coordinates may be
 chosen so that $C:y=x^2$ and $J:x=0$. Write
\[
 A=\lambda x+O(2),\qquad B=\mu y+d x^2+O(xy,y^2,\mathfrak m^3),
 \qquad \lambda\mu\ne0.
\]
For $D=A^2B_x+AB(B_y-A_x)-B^2A_y$, the coefficient of $t^3$ in
 $D(t,t^2)$ is
\[
 \lambda\bigl((\lambda+\mu)d+\mu(\mu-\lambda)\bigr).
\]
It vanishes since $C\mid D$. In particular $\lambda+\mu\ne0$.
If $n\ge2$ is the normal order of the straight inflection component
 $J$, restriction of $D$ to $J$ gives the resonance
 $\lambda=-(n-1)\mu$; this is \eqref{eq:resonance}, since here $\mu$ is the eigenvalue along $J$ (the
 direction $\partial_y$) and $\lambda$ the transverse one. The case $n=2$ would give
 $\lambda+\mu=0$, a contradiction.
\end{proof}

\begin{lemma}\label{clm:old-smooth-node}
A point of $C\cap J$ on the old arrangement lies on precisely one old
 line $H$, which is tangent to $C$ there. Such a point forces $J$ to
 be invariant.
\end{lemma}
\begin{proof}
An old vertex is impossible because the residual has multiplicity two
 while $R\in\mathfrak m I_Z^2$ at an old vertex. At a smooth old point,
 use $H:y=0$ and $J:x=0$ and write
 $X=f\partial_x+yg\partial_y$. The reduced logarithmic zero makes
 $f_1,g_1$ independent, and
\[
 \operatorname{in}_2(D/y)=-\det(f_1,g_1)y f_1.
\]
Consequently $H$ is tangent to $C$ and
 $f_1=a x$, $g_1=cx+dy$, with $ad\ne0$.
Suppose $J$ is non-invariant and put
 $h=\operatorname{ord}f(0,y)\ge2$. In the physical field
 $A=f$, $B=yg$ one has $B(0,y)=d y^2+O(y^3)$,
 $A_x(0)=a\ne0$, $B_x(0,y)=O(y)$, and $B_y(0,y)=O(y)$.
 On $J$, the term $-ABA_x$ in
\[
 D=A^2B_x+AB(B_y-A_x)-B^2A_y
\]
 has nonzero leading term of order $h+2$. The other three terms have
 orders at least $2h+1,h+3,h+3$, all strictly larger. This contradicts
 $D|_J=0$. Hence $f|_J=0$ and $J$ is invariant.
\end{proof}

\begin{proposition}\label{clm:both-invariant}
Both $L$ and $M$ are invariant.
\end{proposition}
\begin{proof}
Suppose $L$ is non-invariant. Lemma~\ref{clm:old-smooth-node}
 puts both its conic nodes outside the old arrangement, and
 Lemma~\ref{clm:noninvariant-node} assigns normal order at least three
 to each. The line-line node is outside the old arrangement, is
 nondegenerate, and contributes normal order at least two. This
 exceeds the degree-seven normal divisor on $L$.

\end{proof}

\begin{proposition}\label{clm:catalog}
Each of $L,M$ meets the old arrangement in exactly five distinct
 points, all its old intersections and three residual nodes accounting
 for the degree-eight tangential zero divisor. The associated map
 $f_J:J\to J$ has degree seven and all its critical points are fixed.
 Every residual node is off the old arrangement.
\end{proposition}
\begin{proof}
An invariant new line has at least five old intersections
 (Lemma~\ref{lem:invline}). At an old conic node the
 tangential field has order two; at each other old point it has order
 at least one; every off-old residual node contributes at least one.
 Therefore its tangential degree-eight divisor gives $r+3\le8$.
 Hence $r=5$ and all these lower bounds are equalities.

Normalize a homogeneous degree-seven representative by
 $\theta(J)=0$. Its restricted pair defines $f_J$. At every tangential
 zero away from an old conic node the physical normal eigenvalue is
 nonzero: this follows from invertibility at off-old residual nodes,
 from scalar nonzero linearization at an old point incident with at
 least two old lines, and from the product of two distinct initial
 inflection directions at a remaining simple old intersection.
 At an old conic node, Lemma~\ref{clm:old-smooth-node} gives one
 nonzero normal eigenvalue and zero tangential eigenvalue. Thus the
 restricted pair has no common zero, and the map has degree seven.
 Its local degree at an old point of multiplicity $m_i$ is at least
 $m_i$. The identity $\sum(m_i-1)=17-5=12=2\cdot7-2$ exhausts
 Riemann--Hurwitz, proving that these are all critical points and that
 every critical point is fixed.

An old conic node would have multiplier one for $f_J$, because its
 tangential eigenvalue vanishes and its normal eigenvalue is nonzero.
 A map of degree at least two with every critical point fixed has no
 parabolic fixed point, by Lemma~\ref{tc4:critical-fixed}.
 Hence every residual node is off the old arrangement.
\end{proof}

\begin{lemma}\label{clm:old-multiplicity}
Every old vertex has multiplicity at most five.
\end{lemma}
\begin{proof}
Let $p$ be an old vertex of
 multiplicity $m\ge6$, represented by a vector $p\in\CC^3$, and let
 $W_p(w)=\det(p,w,\theta(w))$, a form of degree eight in $w\in\CC^3$. The product
 $Q_p$ of its old lines divides $W_p$. The form $W_p/Q_p$, of degree
 $8-m\le2$, vanishes at all five reduced points of $Z$; the preceding proposition
 ensures $Q_p$ is nonzero there. Since each of $L,M$ contains
 three of these points, a nonzero such form must be a multiple of
 $LM$. Thus $W_p$ vanishes identically on both lines, also if the
 quotient is zero. On an invariant line $J$, the foliation is nonzero
 at its generic point by saturation; hence $W_p|_J=0$ implies $p\in J$.
 This forces $p=L\cap M$, contrary to that node lying off the old
 arrangement.
\end{proof}

\begin{lemma}[Two-line incidence certificate]\label{clm:incidence}
Let $f_I$ be one of the maps $A$, $B$ of Proposition~\ref{cat:multipliers} and $f_J$ one of the maps of
type $C$.  Put $f_I$ on a line $I$ and $f_J$ on a line $J$, choose a noncritical fixed point $p$ of each
map and identify the two points $p$ with $I\cap J$; choose a second noncritical fixed point of $f_I$, and
mark the four critical points of $f_I$ and this point with the multiplicities $5,3,3,5,1$ (type $A$) or
$4,4,4,4,1$ (type $B$), and the five critical points of $f_J$ with their local degrees $4,3,3,3,4$.
Let $\mathcal A$ be a set of seventeen lines, each joining a marked point of $I$ to a marked point of $J$,
such that every marked point lies on exactly as many lines of $\mathcal A$ as its mark.  Then some line of
$\mathcal A$ meets the other sixteen in at least $13$ distinct points.
\end{lemma}
\begin{proof}[Computer-assisted proof]
The two lines can be placed at $\{y=0\}$ and $\{y=z\}$, with $p=(1:0:0)$; in the coordinate
$u=1/(x-p)$ on each line the marked points are $(u_i:0:1)$ and $(v_j:1:1)$, and independent affine
changes of the two coordinates are induced by projective transformations preserving the two lines:
$(x:y:z)\mapsto(Ax+By+Cz:\lambda y:\lambda y+e(z-y))$, with $A\lambda e\ne0$, fixes $p$ and induces
$u\mapsto(Au+C)/e$ on $\{y=0\}$ and $v\mapsto(Av+B+C)/\lambda$ on $\{y=z\}$.  A set
$\mathcal A$ is a $5\times5$ matrix of zeros and ones with the marks as row and column sums: there are
$18$ such matrices for type $A$ and $84$ for type $B$.  The choices of the map on $I$ ($A$, or one of the
two Galois conjugate maps $B$; these are Möbius equivalent, and both are included), of $p$ on $I$ and on
$J$ and of the extra point give $108$ marked pairs and $6\,696$ sets $\mathcal A$.  (The map $C$ is defined
over $\QQ$; relabelling its two conjugate critical points of local degree $3$ gives the same marked set.)

All coordinates lie in the splitting field $E$ of the fixed-point polynomials of the maps $A$, $C$ and
the two conjugate maps $B$.  We chose
a prime $\ell$ which splits completely in $E$ and such that the fixed points of each map stay distinct
modulo $\ell$ (namely $\ell=10\,009$), and computed in $\mathbb F_\ell$ over all choices of the relevant
roots modulo $\ell$, a set which contains the reductions coming from every embedding of $E$.  For every
marked pair and every matrix some line has at least $13$ distinct intersection points modulo $\ell$: the
largest number of distinct intersection points on a line is $13$, $14$, $15$ in $888$, $3\,912$, $1\,896$
cases respectively.  Since the coordinates are integral at a prime of $E$ above
$\ell$ and points which are distinct modulo that prime are distinct, the same holds in characteristic
zero.
\end{proof}

\begin{theorem}\label{clm:no-defect-five}
In the critical seventeen-line setting with defect five, the reduced
 residual quartic cannot be a smooth conic and two lines.
\end{theorem}
\begin{proof}
By the preceding propositions both lines are invariant, all five
 residual nodes are off the old arrangement, and old multiplicities
 are at most five, and by Proposition~\ref{clm:catalog} the two lines carry maps $f_L$, $f_M$ as in
 Lemma~\ref{tf4:no-two-catalog-lines} below.  At their common node $p=L\cap M$, a nondegenerate
 singularity with eigenvalue $\lambda$ along $L$ and $\mu$ along $M$, one has $f_L'(p)=1-\lambda/\mu$ and
 $f_M'(p)=1-\mu/\lambda$, as in the proof of Theorem~\ref{tc4:excluded}; so
 $(1-f_L'(p))(1-f_M'(p))=1$, and Lemma~\ref{tf4:no-two-catalog-lines} excludes the configuration.
\end{proof}

\begin{lemma}[Two invariant catalog lines]\label{tf4:no-two-catalog-lines}
Let $\mathcal A$ be an arrangement of seventeen lines with
 $h(H)\le8$ for every $H\in\mathcal A$. Suppose that two new lines
 $I,J$, meeting at $p\notin\mathcal A$, have the following properties.
 Each meets the old arrangement in exactly five distinct points and
 carries a degree-seven rational map whose eight fixed points are
 these five points and three further distinct points, including $p$.
 Its critical divisor is
\[
 \sum_{q\in I\cap\mathcal A}(m_q-1)[q]
 \quad\hbox{or}\quad
 \sum_{q\in J\cap\mathcal A}(m_q-1)[q],
\]
 respectively, where $m_q\le5$ is the old multiplicity. Finally,
 their multipliers at $p$ obey
\[
 (1-f_I'(p))(1-f_J'(p))=1.
\]
Then this configuration is impossible.
\end{lemma}
\begin{proof}
Both maps belong to the catalog of Proposition~\ref{cat:multipliers}, since their critical points are
 fixed and have local degree at most five.  By part~(b) of that proposition one of them, say $f_I$, is of
 type $A$ or $B$ and the other of type $C$, and $p$ is a noncritical fixed point of both.  The old points
 of $I$ are its four critical points and one noncritical fixed point other than $p$, with the old
 multiplicities as marks; the old points of $J$ are its five critical points.  Since $p$ is off the
 arrangement, every old line is the join of an old point of $I$ and an old point of $J$.  So the old
 arrangement is one of the sets of Lemma~\ref{clm:incidence}, and some old line has at least $13$
 vertices, contradicting $h(H)\le8$.
\end{proof}

\subsubsection{Four lines in general position: the full lines}\label{fi5:section}

Suppose throughout this section that
\[
 R=L_1L_2L_3L_4,\qquad
 Z=\{q_{13},q_{14},q_{23},q_{24},q_{34}\},\qquad q_{ij}=L_i\cap L_j,
\]
where the four residual lines are in general position.  The lines
$L_3,L_4$ will be called full, and $L_1,L_2$ deficient.

\begin{lemma}\label{fi5:marked}
The five marked nodes lie off the old arrangement.  Their physical
singularities are nondegenerate, with the two residual directions as
their eigendirections.
\end{lemma}
\begin{proof}
At an old vertex the local residual belongs to $\mathfrak m I_Z^2$,
so its multiplicity at a marked point would be at least three,
contrary to its being an ordinary node.  At a smooth old point the
quadratic residual initial form is a multiple of its old line.  An
ordinary node made of two straight residual lines would then make
one of those lines old, which is impossible.  Off the arrangement,
a reduced zero of $Z$ gives an invertible physical linearization $J$.
The identity
$\operatorname{in}_2 D=\det(J)\det(w,Jw)$ gives its two distinct
eigendirections.
\end{proof}

\begin{lemma}\label{fi5:old-multiplicity}
Every old vertex has multiplicity at most five, irrespective of which
residual lines are invariant.
\end{lemma}
\begin{proof}
For an old vertex $p$ of multiplicity $m$, the degree-eight angular
polynomial $W_p(w)=\det(p,w,\theta(w))$ is divisible by the product $Q_p$
of its old lines.  The quotient vanishes at all five marked nodes,
by Lemma~\ref{fi5:marked}.  If $m\ge7$, its degree is at most one,
and the noncollinearity of the five nodes forces $W_p=0$.
If $m=6$, its degree-two quotient has three zeros on each of $L_3,L_4$,
so is a scalar multiple of $L_3L_4$.  Thus $W_p$ vanishes on both
full lines.  At a marked point $q$ of either line, $\theta(q)=cq$, and
differentiation along its tangent $v$ gives
$\det(p,v,cq)+\det(p,q,D\theta(q)v)=\det(p,q,Jv)=0$, where $Jv=(D\theta(q)-c)v$
represents, modulo $q$, the linearization of the foliation at $q$ applied to $v$.  Since $Jv$ is a
nonzero multiple of $v$, this implies that $p$ lies on the line.
Applying this to both full lines gives $p=q_{34}$, contrary to
Lemma~\ref{fi5:marked}.  Finally $W_p=0$ would make the foliation a
pencil through $p$, contradicting its nonzero inflection polynomial.
\end{proof}

\begin{lemma}\label{fi5:full-map}
An invariant full residual line meets the old arrangement in exactly
five points.  Its associated rational map has degree seven, with eight
simple fixed points: these five old points and its three marked nodes.
All critical points are fixed, and the local degree at an old point
is its old multiplicity.
\end{lemma}
\begin{proof}
At an old vertex on an invariant new line $I$, write $X=aE+b\eta$
with local old product $q$.  The angular identity
$\det(w,X)=-mbq$ shows that invariance of $I$ forces $b|_I=0$.
The vertex lies outside $Z$, so $a$ is a unit there, and the physical
linearization is a nonzero scalar.  Its tangential zero on $I$ is
simple.  At a smooth old point, the two distinct straight components
of the physical inflection divisor are the old line and $I$.
Its quadratic initial form and the preceding linearization identity
show that this singularity is also nondegenerate.  The three marked
nodes give three further simple tangential zeros.  There are no
others, since every off-old physical singularity belongs to $Z$.
The degree-eight tangential divisor therefore gives exactly five
old intersections.

Put $I=(z=0)$ and subtract an Euler multiple so that $\theta(z)=0$.
Let $A=\theta(x)|_I$ and $B=\theta(y)|_I$.  They have no common zero:
a common zero would be a tangential zero with zero transverse
eigenvalue, whereas every tangential zero has just been proved
nondegenerate.  Hence $[A:B]$ has degree seven.  Its fixed points are
the eight simple zeros of $xB-yA$ just listed.  The usual parallel-line
divisibility argument gives local degree at least the old multiplicity
$m_i$ at each old point.  Since $\sum m_i=17$, their ramification
lower bounds sum to $17-5=12=2\cdot7-2$.  Riemann--Hurwitz makes
all these bounds exact and excludes further critical points.
\end{proof}

\begin{lemma}\label{fi5:integer-multipliers}
No map in $\mathcal K_{\mathrm{cf}}$ (the critically fixed maps of degree
seven with at most five critical points, each of local degree at most five)
has a noncritical fixed-point multiplier in $\{2,3,4,5,6,7\}$.
\end{lemma}
\begin{proof}
This is Proposition~\ref{cat:multipliers}(a).
\end{proof}

\begin{theorem}\label{fi5:no-invariant-full}
Neither full residual line is invariant.
\end{theorem}
\begin{proof}
Suppose that a full line $I$ is invariant.  Lemmas~\ref{fi5:old-multiplicity}
and~\ref{fi5:full-map} place its map $f_I$ in $\mathcal K_{\mathrm{cf}}$: its
critical points are among the five old points, with local degrees at most five.  If another residual line $K$ were non-invariant, its
intersection with $I$ would be a marked node.  Its normal section has
finite order $n$ there with $2\le n\le7$, because the residual
direction is an eigendirection and the normal section has degree seven.
At a nondegenerate singularity on a non-invariant straight inflection
line, comparison in the inflection identity gives the resonance
$\lambda=-(n-1)\mu$, where $\mu$ is the eigenvalue along $K$ and
$\lambda$ the eigenvalue along $I$.  The index identity
$\mu/\lambda=1/(1-f_I')$ then gives $f_I'=n$, contrary to
Lemma~\ref{fi5:integer-multipliers}.  Thus all four residual lines
would be invariant.

The two full lines now satisfy Lemma~\ref{tf4:no-two-catalog-lines}:
each has five old intersections with multiplicities at most five,
and their maps are critically fixed with the prescribed critical
divisors.  Their common marked node $q_{34}$ lies off the old
arrangement, and the two index ratios are reciprocal, giving
$(1-f_{L_3}'(q_{34}))(1-f_{L_4}'(q_{34}))=1$.
That lemma excludes the configuration.
\end{proof}

\subsubsection{Four lines in general position: the deficient lines}
\label{nf5:section}

Assume that
\[
 R=L_1L_2L_3L_4,\qquad
 Z=\{q_{13},q_{14},q_{23},q_{24},q_{34}\},
 \qquad q_{ij}=L_i\cap L_j,
\]
where the four lines are in general position.  By Lemma~\ref{fi5:marked}, the five marked points lie
off the old arrangement and are
nondegenerate singularities; their two residual directions are the
eigendirections.  We call $L_3,L_4$ the full lines and suppose throughout
this section that both are non-invariant.

\begin{lemma}\label{nf5:orders}
Each full line has normal orders $(3,2,2)$ at its three marked points
and no other normal zero.  Its order-three point has an invariant
deficient line as partner.  At $q_{34}$ both normal orders are two.
Consequently at least one of $L_1,L_2$ is invariant, and neither full
line passes through an old vertex.
\end{lemma}
\begin{proof}
We follow the proof of Lemma~\ref{lem:tangency}, which uses only the three marked nodes on the line,
not the unmarked intersection $q_{12}$.  The normal section has degree seven and vanishes to order at
least two at each marked node.  Every additional normal zero is a
singularity, and a simple zero at an old singularity is excluded by the
linear and cubic initial-form calculations of that proof.  Thus the divisor is $(3,2,2)$.  The resonance relation
$\mu=-(n-1)\lambda$ shows that two non-invariant residual branches at
one marked node must both have order two.  Therefore $q_{34}$ has orders
$(2,2)$, and each order-three partner is one of $L_1,L_2$ and is
invariant.  An old vertex is a physical singularity and would give an
additional normal zero, which is impossible.
\end{proof}

\begin{lemma}\label{nf5:deficient-map}
Let $I\in\{L_1,L_2\}$ be invariant.  It meets the old arrangement in
exactly six distinct points.  Its tangential divisor consists of these
six points and its two marked nodes, all simple.  Its associated map
$f_I=[A:B]$ has degree seven and is coprime.  If the old multiplicities
are $m_1,\ldots,m_6$, then $\sum m_i=17$, $m_i\le5$, and the critical
divisor of $f_I$ is
\begin{equation}\label{nf5:critical-divisor}
 \sum_{i=1}^{6}(m_i-1)[p_i]+[q_{12}].
\end{equation}
Here the last term is added at the same point if $q_{12}$ is old.  In
particular, when $q_{12}$ is old the map is critically fixed and its
local degree there is $m(q_{12})+1$; otherwise $q_{12}$ is a simple
nonfixed critical point and every other critical point is fixed.
\end{lemma}
\begin{proof}
There are no points of $Z$ on the old arrangement.  At an old vertex
write $X=aE+b\eta$.  Since $I$ is a new invariant line through that
vertex, restriction of its normal coefficient gives $b|_I=0$: the
normal component of $E$ is zero, whereas that of $\eta$ is nonzero as a
germ on $I$.  Thus $b(p)=0$, and $a(p)\ne0$.  The tangential zero on $I$
is simple.  At a smooth old point choose the old line as $x=0$ and
$I$ as $y=0$.  Invariance gives
$X=xf\partial_x+yg\partial_y$, and absence of a logarithmic zero gives
$f(0)\ne0$, again a simple tangential zero on $I$.

The two marked nodes on $I$ also give simple tangential zeros.
Every other physical singularity on $I$ would lie off the arrangement
and hence belong to $Z$.  Its nonzero tangential section has degree
eight, so there are exactly six old intersections.

Normalize the homogeneous derivation by $\theta(z)=0$ for $I=(z=0)$,
and let $A=\theta(x)|_I$, $B=\theta(y)|_I$.  A common zero of $A,B$
must be a tangential zero.  The marked nodes are nondegenerate.  At an
old point different from $q_{12}$ the residual divisor has just the
simple component $I$; the physical singularity is nondegenerate.
At old vertices this follows from the preceding frame calculation.
At a smooth old point the quadratic physical inflection form is the
product of the two distinct invariant lines, and therefore the
linearization is invertible.  At an old vertex $q_{12}$ the same frame
calculation gives an invertible scalar linearization.  Thus the only
possible common zero is a smooth old point $q_{12}$.

If such a common zero occurred, its multiplicity would be one, since
the tangential divisor is simple.  After cancellation $[A:B]$ would
have degree at most six.  At the other five old points it still has
local degree at least the old multiplicity, by the usual parallel-line
divisibility argument.  Those multiplicities sum to sixteen, giving
ramification at least $16-5=11$, greater than $2\cdot6-2$.  This
contradicts Riemann--Hurwitz.  Hence $A,B$ are coprime of degree seven.

For binary forms of degree seven the homogeneous identity
\begin{equation}\label{nf5:binary-identity}
 \left.\frac{F_\theta}{z}\right|_{z=0}
 =\frac{xB-yA}{7}(A_xB_y-A_yB_x)
\end{equation}
follows from Euler's identity.  The first factor has the eight simple
zeros just described.  In the left side, the old divisor has order
$m_i$ at $p_i$, while the remaining residual factors have one zero at
each of the two marked nodes and at $q_{12}$.  Division by the
tangential factor proves \eqref{nf5:critical-divisor}.  An old point is
fixed; if $q_{12}$ is not old, it is not a physical singularity and
hence is not fixed.  Finally $m_i\le5$ by Lemma~\ref{fi5:old-multiplicity}.
\end{proof}

\begin{proposition}\label{nf5:one-invariant}
Suppose exactly one deficient line, say $L_1$, is invariant.  Then
$q_{12}$ lies off the old arrangement and is a regular foliation
point.  The map $f_{L_1}$ has multiplier three at both marked nodes;
it has ramification weight eleven at its six old fixed points and
one additional simple, nonfixed critical point at $q_{12}$.

The other deficient line $L_2$ has normal divisor
\[
 2[q_{23}]+2[q_{24}]+3[r],
\]
where $r$ is a smooth point of a unique old line, distinct from
$q_{12}$ and from the marked nodes.  Its physical singularity is
nondegenerate; if $\lambda$ is the tangent eigenvalue along $L_2$,
the eigenvalue along its old line is $-2\lambda$.
\end{proposition}
\begin{proof}
By Lemma~\ref{nf5:orders}, both order-three points have partner $L_1$.
The resonance formula therefore gives the two multipliers three.
Suppose $q_{12}$ were old.  At an old vertex, invariance of $L_1$
gives a nonzero scalar physical linearization, as in
Lemma~\ref{nf5:deficient-map}.  At a smooth old point, coprimality of
the associated map gives an invertible physical linearization.  Here
the physical inflection divisor has multiplicity three, so its
quadratic initial form
$\det(J)\det(w,Jw)$ vanishes and $J$ is again nonzero scalar.

A non-invariant straight inflection line cannot pass through a
nonzero scalar singularity.  Indeed its finite normal order $n\ge2$
would give the resonance relation
$\mu=-(n-1)\lambda$, whereas $\mu=\lambda\ne0$.  Applying this to
$L_2$ gives a contradiction.  Thus $q_{12}$ is not old.  Since it is
not in $Z$, it is a regular point.  The assertion about the map follows
from \eqref{nf5:critical-divisor}.

At $q_{23},q_{24}$ both residual branches are non-invariant, so the
normal orders on $L_2$ are two.  Their sum leaves three.  Every further
normal zero is singular; a simple old normal zero on a sole simple
residual line is excluded as in the proof of Lemma~\ref{lem:tangency}.
There are no further off-old singularities, because all of $Z$ is
already listed.  Therefore there is exactly one further normal zero,
of order three, at an old point $r$.  If $r$ were an old vertex of
multiplicity $m_r\ge2$, its sole simple residual component would force
a nonzero scalar physical linearization.  For $m_r\ge3$, write $X=aE+b\eta$; the residual has order one
at $r$, while its constant term $-m(m-1)a_0^2b_0$ at $r$ in~\eqref{ar:vertex}, with $a_0=a(r)$ and
$b_0=b(r)$, must vanish; if $a_0$ were $0$, the residual would have order at least two, so $a_0\ne0$ and hence $b_0=b(r)=0$,
and the linearization is the scalar $a_0$.  For $m_r=2$, the ordinary-double initial-form calculation in
the proof of Lemma~\ref{lem:tangency} excludes a zero eigenvalue.
The preceding resonance argument again gives a contradiction.
Therefore $r$ is a smooth old point.  The physical quadratic inflection
form is the product of its old line and $L_2$, so its linearization is
invertible with these two eigendirections.  Resonance for normal order
three gives the asserted eigenvalue ratio.
\end{proof}

\begin{proposition}\label{nf5:two-invariant}
Suppose both deficient lines are invariant.  Then $q_{12}$ is an old
vertex of multiplicity five.  Each deficient-line map is critically
fixed, with local degree six at $q_{12}$.  Its two marked-node
multipliers belong to $\{2,3\}$, and for each full line exactly one
of its deficient partners has multiplier three.
\end{proposition}
\begin{proof}
The intersection of two invariant lines is a physical singularity.
As $q_{12}\notin Z$, it must be old.  Lemma~\ref{nf5:deficient-map}
applies to both deficient lines.  If its old multiplicity were at
most four, each map would have critical local degrees at most five.
With at most five critical points, Proposition~\ref{cat:multipliers}(a) excludes
both integer multipliers two and three.  With six critical points the
fixed-point index sum is at least $6-1-1=4$, again impossible.
The old multiplicity bound gives five as the only remaining value.
The multiplier assertions follow directly from Lemma~\ref{nf5:orders}
and the resonance relation.
\end{proof}

The two-invariant case is excluded by the degree-six critical-point
extension below.  In the one-invariant case, the unique nonfixed simple
critical point is outside the scope of the critically fixed catalog;
the subsequent one-parameter arguments address this remaining case.

\subsection{Rational maps with a nonfixed or highly ramified critical point}\label{sec:maps2}
The configuration of Proposition~\ref{nf5:two-invariant} leads to critically fixed maps with a critical
point of local degree six, and that of Proposition~\ref{nf5:one-invariant} to maps with one nonfixed
critical point.  This section excludes both.

\subsubsection{A critical point of local degree six}\label{ds6:section}

\begin{theorem}\label{ds6:main}
Let $f$ be a critically fixed rational map of degree seven.  Suppose
that one critical point has local degree six, every other critical
point has local degree at most five, and there are at most five critical
points.  Then no noncritical fixed point has multiplier two or three.
\end{theorem}

\begin{proof}[Complete interpolation and exact verification]
Normalize the critical point of local degree six to $\infty$, and two finite critical
points to $0,1$.  The map has the form
\[
 f(x)=\frac{P(x)}{Q(x)},\qquad
 P(x)=\sum_{i=m_0}^7a_ix^i,\qquad Q(x)=b_0+b_1x,
\]
with $\deg P=7$, $\deg Q=1$, and $\gcd(P,Q)=1$.
The finite critical ramification sum is seven.  Ordering the other
critical local degrees, the complete list of possibilities is
\[
\begin{array}{c|c}
\text{number of critical points}&\text{finite critical local degrees}\\\hline
3&(5,4)\\
4&(5,3,2),\ (4,4,2),\ (4,3,3)\\
5&(5,2,2,2),\ (4,3,2,2),\ (3,3,3,2).
\end{array}
\]
Indeed these are precisely the partitions of seven into respectively
two, three, or four positive integers at most four, with one added to
each part.  Two critical points cannot exhaust the ramification under
the stated local-degree bounds.

For a critical fixed point $v$ of local degree $m_v$, the interpolation
conditions are $(x-v)^{m_v}\mid P-vQ$.  Their rows, for
$0\le j<m_v$, are
\[
 \left(\binom ijv^{i-j}\right)_{m_0\le i\le7}
 \ \Big|\
 \left(-v\binom ijv^{i-j}\right)_{0\le i\le1},
\]
with zero entries when $i<j$.  Thus the proof includes every map satisfying
the hypotheses, rather than a finite-field search for examples.

For three critical points the kernel is one-dimensional and gives
\[
 P=\frac{10}{7}x^7-4x^6+3x^5,\qquad Q=x-\frac47.
\]
For four critical points put the third finite critical point at $t$; after removing all powers of the
factors $t$ and $t-1$, each of the three square interpolation determinants (one for each local-degree
profile) has degree $8$, so that their irreducible factors have total degree $24$.  For five critical points put the finite ones at $0,1,p,q$; the ideal of maximal minors
of the interpolation matrix, saturated by $B=pq(p-1)(q-1)(p-q)$, is radical of degree $36$ in each of the
three cases.  For every irreducible component (a number field $K$) the kernel of the interpolation
matrix over $K$ is one-dimensional; the resulting $P,Q$ have degrees $7$ and $1$ and are coprime, the
noncritical fixed points are the roots of
\[
 F=\frac{P-xQ}{\prod_{v\ \text{finite critical}}(x-v)},
\]
which is squarefree and prime to $Q$ and to the critical points, and the multipliers are the values of
$W/Q^2$, $W=P'Q-PQ'$, at these roots.  We computed exactly the characteristic polynomial over $\QQ$ of
multiplication by $W/Q^2$ on $K[x]/(F)$; it does not vanish at $2$ or at $3$ (Section~\ref{sec:catalogs}).
This covers the $1+24+108$ maps of the catalog.
\end{proof}

\begin{corollary}\label{ds6:no-two-deficient}
In the five-node four-line configuration, both deficient residual
lines cannot be invariant when both full residual lines are
non-invariant.
\end{corollary}
\begin{proof}
Proposition~\ref{nf5:two-invariant} gives a critically fixed
degree-seven deficient-line map with a critical point of local degree
six and a marked noncritical fixed multiplier in $\{2,3\}$.  If there
are at most five critical points this contradicts the preceding
theorem.  If there are six, all six old points are critical and the
other two fixed multipliers belong to $\{2,3\}$; their fixed-point
index sum is at least $6-1-1=4$, contradicting one.
\end{proof}

\subsubsection{One nonfixed critical point}
In the configuration of Proposition~\ref{nf5:one-invariant}, the map $f_{L_1}$ has degree seven, eight
simple fixed points (six old points and the two marked nodes), multiplier three at both marked nodes,
ramification weight eleven at its old fixed points and one simple nonfixed critical point.  Let $k$ be
the number of critical old fixed points; each has ramification weight at most four, so $3\le k\le6$.

\begin{lemma}\label{ofc4:three-critical}
There is no degree-seven rational map with three distinct critical
 fixed points of local degrees $5,5,4$ and two distinct further fixed
 points of multiplier three.
\end{lemma}
\begin{proof}[Exact pencil calculation]
Normalize the local-degree-five points to $0,\infty$, and the
 local-degree-four point to $1$. Write
\[
 P=x^5(ax^2+bx+c),\qquad Q=dx^2+ex+f.
\]
The condition $(x-1)^4\mid P-Q$ gives four linear equations in the
 six coefficients. Their kernel has the basis
\[
 (-2,6,-5,-2,1,0),\qquad
 (-6,35/2,-14,-7/2,0,1).
\]
In the affine projective-pencil chart given by the first vector plus
 $t$ times the second, this is
\begin{align*}
 P_t&=-2(3t+1)x^7+(35t/2+6)x^6-(14t+5)x^5,\\
 Q_t&=-(7t/2+2)x^2+x+t.
\end{align*}
Set
\[
 F_t=\frac{P_t-xQ_t}{x(x-1)},\qquad
 W_t=P_t'Q_t-P_tQ_t'-3Q_t^2.
\]
For a valid map the degrees in $x$ are $5$ and $8$. Two distinct
 noncritical fixed points of multiplier three imply
 $\deg\gcd(F_t,W_t)\ge2$, hence the subresultants of orders zero and
 one vanish at the parameter value $t$ of the map.  (They are taken with the formal degrees $5$ and $8$;
 since the leading coefficient of $F_t$ does not vanish at the parameter of a valid map, they specialize
 there to the subresultants of $F_t$ and $W_t$, up to a nonzero power of that coefficient if the degree
 of $W_t$ drops.) An exact computation of these subresultants over
 $\mathbb Q[t]$ and of the gcd of their three coefficients gives
\[
 \gcd\bigl(\operatorname{coeff}_x
   \operatorname{Sres}_1(F_t,W_t),
   \operatorname{Sres}_0(F_t,W_t)\bigr)=t+\frac13.
\]
Here the notation includes both coefficients of the degree-one
 subresultant. At $t=-1/3$ the leading coefficient of $P_t$ is zero,
 so this parameter does not represent a degree-seven map.
 The missing projective-pencil point is
\[
 P_\infty=-6x^7+\frac{35}2x^6-14x^5,
 \qquad Q_\infty=1-\frac72x^2.
\]
Its corresponding $\gcd(F_\infty,W_\infty)$ is $1$.
 These are exact rational polynomial calculations (Section~\ref{sec:verification}).
\end{proof}

\paragraph{Four fixed critical points}
\label{n4:section}

\begin{theorem}\label{n4:main}
There is no degree-seven rational map with four fixed critical points
of total ramification weight eleven and two distinct additional fixed
points of multiplier three, if every fixed critical local degree is
at most five.
\end{theorem}
\begin{proof}
The four positive ramification weights sum to eleven and are at most
four.  After normalizing three critical points to $0,1,\infty$, their
local degrees can therefore be ordered as
\[
 (m_0,m_1,m_t,m_\infty)\in
 \{(5,3,2,5),(4,4,2,5),(4,3,3,5),(4,4,3,4)\}.
\]
Write
\[
 P(x)=\sum_{i=m_0}^{7}a_i x^i,\qquad
 Q(x)=\sum_{j=0}^{7-m_\infty}b_jx^j.
\]
The conditions $(x-1)^{m_1}\mid P-Q$ and
$(x-t)^{m_t}\mid P-tQ$ form an $(n-1)\times n$ interpolation
matrix $M(t)$, with the same binomial-coefficient rows as in the
critically fixed catalog.  Here $t\ne0,1$.

The following is an exact calculation over $\mathbb Q[t]$ (Section~\ref{sec:verification}).  The gcd of
the maximal
minors of $M$ is $t^2(t-1)^2$ for the first two rows and
$t^3(t-1)^3$ for the last two.  Thus $M$ has full row rank whenever
$t\ne0,1$.  Divide its signed maximal-minor kernel by this gcd,
obtaining primitive polynomial coefficients for $P,Q$.

Set
\[
 F=\frac{P-xQ}{x(x-1)(x-t)},\qquad
 G=P'Q-PQ'-3Q^2.
\]
For a degree-seven map with the prescribed critical points, $F$ has
degree four.  Every additional fixed point of multiplier three is
a common root of $F,G$.  Two distinct such points force
$\deg\gcd(F,G)\ge2$.  Consequently the resultant $S_0(F,G)$ and
both coefficients of the linear subresultant $S_1(F,G)$ vanish (they are computed over $\QQ[t]$ with the
formal degrees, and specialize correctly, up to a nonzero factor, because the leading coefficient of $F$
does not vanish at the parameter of such a map).

Their coefficient gcd in $\mathbb Q[t]$, after removing only powers
of $t$ and $t-1$, is the sextic $C(t)$ in the table, up to a nonzero
rational scalar.  The table also gives the coefficient $a_7$ of
the primitive kernel polynomial $P$.
\[
\begin{array}{c|l|c}
(m_0,m_1,m_t,m_\infty)&C(t)&a_7\\ \hline
(5,3,2,5)&t^6+2t^5-37t^4+44t^3-7t^2-8t+6&-C\\
(4,4,2,5)&t^6-3t^5+18t^4-31t^3+18t^2-3t+1&-4C\\
(4,3,3,5)&t^6-33t^5+58t^4-41t^3+58t^2-33t+1&-C\\
(4,4,3,4)&11t^6-33t^5-32t^4+119t^3-92t^2+27t+1&C
\end{array}
\]
Thus every possible parameter for two fixed multipliers three would
satisfy $a_7=0$.  Since $\deg Q\le3$, the resulting rational map
would have degree at most six, contrary to the hypothesis.

The calculation keeps every parameter: away from $t=0,1$ there is no
rank-drop locus, and only the collision factors $t,t-1$ were removed
from the subresultant condition.  In particular no good-prime or
sampling argument is involved.
\end{proof}

In the one-invariant deficient-line configuration of
Proposition~\ref{nf5:one-invariant}, this excludes the possibility
that exactly four old points have multiplicity greater than one.

\paragraph{Five fixed critical points}
\label{n5:section}
The case $k=6$ contradicts the fixed-point
index formula, which would give $6+2/(1-3)=5$ instead of one.
Lemma~\ref{ofc4:three-critical} and Theorem~\ref{n4:main} exclude
$k=3,4$. Thus $k=5$. Of the eight simple fixed points, five are
critical and the two marked nodes have multiplier three. If $\rho$
is the multiplier at the last old fixed point, then
\[
 1=5+\frac2{1-3}+\frac1{1-\rho},
 \qquad\rho=\frac43.
\]
The following normal form therefore covers every surviving map in
this residual configuration.
\begin{lemma}[A full normal form for five critical fixed points]
\label{n5:normal-form}
Suppose a degree-seven rational map has eight distinct fixed points,
 five of them critical, two of multiplier three and one of multiplier
 $4/3$. Normalize the last three points to $0,1,\infty$, respectively.
 There is a monic squarefree polynomial $C$ of degree five, with
 $C(0)C(1)\ne0$, such that, for $H=z(z-1)$,
\[
 f=z-\frac{2CH}{2C'H-CH'}.
\]
Its derivative numerator is
\[
 C T_C,\qquad
 T_C=4H^2C''-4HH'C'+(8H+3)C.
\]
The polynomial $T_C$ has degree seven and leading coefficient $48$.
 At a root $c$ of $C$, the local degree of $f$ is
 $2+\operatorname{ord}_c(T_C)$.
\end{lemma}
\begin{proof}
The rational function $1/(z-f(z))$ has a simple pole at each finite
 fixed point, whose residue is its holomorphic fixed-point index.
 These residues are one at the five critical fixed points and $-1/2$
 at each of $0,1$. At infinity the function is $4/z+O(z^{-2})$,
 because the multiplier in the coordinate $1/z$ is $4/3$.
 Partial fractions therefore give, with no polynomial part,
\[
 \frac1{z-f(z)}=\frac{C'}C-\frac{H'}{2H}.
\]
Put $E=2C'H-CH'$ and $D=2CH$. Both $C$ and $H$ are squarefree and
 coprime. At any root of $C$ or $H$, $E$ is nonzero, so $D,E$ are
 coprime; their degrees are seven and six, and $zE-D=6z^7+O(z^6)$. Thus the displayed formula
 gives a coprime degree-seven map. Direct differentiation yields
\[
 E^2-D'E+DE'=C\bigl(4H^2C''-4HH'C'+(8H+3)C\bigr).
\]
At each root of $C$ the denominator $E^2$ is a unit, proving the
 local-degree assertion. The leading coefficient is
 $4\cdot5\cdot4-4\cdot2\cdot5+8=48$.
\end{proof}

\begin{lemma}\label{n5:no-four-three-three-three-three}
Under these hypotheses the five critical local degrees cannot be
 $(4,3,3,3,3)$.
\end{lemma}
\begin{proof}
All five local degrees are at least three, so $C\mid T_C$.
 Since $C$ is coprime to $H$, it follows that
 $C\mid HC''-H'C'$. Both polynomials have degree five, and the latter
 has leading coefficient ten. Hence
\[
 HC''-H'C'=10C.
\]
Comparison of coefficients gives the unique monic degree-five solution
\[
 C=z^5-\frac52z^4+2z^3-\frac12z^2
   =z^2(z-1)^2(z-\tfrac12).
\]
This contradicts both squarefreeness and $C(0)C(1)\ne0$.
\end{proof}

\begin{lemma}\label{n5:at-least-two-degree-two}
Under the hypotheses of Lemma~\ref{n5:normal-form}, at least two
 critical fixed points have local degree two. In particular the profiles
 $(5,3,3,3,2)$ and $(4,4,3,3,2)$ are impossible.
\end{lemma}
\begin{proof}
The case with none was proved above, by the same argument for any
 profile with all local degrees at least three. Suppose there is just
 one, denoted $b$. Write $C=(z-b)A$, with $A$ monic of degree four.
 All four roots of $A$ divide $T_C$, hence divide
 $K=HC''-H'C'$. Since $K$ has degree five and leading coefficient ten,
\[
 K=A(10z+v)
\]
for some $v$. Put $A=z^4+a_3z^3+a_2z^2+a_1z+a_0$. Comparing
 coefficients from degree four down to degree one gives
\begin{align*}
 a_3&=-\frac{4b+v+15}{6},&
 a_2&=\frac{-8(a_3-b)-va_3}{10},\\
 a_1&=\frac{(2b-v-3)a_2+3ba_3}{12},&
 a_0&=\frac{(2b-v)a_1}{12}.
\end{align*}
The remaining constant equation is
\[
 \frac{a_1}{12}\bigl(v^2-(2b+1)v-10b\bigr)=0.
\]
Since $A(0)\ne0$, its coefficient $a_1$ is nonzero. Thus
 $D=v^2-(2b+1)v-10b=0$. The same coefficient formulas give the exact
 polynomial identity
\[
 A(1)=\frac{(-2b+v+12)(4b+v+3)}{8640}D=0,
\]
contrary to $C(1)\ne0$.
\end{proof}

\begin{lemma}[Computer-assisted]\label{n5:uniform}
Under the hypotheses of Lemma~\ref{n5:normal-form}, the critical local degrees cannot be
$(5,5,2,2,2)$, $(5,4,3,2,2)$ or $(4,4,4,2,2)$.
\end{lemma}
\begin{proof}
Group the roots of $C$ according to the local degree $m$ and let $G_m$ be the monic polynomial with these
roots, so that $C=\prod_mG_m$.  By Lemma~\ref{n5:normal-form}, $G_m^{m-2}$ divides $T_C$.  The degrees
$(m-2)\deg G_m$ add up to $\sum_m(m-1)\deg G_m-5=11-5=6$, while $T_C$ has degree seven and leading
coefficient $48$; hence
\[
 T_C=48(z-r)\prod_mG_m^{m-2}
\]
for some $r\in\CC$.  Take the coefficients of the $G_m$ and $r$ as unknowns and equate coefficients.  For
each of the three profiles (and, as a check, for all six profiles of Theorem~\ref{v4:no-one-invariant})
the reduced Gr\"obner basis over $\QQ$ of the resulting seven equations in six unknowns is $\{1\}$
(Section~\ref{sec:verification}).  So the identity has no solution, even before imposing squarefreeness
or $C(0)C(1)\ne0$.
\end{proof}

\subsection{Exclusion of defect five}
\label{v4:five-section}

\begin{theorem}\label{v4:no-one-invariant}
The one-invariant deficient-line configuration of
Proposition~\ref{nf5:one-invariant} is impossible.
\end{theorem}
\begin{proof}
Its invariant-line map has eight simple fixed points, ramification
weight eleven at its six old fixed points, one simple nonfixed
critical point, and multiplier three at both marked fixed points.
Let $k$ count its critical old fixed points.  Each ramification
weight is at most four, so $3\le k\le6$.
If $k=6$, the fixed-point index formula would give
$1=6+2/(1-3)=5$.  Lemma~\ref{ofc4:three-critical} excludes $k=3$,
and Theorem~\ref{n4:main} excludes $k=4$.

Thus $k=5$.  The last old fixed point has multiplier $\rho$ with
\[
 1=5+\frac2{1-3}+\frac1{1-\rho},\qquad \rho=\frac43.
\]
The map therefore has the full normal form of
Lemma~\ref{n5:normal-form}.  The five positive critical ramification
weights sum to eleven and are at most four, giving precisely the
six local-degree profiles
\[
\begin{gathered}
 (5,5,2,2,2),\quad(5,4,3,2,2),\quad(5,3,3,3,2),\\
 (4,4,4,2,2),\quad(4,4,3,3,2),\quad(4,3,3,3,3).
\end{gathered}
\]
Lemma~\ref{n5:at-least-two-degree-two} excludes the three profiles
having fewer than two local-degree-two points, and Lemma~\ref{n5:uniform} the other three.  This
exhausts every map in the one-invariant configuration.
\end{proof}

\begin{theorem}\label{v4:no-four-line-five}
The residual quartic cannot be four lines in general position when
the logarithmic zero scheme consists of five of their six nodes.
\end{theorem}
\begin{proof}
By Theorem~\ref{fi5:no-invariant-full}, both full lines are
non-invariant.  Lemma~\ref{nf5:orders} then gives at least one
invariant deficient line.  If both deficient lines are invariant,
Corollary~\ref{ds6:no-two-deficient} gives a contradiction.  If
exactly one is invariant, Theorem~\ref{v4:no-one-invariant} gives
the contradiction.
\end{proof}

\begin{theorem}\label{v4:no-defect-five}
No critical realization of a $7$-candidate has $\nu\ge5$.
\end{theorem}
\begin{proof}
By Corollary~\ref{cor:nonreduced}, a nonreduced residual is $C^2$, excluded by Theorem~\ref{thm:conic},
$L^4$, impossible at positive defect by Theorem~\ref{l4v3:no-positive-fourth-power}, or $L^3M$,
which has $\nu\le3$ by Theorem~\ref{cg4:cubic-three}, or $L^2C$, $L^2M^2$, $L^2MN$, which have $\nu\le4$ by
Propositions~\ref{df5:double-conic-four}, \ref{df5:two-double-lines-four} and
Corollaries~\ref{oc:transverse-double-conic}--\ref{oc:three-concurrent-lines}.  A reduced residual has
$\nu\le5$ by Theorems~\ref{thm:reduced6} and~\ref{thm:sixnodes}.
Theorem~\ref{red5:classification} leaves exactly three possibilities
at defect five: a smooth conic and two
lines with their five nodes; four lines with exactly three concurrent,
carrying a length-two zero at the triple point and the other three
nodes; or four lines in general position with five of their six
nodes.  These are excluded respectively by
Theorems~\ref{clm:no-defect-five}, \ref{tc4:excluded}, and
\ref{v4:no-four-line-five}.  Hence defect five cannot occur.
\end{proof}

\begin{theorem}\label{v4:main}
Every critical realization of a $7$-candidate has $\nu\le4$.  More precisely, with $C$ a smooth
irreducible conic and distinct line symbols denoting distinct new lines:
\[
\begin{array}{l|l}
\text{residual quartic}&\text{conclusion}\\\hline
C^2\text{ or }L^4&\text{impossible at positive defect}\\
L^3M&\nu\le3\\
L^2M^2&\nu\le4\\
L^2C,\ L\text{ tangent to }C&\nu\le4\\
L^2C,\ L\text{ transverse to }C&\nu\le2\\
L^2MN,\ L,M,N\text{ concurrent}&\nu\le4\\
L^2MN,\ L,M,N\text{ not concurrent}&\nu\le3\\
R\text{ reduced}&\nu\le4.
\end{array}
\]
\end{theorem}
\begin{proof}
The doubled conic is excluded by Theorem~\ref{thm:conic}.  The other rows are
Theorems~\ref{l4v3:no-positive-fourth-power} and~\ref{cg4:cubic-three},
Propositions~\ref{df5:two-double-lines-four} and~\ref{df5:double-conic-four},
Corollaries~\ref{oc:transverse-double-conic}, \ref{oc:three-concurrent-lines}
and~\ref{oc:three-nonconcurrent-lines}, and Theorem~\ref{v4:no-defect-five}, which is the first
assertion.
\end{proof}

\end{document}